\documentclass[reqno]{amsart}
\usepackage{amsmath,amssymb}

\newcommand{\ep}{\ensuremath{\epsilon}}

\newcommand{\CC}{\ensuremath{\mathbb{C}}}

\newcommand{\NN}{\ensuremath{\mathbb{N}}}

\newcommand{\RR}{\ensuremath{\mathbb{R}}}
\newcommand{\TT}{\ensuremath{\mathbb{T}}}
\newcommand{\ZZ}{\ensuremath{\mathbb{Z}}}

\numberwithin{equation}{section}
\theoremstyle{plain}
\newtheorem{thm}{Theorem}[section]
\newtheorem{cor}[thm]{Corollary}
\newtheorem{lem}[thm]{Lemma}
\newtheorem{prop}[thm]{Proposition}
\newtheorem{Def}[thm]{Definition}
\theoremstyle{remark}
\newtheorem{rem}[thm]{Remark}

\begin{document}
\title{The big $q$-Jacobi function transform}
\author{Erik Koelink}
\author{Jasper V. Stokman}
\date{April 21, 1999}
\address{Technische Universiteit Delft, Werkeenheid Algemene Wiskunde,
Faculteit Informatietechnologie en
Systemen, Postbus 5031, 2600 GA Delft, The Netherlands}
\email{koelink@twi.tudelft.nl}
\address{Universit{\'e} Louis Pasteur, Institut de Recherche Math{\'e}matique
Avanc{\'e}e, 7 rue Ren{\'e} Descartes, 67084 Strasbourg, France}
\email{stokman@math.u-strasbg.fr}
%%%%%%%%%%%%%%%%%%%%%%%%%%%%%%%%%%%%%%%%%%%%%%%%%%%%%%%%%%%%%%%%%
%%                                                             %%
%%                        Abstract                             %%
%%                                                             %%
%%%%%%%%%%%%%%%%%%%%%%%%%%%%%%%%%%%%%%%%%%%%%%%%%%%%%%%%%%%%%%%%%
\begin{abstract}
We give a detailed description of the resolution of
the identity of a second order $q$-difference operator 
considered as an unbounded self-adjoint operator on 
two different Hilbert spaces. 
The $q$-difference operator and the 
two choices of Hilbert spaces naturally arise from harmonic analysis
on the quantum group $SU_q(1,1)$ and $SU_q(2)$.
The spectral analysis
associated to $SU_q(1,1)$ leads to the big $q$-Jacobi function transform
together with its Plancherel measure and inversion formula.
The dual orthogonality relations give a one-parameter family
of non-extremal orthogonality measures for 
the continuous dual $q^{-1}$-Hahn polynomials
with $q^{-1}>1$, and explicit sets of functions which complement these
polynomials to orthogonal bases of the associated Hilbert spaces.
The spectral analysis associated to $SU_q(2)$ leads to a functional
analytic proof of the orthogonality relations and quadratic norm
evaluations for the big $q$-Jacobi polynomials. 
\end{abstract}
%%%%%%%%%%%%%%%%%%%%%%%%%%%%%%%%%%%%%%%%%%%%%%%%%%%%%%%%%%%%%

\maketitle
%%%%%%%%%%%%%%%%%%%%%%%%%%%%%%%%%%%%%%%%%%%%%%%%%%%%%%%%%%%%%%%%
%%                                                            %%
%%                    Introduction                            %%
%%                                                            %%
%%%%%%%%%%%%%%%%%%%%%%%%%%%%%%%%%%%%%%%%%%%%%%%%%%%%%%%%%%%%%%%%
\section{Introduction}\label{section1}

Harmonic analysis on the compact quantum group $SU_q(2)$ has been
studied by several authors, see for example \cite{VS}, \cite{NM}, \cite{K2},
\cite{Koel}. 
One possible approach is studying the restriction of the Haar functional and of 
the action of the quantum Casimir 
to ``functions'' in $SU_q(2)$ of a given, fixed bi-$T$-type, 
where $T$ is the 
standard compact torus in $SU(2)$. This reduces the harmonic analysis 
to the spectral analysis of 
an explicit second order $q$-difference operator considered as 
unbounded linear operator on an explicit Hilbert space. 
The Hilbert space is the $L^2$-space corresponding to the orthogonality
measure of the little $q$-Jacobi polynomials, while
the second order
$q$-difference operator is diagonalized by the little $q$-Jacobi polynomials,
so this leads to an interpretation of (the orthogonality 
relations of) the little $q$-Jacobi polynomials on $SU_q(2)$.

This approach was employed by  
Kakehi, Masuda and Ueno \cite{KMU} and  Kakehi \cite{K} for the development of 
harmonic analysis on the non-compact quantum
group $SU_q(1,1)$, see also Vaksman and Korogodsky \cite{KV}. 
The restriction of the Casimir to functions on
$SU_q(1,1)$ of a given, fixed  
bi-$T$-type leads to the same second order $q$-difference operator as for
$SU_q(2)$, while the restriction of a Haar functional on $SU_q(1,1)$
to the given bi-$T$-type
leads to a $L^2$-space with respect to a discrete measure with 
unbounded support of the form $[0,\infty(z))_q=\{zq^k\}_{k\in\ZZ}$ ($z\not=0$),
where $0<q<1$. The spectral analysis of 
the second order $q$-difference operator considered as unbounded
linear operator on this Hilbert space was developed in \cite{KMU}, \cite{K}. 
It leads to the Plancherel formula and the inversion formula for
the little $q$-Jacobi function transform.

It is well known that harmonic analysis on $SU_q(2)$ can be generalized by
replacing the role of the torus $T$ in $SU_q(2)$ by a ``conjugate''
$T_t$ ($t\in \RR$), which is defined in terms
of twisted primitive elements in the quantized universal enveloping
algebra $U_q\bigl(su(2)\bigr)$, see \cite{K2}. 
In particular, considering the harmonic analysis on $SU_q(2)$
with respect to left $T$-types and right $T_t$-types leads to
an interpretation of 
the big $q$-Jacobi polynomials on $SU_q(2)$, see \cite{NM}. 
Furthermore, the restriction of the Haar functional 
to functions of a given left $T$-type and right $T_t$-type can be
computed directly and identified with the orthogonality measure
of the big $q$-Jacobi polynomials, see \cite{KV} for the spherical case.

In this paper we develop the functional analytic aspects 
related to the harmonic analyis on $SU_q(1,1)$ 
with respect to left $T$-types and right
$T_t$-types. Concretely, we give a detailed description of 
the spectral properties of the second order $q$-difference operator $L$
which is diagonalized by the big $q$-Jacobi polynomials, considered as
an unbounded linear operator on a one-parameter family of Hilbert spaces.
The Hilbert spaces are $L^2$-spaces corresponding to explicit 
discrete measures with unbounded support $[-1,\infty(z))_q=\{-q^k\}_{k\in\ZZ_+}\cup
\{zq^k\}_{k\in\ZZ}$ ($z>0$). They can be interpreted 
for specific parameter values as the 
restriction of a Haar functional on $SU_q(1,1)$ to
functions of a given left $T$-type and right $T_t$-type.

We construct a one-parameter family of dense domains for
which $L$ is self-adjoint. 
With respect to one of these domains, we explicitly compute the 
resolution of the identity of $L$. 
This leads to the big $q$-Jacobi function transform, its Plancherel
formula and its inversion formula. Furthermore, we show that the corresponding
dual orthogonality relations lead to a one-parameter family of non-extremal
orthogonality measures for the continuous dual $q^{-1}$-Hahn polynomials, 
and explicit sets of functions complementing these 
polynomials to orthogonal bases of the corresponding Hilbert spaces.

The little (respectively big) $q$-Jacobi function transform 
is a $q$-analogue of the Jacobi function transform
on the interval $[0,\infty)$ (respectively $[-1,\infty)$), 
see e.g. Braaksma and Meulenbeld \cite{BM}, Koornwinder
\cite{K1} and references given there. In the classical setting, these
two transforms are related by a dilation of the geometric
parameter in which the end-point $0$ of the interval $[0,\infty)$ is
mapped onto the end-point $-1$ of the interval $[-1,\infty)$.
The key point in our study of the big $q$-Jacobi function transform is
to show that the limiting
point $0$ of the $q$-interval $[-1,\infty(z))_q$ 
does not play a special role in the spectral analysis, while the 
role of the end-point $-1$ is similar to the role of the 
limiting point $0$ for the little $q$-Jacobi function transform, see
\cite{KMU}, \cite{K}.

We will report elsewhere in detail about the 
connection between the big $q$-Jacobi function transform and harmonic
analysis on $SU_q(1,1)$. Furthermore, we will report elsewhere on
the Plancherel formula and the inversion formula 
for the Askey-Wilson function transform, which  
is associated to harmonic analysis on $SU_q(1,1)$ with
respect to left $T_s$-types and right $T_t$-types ($s,t\in\RR$).

The organization of the paper is as follows.

In section 2 we introduce the second order $q$-difference operator $L$ 
and the Hilbert spaces 
on which we consider $L$ as an unbounded linear operator.
Furthermore, we derive dense domains 
for which $L$ is self-adjoint. The domains are given explicitly in terms
of continuity and differentiability conditions of the functions at the origin.

In section 3 we derive criteria when it is possible to extend
an arbitrary eigenfunction of $L$ on the positive real axis 
to a global eigenfunction in such a way that the solution 
is ``continuously differentiable'' at the origin. We furthermore give,
in terms of basic hypergeometric series, 
two explicit eigenfunctions which are continuously differentiable at the origin.
One of them is
the spherical function for the big $q$-Jacobi function transform.

In section 4 we introduce the asymptotically free eigenfunction of $L$ 
on the positive real axis and we give the 
corresponding $c$-function expansion of
the spherical function on the positive real axis.

In section 5 we extend for generic parameters the asymptotically free solution
to a global eigenfunction of $L$ in such a way, that the global eigenfunction
is continuously differentiable at the origin. This leads to the $c$-function
expansion of the spherical function on the whole support of the measure. 

In section 6 we define the Green function of $L$ in terms
of the spherical function and the extended asymptotically free eigenfunction.

In section 7 we use the Green function to compute the continuous 
contribution to the resolution of the identity of $L$. Furthermore,
we derive the Plancherel formula and the inversion
formula for the continuous part of the big $q$-Jacobi function transform.

In section 8 we derive the discrete contribution of the resolution of the 
identity of $L$ and we derive the orthogonality relations and the quadratic
norm evaluations of the corresponding spherical functions, which are square
integrable in these cases.
We furthermore give a precise description of the spectrum
of $L$.

In section 9 we state and prove the main results of the paper. We
derive the Plancherel formula and the inversion formula for the big $q$-Jacobi 
function transform. Furthermore, we show that the dual orthogonality relations
lead to a one-parameter family of non-extremal orthogonality measures
for the continuous dual $q^{-1}$-Hahn polynomials, as well as explicit sets
of functions which complement the polynomials to orthogonal bases
of the corresponding Hilbert spaces.

Finally we derive in section 10 a functional analytic proof of 
the orthogonality relations and the quadratic norm evaluations 
for the big $q$-Jacobi polynomials. Proofs will only be sketched 
in this section,
we mainly emphasize the differences
between the compact setting, which 
corresponds to $SU_q(2)$ and the big $q$-Jacobi
polynomials, and the non-compact setting, which corresponds to 
$SU_q(1,1)$ and the big $q$-Jacobi functions.

{\it Notations:} We assume throughout the paper that $0<q<1$ is fixed.
We follow the notation of Gasper and Rahman
\cite{GR} concerning $q$-shifted factorials and basic
hypergeometric series. We write
\[\theta(x_1,\ldots, x_r)=\theta(x_1)\ldots\theta(x_r)
\]
for products of (renormalized) Jacobi theta-products 
$\theta(x)=\bigl(x,q/x;q\bigr)_{\infty}$. We write $\ZZ_+=\{0,1,2,\ldots,\}$
and $\NN=\{1,2,\ldots\}$.

{\it Acknowledgements:} The research for this paper started when both 
authors were affiliated to the University of Amsterdam, the first author
as post-doc supported by the Netherlands Organization for 
Scientific Research (NWO)
under project number 610.06.100. Part of the research was done
while the second author was a post-doc at the Centre de Math{\'e}matiques
de Jussieu, Universit{\'e} Paris VI Pierre et Marie Curie, Paris, France, 
supported by the EC TMR network ``Algebraic Lie Representations'', Grant no.
ERB FMRX-CT97-0100. The second author is also supported by
a NWO-TALENT stipendium of
the Netherlands Organization for Scientific Research (NWO).

We thank Tom H. Koornwinder for sharing with us 
his private notes on the little
$q$-Jacobi function transform and for his interest and valuable comments
on the subject.

%%%%%%%%%%%%%%%%%%%%%%%%%%%%%%%%%%%%%%%%%%%%%%%%%%%%%%%%%%%%%%%%%%%
%%                                                               %%
%%            The second order $q$-difference operator           %%
%%                                                               %%
%%%%%%%%%%%%%%%%%%%%%%%%%%%%%%%%%%%%%%%%%%%%%%%%%%%%%%%%%%%%%%%%%%%

\section{The second order $q$-difference operator}\label{sdomain}
In this section we consider the second order $q$-difference operator
$L$ associated with the big $q$-Jacobi polynomials as an unbounded operator
on a one-parameter family of Hilbert space.
We determine suitable domains of definition for which $L$ is self-adjoint.

The second order $q$-difference operator $L$ depends on three parameters
$(a,b,c)$. The corresponding one-parameter family of
Hilbert spaces thus depends on four parameters $(a,b,c,z)$, where
the parameter $z$ labels the Hilbert spaces.
In sections \ref{sdomain}--\ref{PID} we assume that the four parameters
$(a,b,c,z)$ satisfy $z>0$ and $(a,b,c)\in V$, where
\begin{equation}\label{con1}
V=\{(a,b,c) \,\, | \,\, a,b,c>0,\quad ab,ac,bc<1\},
\end{equation}
unless explicitly stated otherwise.
Observe that two of the three parameters
$a,b,c$ take their values in the open interval $(0,1)$.

We define the $q$-interval $I$ by
\[ I=[-1,\infty(z))_q=\{ -q^k \, | \, k\in\ZZ_+\}\cup
\{zq^k \, | \, k\in\ZZ\},
\]
which is regarded as a discrete $q$-analogue
of the interval $[-1,\infty)$. We use the notation $(-1,\infty(z))_q$
for the $q$-interval $I\setminus \{-1\}=[-q,\infty(z))_q$.
Let ${\mathcal{F}}(I)$ be the
linear space of complex-valued functions $f: I\to\CC$, and define
a linear operator $L\in {\hbox{End}}_{\CC}\bigl({\mathcal{F}}(I)\bigr)$
by
\begin{equation}\label{equiv2}
L=A(\cdot)\bigl(T_q-\hbox{Id}\bigr)+B(\cdot)\bigl(T_{q^{-1}}-\hbox{Id}\bigr),
\end{equation}
with the $q$-shift operators defined by 
$\bigl(T_{q^{\pm 1}}f\bigr)(x)=f(q^{\pm 1}x)$, and with
\begin{equation}\label{ABbig}
A(x)=a^2\left(1+\frac{1}{abx}\right)\left(1+\frac{1}{acx}\right),\quad
B(x)=\left(1+\frac{q}{bcx}\right)\left(1+\frac{1}{x}\right).
\end{equation}
Here $\bigl(Lf\bigr)(-1)$ is by definition given by 
\begin{equation}\label{special}
\bigl(Lf\bigr)(-1)=A(-1)\bigl(f(-q)-f(-1)\bigr),\quad f\in
{\mathcal{F}}(I),
\end{equation}
which is formally compatible with the fact that $B(-1)=0$. 

In the next lemma
we rewrite the operator $L$ as a second order operator in the
$q$-difference operator $D_q$, which is defined by
\[\bigl(D_qf\bigr)(x)=\frac{f(x)-f(qx)}{(1-q)x}.
\]
%%%%%%%%%%%%%%%%%%%%%%%%%%%%%%%%%%%%%%%%%%%%%%%%%%%%%%%%%%%%%%%%%%
\begin{lem}\label{self-adjointform}
Let $f\in {\mathcal{F}}(I)$, then
\begin{equation}\label{equiv1}
\bigl(Lf\bigr)(x)=
\begin{cases}
{\displaystyle{p(x)\bigl(D_q(r(\cdot)D_qf)\bigr)(q^{-1}x)}},\qquad &x\in 
(-1,\infty(z))_q,\\
{\displaystyle{\frac{-qp(x)r(x)}{(1-q)x}}}\bigl(D_qf\bigr)(x),
\qquad &x=-1,
\end{cases}
\end{equation}
where
\begin{equation}\label{pr}
\begin{split}
p(x)&= \frac{\bigl(-abx,-acx;q\bigr)_{\infty}}
{\bigl(-bcx, -qx;q\bigr)_{\infty}},\\
r(x)&= \frac{(1-q)^2}{qbc}
\frac{\bigl(-bcx, -qx;q\bigr)_{\infty}}
{\bigl(-qabx,-qacx;q\bigr)_{\infty}}.
\end{split}
\end{equation}
\end{lem}
%%%%%%%%%%%%%%%%%%%%%%%%%%%%%%%%%%%%%%%%%%%%%%%%%%%%%%%%%%%%%%%%%%%%%
\begin{proof}
Observe that $p$ and $r$ are well defined as functions on the
$q$-interval $I$ since $(a,b,c)\in V$. 
The formula for $x\in (-1,\infty(z))_q$ follows by observing that 
$p(\cdot)T_{q^{-1}}D_qr(\cdot)D_q$ is of the form \eqref{equiv2} with
$A$ and $B$ given by
\[
A(x)=\frac{qp(x)r(x)}{(1-q)^2x^2}, \quad
B(x)=\frac{q^2p(x)r(q^{-1}x)}{(1-q)^2x^2}.
\]
By inserting the explicit functions $p$ and $r$, we see that the
corresponding $A$ and $B$
coincide with the ones given by \eqref{ABbig}.
The formula for $Lf$ in the point $x=-1$ follows now 
immediately from the definition \eqref{special} of
$\bigl(Lf\bigr)(-1)$.
\end{proof}
%%%%%%%%%%%%%%%%%%%%%%%%%%%%%%%%%%%%%%%%%%%%%%%%%%%%%%%%%%%%%%%%%

We have that $p(x)>0$ for all $x\in I$ since $(a,b,c)\in V$. 
We can thus define the Hilbert space
\begin{equation}
{\mathcal{H}}=\{\,\, 
f\in {\mathcal{F}}(I) \,\,\,\,\, | \,\,\,\,\, \|f\|^2=\langle
f,f\rangle <\infty \,\,\},
\end{equation}
where
\begin{equation}
\langle f,g\rangle = \int_{-1}^{\infty(z)}f(x){\overline{g(x)}}
\frac{d_qx}{p(x)},
\end{equation}
with the Jackson $q$-integral for $\alpha,\gamma\in \CC^*$ defined by
\begin{equation}\label{Jackson}
\begin{split}
\int_{\alpha}^{\beta}f(x)d_qx&=\int_{0}^{\beta}f(x)d_qx-
\int_{0}^{\alpha}f(x)d_qx,\quad \beta=\gamma,\infty(\gamma)\\
\int_{0}^{\gamma}f(x)d_qx&=(1-q)\sum_{n=0}^{\infty}f(\gamma q^n)\gamma
q^n,\\
\int_{0}^{\infty(\gamma)}f(x)d_qx&=
(1-q)\sum_{n=-\infty}^{\infty}f(\gamma q^n)\gamma q^n,
\end{split}
\end{equation}
for functions $f$ such that the sums converge absolutely.
If $\alpha=q^k\gamma$ for certain $k\in \ZZ_+$, then the
$q$-integral from $\alpha$ to $\gamma$ can be defined for arbitrary
functions $f$, since the $q$-integral reduces then to a finite sum.

We regard $L$ as an unbounded operator on the Hilbert space
${\mathcal{H}}$. In the remainder of this section, we define
a one-parameter family of dense subspaces ${\mathcal{D}}\subset{\mathcal{H}}$ 
such that $L$, with domain of definition ${\mathcal{D}}$, is self-adjoint.
In order to determine the subspaces, it is convenient to 
consider the symmetry of $L$ with respect to a truncated 
version of the inner product $\langle
.,.\rangle$, which is defined for $k\in \ZZ_+$ and $l,m\in \ZZ$ with $l<m$ by
\begin{equation*}
\begin{split}
\langle f,g \rangle_{k;l,m}&=\left(\int_{-1}^{-q^{k+1}}+
\int_{zq^{m+1}}^{zq^l}\right)f(x){\overline{g(x)}}\frac{d_qx}{p(x)}\\
&=
\sum_{n=0}^kf(-q^n){\overline{g(-q^n)}}\frac{(1-q)q^n}{p(-q^n)}+
\sum_{n=l}^mf(zq^n){\overline{g(zq^n)}}\frac{(1-q)zq^n}{p(zq^n)}.
\end{split}
\end{equation*}
Observe that for all $f,g\in {\mathcal{H}}$, we have
\begin{equation}\label{lim}
\lim_{\stackrel{{\scriptstyle{k,m\to\infty}}}{l\to -\infty}}
\langle f,g\rangle_{k;l,m}=\langle f,g\rangle.
\end{equation}
We have now the following lemma. 
%%%%%%%%%%%%%%%%%%%%%%%%%%%%%%%%%%%%%%%%%%%%%%%%%%%%%%%%%%%%%%%
\begin{lem}\label{stoktermen}
Let $f,g\in {\mathcal{F}}(I)$,
$k\in \ZZ_+$ and $l,m\in \ZZ$ with $l<m$, then
\begin{equation*}
\langle Lf,g\rangle_{k;l,m}-\langle f,Lg\rangle_{k;l,m}=
W(f,{\overline{g}})(zq^{l-1})-W(f,{\overline{g}})(zq^m)
+W(f,{\overline{g}})(-q^k),
\end{equation*}
with ${\overline{f}}(x)={\overline{f(x)}}$ and with
the Wronskian $W(f,g)\in {\mathcal{F}}(I)$ defined by
\begin{equation}\label{Wronskian}
\begin{split}
W(f,g)(x)&=
\frac{qr(x)}{(1-q)x}\bigl(f(x)g(qx)-f(qx)g(x)\bigr)\\
&= qr(x)\bigl((D_qf)(x)g(x)-f(x)(D_qg)(x)\bigr)
\end{split}
\end{equation}
for all $x\in I$.
\end{lem}
%%%%%%%%%%%%%%%%%%%%%%%%%%%%%%%%%%%%%%%%%%%%%%%%%%%%%%%%%%%%%%%%%%%%%%%%%%%%%%
\begin{proof}
Since $p$ and $r$ (respectively $A$ and $B$) are real valued on the
$q$-interval $I$, we may restrict the proof to real valued functions $f$ and
$g$. Then we derive by a direct
computation using Lemma \ref{self-adjointform} that
\begin{equation*}
\begin{split}
 \Bigl\{\bigl(Lf\bigr)(x)g(x)-
f(x)\bigl(Lg\bigr)(x)\Bigr\}&\frac{(1-q)x}{p(x)}\\
=&\begin{cases}
W(f,g)(q^{-1}x)-W(f,g)(x),\,\,\, &x\in (-1,\infty(z))_q,\\
-W(f,g)(-1), &x=-1.
\end{cases}
\end{split}
\end{equation*}
The lemma is now an easy consequence of this formula since the finite
sums become telescoping.
\end{proof}
%%%%%%%%%%%%%%%%%%%%%%%%%%%%%%%%%%%%%%%%%%%%%%%%%%%%%%%%%%%%%%%%

In the remainder of this section, we use some standard terminology on
unbounded linear operators, see  \cite[Chapter XII]{DS} and
\cite[Chapter 13]{R}.
Let ${\mathcal{D}}\subset {\mathcal{H}}$ be a dense linear space satisfying 
$L\bigl({\mathcal{D}}\bigr)\subset {\mathcal{H}}$.
Then ${\mathcal{D}}$ may be considered as a domain of definition for
the unbounded operator $L$ on ${\mathcal{H}}$. Lemma \ref{stoktermen}
and \eqref{lim} then show that $(L,{\mathcal{D}})$ is a densely defined
symmetric operator
if and only if for all $f,g\in {\mathcal{D}}$, 
the three limits
\[ \lim_{l\to -\infty}W(f,{\overline{g}})(zq^{l-1}), \quad
\lim_{m\to\infty}W(f,{\overline{g}})(zq^m),\quad
\lim_{k\to\infty}W(f,{\overline{g}})(-q^k)
\] 
exist, and  
\[
\lim_{l\to-\infty}W(f,{\overline{g}})(zq^{l-1})
-\lim_{m\to\infty}W(f,{\overline{g}})(zq^m)+\lim_{k\to\infty}
W(f,{\overline{g}})(-q^k)=0.
\]
The following lemma
determines the behaviour of
the Wronskian in the limit to infinity.

%%%%%%%%%%%%%%%%%%%%%%%%%%%%%%%%%%%%%%%%%%%%%%%%%%%%%%%%%%%%%%%%%%%
\begin{lem}\label{Winfty}
Let $f,g\in {\mathcal{H}}$, then $\lim_{m\to \infty}W(f,g)(zq^{-m})=0$.
\end{lem}
%%%%%%%%%%%%%%%%%%%%%%%%%%%%%%%%%%%%%%%%%%%%%%%%%%%%%%%%%%%%%%%%%%
\begin{proof}
Using the formula
\begin{equation}\label{qshifttransform}
\bigl(xq^{1-m};q\bigr)_{\infty}=(-x)^mq^{-m(m-1)/2}\bigl(1/x;q\bigr)_m
\bigl(qx;q\bigr)_{\infty},\qquad x\in\CC^*,
\end{equation}
it follows that the behaviour of the 
weights of $\langle .,.\rangle$ at infinity
is given by 
\begin{equation}\label{pinfty}
\frac{(1-q)zq^{-m}}{p(zq^{-m})}=K
a^{-2m}\left(1+{\mathcal{O}}(q^m)\right),\qquad m\to\infty,
\end{equation}
where $K$ is the positive constant
\begin{equation}\label{Kconstant}
K=(1-q)z\frac{\theta(-bcz,-qz)}
{\theta(-abz,-acz)}.
\end{equation}
This implies that $\lim_{m\to\infty}a^{-m}f(zq^{-m})=0$ for 
all $f\in {\mathcal{H}}$. 
The proof follows now from the definition of the Wronskian 
\eqref{Wronskian} and the fact that
\begin{equation}\label{rinfty}
\frac{qr(zq^{-m})}{(1-q)zq^{-m}}=Ka^{2-2m}\left(1+{\mathcal{O}}(q^m)\right),
\qquad m\to\infty.
\end{equation}
\end{proof}
%%%%%%%%%%%%%%%%%%%%%%%%%%%%%%%%%%%%%%%%%%%%%%%%%%%%%%%%%%%%%%%%%%

It follows from Lemma \ref{Winfty}
that the domains ${\mathcal{D}}\subset {\mathcal{H}}$ 
for which $(L,{\mathcal{D}})$ is symmetric are determined by
vanishing properties of the Wronskian at the origin. Before we define
suitable domains of definition for $L$ explicitly, 
we first introduce some convenient notations.
For $f\in {\mathcal{F}}(I)$, we define
\begin{equation}
\begin{split}
f(0^+)&=\lim_{k\to\infty}f(zq^k),\quad
f(0^-)=\lim_{k\to\infty}f(-q^k),\\
f'(0^+)&=\lim_{k\to\infty}\bigl(D_qf\bigr)(zq^k),\quad
f'(0^-)=\lim_{k\to\infty}\bigl(D_qf\bigr)(-q^k)
\end{split}
\end{equation}
provided that the limits exist. {}From now on we tacitly assume that the limits 
exist whenever we write $f(0^+)$, $f(0^-)$ etc. 
Let $\TT=\{ \alpha\in \CC \, | \, |\alpha|=1 \}$ 
be the unit circle in the complex plane.
%%%%%%%%%%%%%%%%%%%%%%%%%%%%%%%%%%%%%%%%%%%%%%%%%%%%%%%%%%%%%%%%%%%%%
\begin{Def}\label{Defdomain}
Let $\alpha\in \TT$. We write ${\mathcal{D}}_{\alpha}\subset {\mathcal{H}}$ 
for the subspace of functions $f\in {\mathcal{H}}$ satisfying 
$Lf\in {\mathcal{H}}$, $f(0^+)=\alpha f(0^-)$ and $f'(0^+)=\alpha f'(0^-)$.
\end{Def}
%%%%%%%%%%%%%%%%%%%%%%%%%%%%%%%%%%%%%%%%%%%%%%%%%%%%%%%%%%%%%%%%%%%%%

Observe that ${\mathcal{D}}_{\alpha}$ contains the functions
with finite
support, hence ${\mathcal{D}}_{\alpha}\subset {\mathcal{H}}$ is dense.
%%%%%%%%%%%%%%%%%%%%%%%%%%%%%%%%%%%%%%%%%%%%%%%%%%%%%%%%%%%%%%%%%%%%%
\begin{lem}\label{suitable}
Let $\alpha\in\TT$.

{\bf (i)} There exists a function $f\in {\mathcal{D}}_{\alpha}$
such that $f(0^+)\not=0$ and $\bigl(D_qf)(zq^k)=0$, 
$\bigl(D_qf\bigr)(-q^l)=0$ for $k,l\in\ZZ_+$.

{\bf (ii)} There exists a function 
$g\in {\mathcal{D}}_{\alpha}$ such that $g(0^+)=0$ and
$g'(0^+)\not=0$. 

{\bf (iii)} ${\mathcal{D}}_{\alpha}\not={\mathcal{D}}_{\beta}$ for
$\alpha,\beta\in \TT$ with $\alpha\not=\beta$.
\end{lem}
%%%%%%%%%%%%%%%%%%%%%%%%%%%%%%%%%%%%%%%%%%%%%%%%%%%%%%%%%%%%%%%%%%%%%%
\begin{proof}
Part {\bf (iii)} follows immediately from {\bf (i)} and 
Definition \ref{Defdomain}.
For the proof of {\bf (i)}, we observe first that the weights of 
$\langle .,.\rangle$ around zero behave like
\begin{equation}\label{asymtozero}
\frac{(1-q)q^k}{p(-q^k)}={\mathcal{O}}(q^k), \quad 
\frac{(1-q)zq^k}{p(zq^k)} =
{\mathcal{O}}(q^k),\qquad k\to\infty.
\end{equation}
Define now the function $f\in {\mathcal{F}}(I)$ by
$f(x)=1$ if $-1\leq x<0$, $f(x)=\alpha$ if $0<x\leq z$ and
$f(x)=0$ if $x>z$.
Then $f\in {\mathcal{H}}$ follows from  \eqref{asymtozero}, and
$Lf\in {\mathcal{H}}$ since 
$(Lf)(x)=0$ if $x\not\in \{z,q^{-1}z\}$. 
Furthermore, $f(0^+)=\alpha=\alpha f(0^-)$ and
$f'(0^+)=0=\alpha f'(0^-)$, hence $f\in {\mathcal{D}}_{\alpha}$.
By construction we have $f(0^+)\not=0$ and $\bigl(D_qf)(zq^k)=0$,
$\bigl(D_qf\bigr)(-q^l)=0$ for $k,l\in\ZZ_+$.

For the proof of {\bf (ii)} we define the function $g\in {\mathcal{F}}(I)$ by
$g(x)=x$ if $-1\leq x<0$, $g(x)=\alpha x$ if $0<x\leq z$, and
$g(x)=0$ if $x>z$. Then $g\in {\mathcal{H}}$. Furthermore, 
$g(0^+)=0=\alpha g(0^-)$ 
and $g'(0^+)=\alpha=\alpha g'(0^-)$.
So it remains to show that $Lg\in {\mathcal{H}}$. 
Using the explicit expression \eqref{equiv2}
for the $q$-difference operator $L$, 
we see that $\bigl(Lg\bigr)(-q^k)={\mathcal{O}}(1)$, 
$\bigl(Lg\bigr)(zq^k)={\mathcal{O}}(1)$
as $k\to\infty$. Combined with \eqref{asymtozero}, 
it follows that $Lg\in {\mathcal{H}}$.
\end{proof}
%%%%%%%%%%%%%%%%%%%%%%%%%%%%%%%%%%%%%%%%%%%%%%%%%%%%%%%%%%%%%%%%%%%%%%%
\begin{lem}\label{symmetryoperator}
Let $\alpha\in \TT$, then
$(L,{\mathcal{D}}_{\alpha})$ is a symmetric operator.
\end{lem}
%%%%%%%%%%%%%%%%%%%%%%%%%%%%%%%%%%%%%%%%%%%%%%%%%%%%%%%%%%%%%%%%%%%%%%%%
\begin{proof}
We have to show that
\[
\langle Lf,g\rangle-\langle f,Lg\rangle=0, \quad 
\forall f,g\in {\mathcal{D}}_{\alpha}.
\]
By Lemma \ref{stoktermen}, Lemma \ref{Winfty}, \eqref{lim} and the
fact that $h, Lh\in {\mathcal{H}}$ for $h\in {\mathcal{D}}_{\alpha}$, 
it suffices to show that $W(f,{\overline{g}})(0^+)=W(f,{\overline{g}})(0^-)$
for all  $f,g\in {\mathcal{D}}_{\alpha}$.
Since $r(0^+)=(1-q)^2/qbc=r(0^-)$, we have for all 
$f,g\in {\mathcal{D}}_{\alpha}$,
\begin{equation*}
\begin{split}
W(f,{\overline{g}})(0^+)=&
\frac{(1-q)^2}{bc}
\bigl(f'(0^+){\overline{g(0^+)}}-f(0^+){\overline{g'(0^+)}}\bigr)\\
=&\frac{(1-q)^2}{bc}\alpha{\overline{\alpha}}
\bigl(f'(0^-){\overline{g(0^-)}}-f(0^-){\overline{g'(0^-)}}\bigr)=
W(f,{\overline{g}})(0^-),
\end{split}
\end{equation*}
as desired.
\end{proof}
%%%%%%%%%%%%%%%%%%%%%%%%%%%%%%%%%%%%%%%%%%%%%%%%%%%%%%%%%%%%%%%%%%%%%

For $\alpha\in \TT$, we write $(L^*, {\mathcal{D}}_{\alpha}^*)$
for the adjoint of the operator $(L,{\mathcal{D}}_{\alpha})$.
Since $(L,{\mathcal{D}}_{\alpha})$ is a symmetric operator,
we have $(L,{\mathcal{D}}_{\alpha})\subset (L^*,{\mathcal{D}}_{\alpha}^*)$.

Recall that $L$ was initially defined as a linear operator on
the linear space ${\mathcal{F}}(I)$ of complex-valued functions
on $I$. In particular, $L$ is well defined  as a linear map 
$L:{\mathcal{D}}_{\alpha}^*\rightarrow {\mathcal{F}}(I)$ by restriction. 
We claim that 
\begin{equation}\label{connect*}
L^*=L|_{{\mathcal{D}}_{\alpha}^*}.
\end{equation}
To prove the claim, we observe that 
\begin{equation*}
\langle Lf,g\rangle=\langle f, Lg\rangle
\end{equation*}
for functions $f,g\in {\mathcal{F}}(I)$ such that $f$ has finite support
(compare with the proof of Lemma \ref{stoktermen}).
For $g\in {\mathcal{D}}_{\alpha}^*$ this implies that
$\langle f, Lg\rangle=\langle f, L^*g\rangle$ for functions
$f$ with finite support. Applying this formula with
a non-zero function $f$ having support in one  
point $x\in I$, we arrive at 
$\bigl(Lg\bigr)(x)=\bigl(L^*g\bigr)(x)$. The claim \eqref{connect*}
follows, since  $x\in I$ can be chosen arbitrarily.
%%%%%%%%%%%%%%%%%%%%%%%%%%%%%%%%%%%%%%%%%%%%%%%%%%%%%%%%%%%%%%%%%%%
\begin{prop}\label{s}
Let $\alpha\in\TT$, then
$(L,{\mathcal{D}}_{\alpha})$ is self-adjoint,
i.e. ${\mathcal{D}}_{\alpha}={\mathcal{D}}_{\alpha}^*$.
\end{prop}
%%%%%%%%%%%%%%%%%%%%%%%%%%%%%%%%%%%%%%%%%%%%%%%%%%%%%%%%%%%%%%%%%%%
\begin{proof}
Since $(L,{\mathcal{D}}_{\alpha})$ is symmetric, it suffices to prove the
inclusion ${\mathcal{D}}_{\alpha}^*\subset {\mathcal{D}}_{\alpha}$.
Let $g\in {\mathcal{D}}_{\alpha}^*$. Then by \eqref{connect*}, 
$Lg=L^*g\in {\mathcal{H}}$.
For any $f\in {\mathcal{D}}_{\alpha}$ we have, using \eqref{connect*}, 
Lemma \ref{stoktermen}, Lemma \ref{Winfty} and \eqref{lim}, that
\[0=\langle Lf,g\rangle-\langle f, L^*g\rangle
=W\bigl(f,{\overline{g}}\bigr)(0^-)
-W\bigl(f,{\overline{g}}\bigr)(0^+).
\]
Now using Lemma \ref{suitable},
one derives from the existence of the limits 
\begin{equation}
\begin{split}
W\bigl(f,{\overline{g}}\bigr)(0^+)&=\frac{(1-q)^2}{bc}\lim_{k\to\infty}
\bigl(\bigl(D_qf\bigr)(zq^k){\overline{g(zq^k)}}-
f(zq^k){\overline{\bigl(D_qg\bigr)(zq^k)}}\bigr)\nonumber\\
W\bigl(f,{\overline{g}}\bigr)(0^-)&=\frac{(1-q)^2}{bc}\lim_{k\to\infty}
\bigl(\bigl(D_qf\bigr)(-q^k){\overline{g(-q^k)}}-
f(-q^k){\overline{\bigl(D_qg\bigr)(-q^k)}}\bigr)\nonumber
\end{split}
\end{equation}
for all $f\in {\mathcal{D}}_{\alpha}$
that $g(0^+)=\alpha g(0^-)$ and $g'(0^+)=\alpha g'(0^-)$. It follows that 
$g\in {\mathcal{D}}_{\alpha}$, as desired.
\end{proof}
%%%%%%%%%%%%%%%%%%%%%%%%%%%%%%%%%%%%%%%%%%%%%%%%%%%%%%%%%%%%%%%%

Observe that the 
operator $L\in \hbox{End}\bigl({\mathcal{F}}(I)\bigr)$
preserves the subspace of functions $f\in {\mathcal{F}}(I)$
with support in $I_-=[-1,0]_q=\{-q^k\}_{k\in\ZZ_+}$, as well as 
the subspace of functions $f\in {\mathcal{F}}(I)$
with support in $I_+=[0,\infty(z))_q=\{zq^l\}_{l\in\ZZ}$.
We denote $L_-$ (respectively $L_+$) for the restriction of $L$ to
functions $f\in {\mathcal{F}}(I)$ with support in $I_-$
(respectively $I_+$) and we write
${\mathcal{H}}={\mathcal{H}}^-\oplus {\mathcal{H}}^+$
for the orthogonal direct sum decomposition 
where ${\mathcal{H}}^{\pm}$ are the functions $f\in {\mathcal{H}}$
with support in $I_\pm$.

We end this section by indicating 
the relation of this splitting of $L$
with the possible choices of domains for $L$.
The following results will not be used in the remainder of the paper, 
so we omit detailed proofs.

We start with the symmetric operator $(L_+,{\mathcal{D}}_{fin}^+)$
on ${\mathcal{H}}_+$, where ${\mathcal{D}}_{fin}^+\subset {\mathcal{H}}_+$ 
is the subspace of functions
with finite support. Its minimal closure is given by $(L,{\mathcal{D}}_0^+)$,
where 
\[{\mathcal{D}}_0^{+}=\{f \in {\mathcal{H}}^{+} \, | \, 
L_{+}f\in {\mathcal{H}}^{+},\,\,\, f(0^{+})=f'(0^{+})=0 \},
\]
and the deficiency indices of $(L,{\mathcal{D}}_0^+)$ are $(1,1)$.

The spectral properties of the self-adjoint extensions
of $(L_{+},{\mathcal{D}}_0^{+})$ 
are highly sensitive with respect to the choice of 
self-adjoint extension. Furthermore, the spectral analysis of $L_+$ with
respect to a fixed
choice of self-adjoint extension seems to result 
in a rather implicit description of
the spectrum and of the resolution of the identity.

In order to obtain better spectral properties, we ``blow up'' the Hilbert
space ${\mathcal{H}}_+$ to 
${\mathcal{H}}={\mathcal{H}}_-\oplus {\mathcal{H}}_+$ and we regard
$L$ now as an unbounded operator on ${\mathcal{H}}$ with domain 
${\mathcal{D}}_{fin}$ given by the functions in ${\mathcal{H}}$ with 
finite support. Its minimal closure is given by $(L,{\mathcal{D}}_0)$,
where
\[
{\mathcal{D}}_0=\{f\in {\mathcal{H}} \, | \, Lf\in {\mathcal{H}},
\,\,\,\, f(0^{+})=f(0^-)=0=f'(0^{-})=f'(0^+) \}.
\]
The closed, symmetric operator 
$(L,{\mathcal{D}}_0)$ has deficiency indices $(2,2)$.
By the increase of deficiency indices we have a larger choice of
self-adjoint extensions for $(L,{\mathcal{D}}_0)$ than for
$(L_+,{\mathcal{D}}_0^+)$. In particular, the self-adjoint extensions
$(L,{\mathcal{D}}_{\alpha})$ ($\alpha\in\TT$) are self-adjoint extensions
for which $(L_+,{\mathcal{D}}_{\alpha}\cap {\mathcal{H}}_+)$
are {\it not} self-adjoint.

It turns out that the spectral analysis of $(L,{\mathcal{D}}_{\alpha})$
is essentially
independent of $\alpha\in \TT$, so we restrict attention in this paper
to the unbounded self-adjoint operator $L$ 
on ${\mathcal{H}}$ with domain of definition
${\mathcal{D}}={\mathcal{D}}_1\subset {\mathcal{H}}$. 
In particular, any function 
$f\in {\mathcal{D}}$ is {\it continuously differentiable}
at the origin, i.e. $f$ satisfies $f(0^+)=f(0^-)$ and  $f'(0^+)=f'(0^-)$.

In the next sections we show that $(L,{\mathcal{D}})$ is 
a self-adjoint extension of $(L,{\mathcal{D}}_0)$ for which the spectral
analysis can be derived in a very explicit manner.
The reason is that for a given generic eigenvalue, one has two
linear independent eigenfunctions of $L$ which are explicitly given 
in terms of basic hypergeometric series and 
which are continuously differentiable at
the origin. The corresponding Wronskian, as well as their Wronskian 
with the asymptotically free solutions of the corresponding eigenvalue
equation, can be computed explicitly. This 
allows us to derive the spectral properties in a very explicit manner.

Let us finally make some remarks on the ``blowing up'' procedure of
the Hilbert space ${\mathcal{H}}_+$ as described in the previous paragraphs.
It is a well known principle in harmonic analysis 
that the spectral properties of the unbounded self-adjoint operator
under consideration
is essentially determined  
by the behaviour of the eigenfunctions at infinity (along the
support of the measure).
The fact that ${\mathcal{H}}_-$ is a $L^2$-space of functions supported
on a compact space implies that the situation at
infinity does not alter
in the enlarged setting. Hence, the spectral properties of the
extension of $L_+$ to the unbounded operator $L$ on 
${\mathcal{H}}$ should be essentially determined by spectral properties
of $L_+$. We will justify this principle in this paper
by deriving most of the spectral properties of $L$ from the spectral
properties of $L_+$.

The blowing up of the Hilbert space ${\mathcal{H}}_+$
is canonical in the sense that the 
extra piece ${\mathcal{H}}_-$ added to
the Hilbert space ${\mathcal{H}}_+$ is essentially determined by the
following three properties:
\begin{enumerate}
\item[\bf{(i)}] 
The support of the measure of ${\mathcal{H}}^-$ is a {\it bounded}
$q$-interval,
\item[\bf{(ii)}]
The operator $(L,{\mathcal{D}}_{fin}^-)$ on ${\mathcal{H}}^-$ is symmetric, 
where ${\mathcal{D}}_{fin}^-$ consists of 
functions in ${\mathcal{H}}^-$
with finite support, 
\item[\bf{(iii)}]
The weight function is continuous at
the origin. 
\end{enumerate} 
Indeed, the condition that the support of the measure of ${\mathcal{H}}^-$
is a bounded  $q$-interval implies that the support is of the form
$[y,0]_q=\{yq^k\}_{k\in\ZZ_+}$ with $B(y)=0$. 
This leads to the possibilities 
$y=-1$ or $y=-bc/q$. Choosing $y=-1$, we derive from
condition {\bf (ii)} that
the weight function corresponding to ${\mathcal{H}}^-$ 
is given by $1/p(\cdot)$ up to a positive constant, while condition {\bf (iii)}
implies that the constant is one.
Since $p(x;a,b,c)=p(bcx/q;a,q/b,q/c)$ and a similar property holds
for $L$, see Remark \ref{symmetryrem}, we may take
$y=-1$ without loss of generality.

%%%%%%%%%%%%%%%%%%%%%%%%%%%%%%%%%%%%%%%%%%%%%%%%%%%%%%%%%%%%%%%%%
%%                                                             %%
%%                   Eigenfunctions of $L$                     %%
%%                                                             %%
%%%%%%%%%%%%%%%%%%%%%%%%%%%%%%%%%%%%%%%%%%%%%%%%%%%%%%%%%%%%%%%%%
\section{Eigenfunctions of $L$}\label{eigenf}

In sections \ref{eigenf}--\ref{PID} we 
consider $L$ with domain of definition ${\mathcal{D}}={\mathcal{D}}_1$.
In particular, any function $f$ in the domain of definition ${\mathcal{D}}$ 
is continuously differentiable at the origin. 
We set $f(0)$ (respectively $f'(0)$) for the common limits
$f(0^\pm)$ (respectively $f'(0^\pm)$).

Furthermore, we need to restrict sometimes 
the choice of parameters $(a,b,c)\in V$ to a dense subdomain $V_z^{gen}$,
which is defined by
\begin{equation}\label{con3}
V_z^{gen}=\{(a,b,c)\in V \, | \, a^2,b^2,c^2,ab,ac,bc,a/b,a/c,a^2b^2c^2z^2
\not\in \{q^k\}_{k\in\ZZ}\}.
\end{equation}

In this section we consider eigenfunctions of the linear operators
$L$ and $L_{\pm}$, and study their
behaviour at the origin.
We write $(-1,0]_q=I_-\setminus \{-1\}=[-q,0]_q$, and define
for $\mu\in\CC$ the linear spaces
\begin{equation}\label{subspaces}
\begin{split}
V_{\mu}^{\pm}&=\{\, f: I_{\pm} \to\CC \, | \, L_{\pm}f=\mu f \,\,
\hbox{on} \,\, I_{\pm} \,\},\\
{\hat{V}}_{\mu}^{-}&=\{\, f: I_{-} \to\CC \, | \, L_{-}f=\mu f \,\,
\hbox{on} \,\, (-1,0]_q \,\},\\
V_{\mu}&=\{ f\in {\mathcal{F}}(I) \, | \, Lf=\mu f \, \hbox{ on }\, 
(-1,\infty(z))_q,\, f\hbox{ cont. differentiable in }  0 \}.
\end{split}
\end{equation}
Observe that $V_{\mu}^-\subset {\hat{V}}_{\mu}^{-}$.
%%%%%%%%%%%%%%%%%%%%%%%%%%%%%%%%%%%%%%%%%%%%%%%%%%%%%%%%%%%%%%
\begin{lem}\label{generaleigenfunction}
Let $\mu\in\CC$.\\
{\bf (i)} $\hbox{dim}\bigl(V_{\mu}^+\bigr)=2$,
$\hbox{dim}\bigl(V_{\mu}^-\bigr)=1$
and  $\hbox{dim}\bigl({\hat{V}}_{\mu}^{-}\bigr)=2$.
\\
{\bf (ii)} If $f_1,f_2\in V_{\mu}^+$ \textup{(}respectively
$f_1,f_2\in {\hat{V}}_{\mu}^{-}$\textup{)}, then 
$W(f_1,f_2)$ is constant on $I_{+}$ 
\textup{(}respectively constant on $I_-$\textup{)}.\\
{\bf (iii)} If $f_1, f_2\in V_{\mu}$, then 
$W(f_1,f_2)\in {\mathcal{F}}(I)$ is constant on $I$.
\end{lem} 
%%%%%%%%%%%%%%%%%%%%%%%%%%%%%%%%%%%%%%%%%%%%%%%%%%%%%%%%%%%%%%%%
\begin{proof}
{\bf (i)} We start with computing $\hbox{dim}\bigl(V_{\mu}^-\bigl)$.
The coefficient $A(x)$ in the expression 
\eqref{equiv2} of $L$ is non-zero for $x\in I_-$ since $(a,b,c)\in V$.
It follows that solutions of $L_-f=\mu f$ on $I_-$ are in one to one
correspondence with solutions $(a_k)_{k\in\ZZ_+}$ 
of a recurrence relation of the form
\begin{equation*}
\mu a_k= \alpha_ka_{k-1}+\beta_ka_k+\gamma_ka_{k+1},\qquad
\gamma_k\not=0,\,\,k\in\ZZ_+  
\end{equation*}
with $a_{-1}=0$ by definition. The correspondence is obtained by associating 
the sequence $a_k=f(-q^k)$ ($k\in\ZZ_+$) to $f\in V_{\mu}^-$. The
coefficients of the recurrence relation are then given by
$\alpha_k=-B(-q^k)$, $\beta_k=-A(-q^k)-B(-q^k)$ and $\gamma_k=A(-q^k)$.
Any solution of such a recurrence relation is uniquely 
determined by $a_0=f(-1)\in\CC$, hence $\hbox{dim}\bigl(V_{\mu}^-\bigl)=1$.

Then $\hbox{dim}\bigl({\hat{V}}_{\mu}^{-}\bigl)=2$ follows
from
the fact that functions $f\in {\hat{V}}_{\mu}^-$ are in one to one 
correspondence with solutions $(a_k)_{k\in\ZZ_+}$ of the 
recurrence relation
\begin{equation*}
\mu a_k=\alpha_ka_{k-1}+\beta_ka_k+\gamma_ka_{k+1},\qquad
\gamma_k\not=0,\,\, k\in \NN.
\end{equation*}

Finally, in order to show that $\hbox{dim}\bigl(V_{\mu}^-\bigl)=2$,
observe that solutions of $L_+f=\mu f$ on $I_+$ are in one to one 
correspondence
with solutions  $(a_k)_{k\in\ZZ}$ of a double infinite recurrence
relation of the form
\[ \mu a_k=\alpha_ka_{k-1}+\beta_ka_k+\gamma_ka_{k+1},\qquad
\alpha_k\not=0,\, \gamma_k\not=0,
\,\,k\in\ZZ,\]
since $A(x)$ and $B(x)$ \eqref{ABbig} in the expression 
\eqref{equiv2} of $L$ are non-zero for $x\in I_+$.
This is a two-dimensional space.\\
{\bf (ii)} Let $f_1, f_2\in {\mathcal{F}}(I)$. 
Using Lemma \ref{self-adjointform}, the product rule
for the $q$-derivative
\begin{equation}\label{productrule}
\bigl(D_q(fg)\bigr)(x)=\bigl(D_qf\bigr)(x)g(x)+f(qx)\bigl(D_qg\bigr)(x),
\end{equation}
and the second equality of \eqref{Wronskian},
we have for all $x\in (-1,\infty(z))_q$,
\begin{equation}\label{consta}
\begin{split}
q^{-1}p(x)&\bigl(D_qW(f_1,f_2)\bigr)(q^{-1}x)\\
&=\bigl(Lf_1\bigr)(x)f_2(x)+p(x)r(q^{-1}x)\bigl(D_qf_1\bigr)(q^{-1}x)
\bigl(D_qf_2\bigr)(q^{-1}x)\\
&-\bigl(Lf_2\bigr)(x)f_1(x)-p(x)r(q^{-1}x)\bigl(D_qf_2\bigr)(q^{-1}x)
\bigl(D_qf_1\bigr)(q^{-1}x)\\
&=\bigl(Lf_1\bigr)(x)f_2(x)-\bigl(Lf_2\bigr)(x)f_1(x).
\end{split}
\end{equation}
If $f_1,f_2\in {\hat{V}}_{\mu}^-$, then it follows from \eqref{consta} that
$\bigl(D_qW(f_1,f_2)\bigr)(x)=0$ for all $x\in (-1,0]_q$ since
$p(x)\not=0$ $(x\in I_-)$, hence $W(f_1,f_2)$ is constant on
$I_-$. A similar argument shows that $W(f_1,f_2)$ is constant on $I_+$
if $f_1,f_2\in V_{\mu}^+$.\\
{\bf (iii)} Let $f_1,f_2\in V_{\mu}$.
The continuously differentiability of
the $f_j$'s at the origin yields that $W(f_1,f_2)(0^+)=W(f_1,f_2)(0^-)$.
Indeed, both sides are equal to $qr(0)\bigl(f_1'(0)f_2(0)-
f_1(0)f_2'(0)\bigr)$ (compare with the proof of
Lemma \ref{symmetryoperator}). The result follows now from {\bf (ii)}.
\end{proof}
%%%%%%%%%%%%%%%%%%%%%%%%%%%%%%%%%%%%%%%%%%%%%%%%%%%%%%%%%%%%%%%

Observe that $f|_{I_-}\in {\hat{V}}_{\mu}^-$ and 
$f|_{I_+}\in V_{\mu}^+$ for $f\in V_{\mu}$.
%%%%%%%%%%%%%%%%%%%%%%%%%%%%%%%%%%%%%%%%%%%%%%%%%%%%%%%%%%%%%%%%
\begin{prop}\label{generaleigenfunctionprop}
We have ${\hbox{dim}}(V_{\mu})\leq 2$. Furthermore, the
following three statements are equivalent:\\
{\bf (i)} The linear map $f\mapsto f|_{I_-}:\,\, 
V_{\mu}\to {\hat{V}}_{\mu}^-$ is a bijection.\\
{\bf (ii)} The linear map $f\mapsto f|_{I_+}: V_{\mu}\to
V_{\mu}^+$ is a bijection.\\
{\bf (iii)} $\hbox{dim}(V_{\mu})=2$.
\end{prop}
%%%%%%%%%%%%%%%%%%%%%%%%%%%%%%%%%%%%%%%%%%%%%%%%%%%%%%%%%%%%%%%%%%
\begin{proof}
Suppose that $\hbox{dim}(V_{\mu})\geq 2$, then we claim that
the map $f\mapsto f|_{I_+}:\,\, V_{\mu}\to V_{\mu}^+$ is injective. 
Suppose that the map is not injective. 
Then there exist
two linearly independent functions $f_1,f_2\in V_{\mu}$ such that
$f_1|_{I_+}, f_2|_{I_+}\in V_{\mu}^+$ are linearly dependent.
The Wronskian $W(f_1,f_2)$ is not identically zero as function on $I$ by the 
linear independence of $f_1$ and $f_2$. Since $W(f_1,f_2)$ is constant on $I$ by 
Lemma \ref{generaleigenfunction}{\bf (iii)}, it follows that
$W(f_1,f_2)$ is not identically zero on $I_+$, which contradicts the
linear dependence of $f_1|_{I_+}$ and $f_2|_{I_+}$.

So if $\hbox{dim}(V_{\mu})\geq 2$, then it follows from
Lemma \ref{generaleigenfunction}{\bf (i)} and from 
the injectivity of the map $f\mapsto f|_{I_+}:\,\, V_{\mu}\to
V_{\mu}^+$ that $\hbox{dim}(V_{\mu})=2$.
This proves that $\hbox{dim}(V_{\mu})\leq 2$ for all $\mu\in \CC$, 
and it proves the implication ${\bf (iii)}\Rightarrow {\bf (ii)}$.
The implication ${\bf (iii)}\Rightarrow {\bf (i)}$ is proved in a similar
manner, while the implications
${\bf (i)}\Rightarrow {\bf (iii)}$ and ${\bf (ii)}\Rightarrow
{\bf (iii)}$ are immediate consequences of Lemma 
\ref{generaleigenfunction}{\bf (i)}.
\end{proof}
%%%%%%%%%%%%%%%%%%%%%%%%%%%%%%%%%%%%%%%%%%%%%%%%%%%%%%%%%%%%%%%%

%%%%%%%%%%%%%%%%%%%%%%%%%%%%%%%%%%%%%%%%%%%%%%%%%%%%%%%%%%%%%%%%%
\begin{cor}\label{extensionpossible}
If $\hbox{dim}(V_{\mu})=2$, then 
any $g\in {\hat{V}}_{\mu}^-$ \textup{(}respectively $g\in V_{\mu}^+$\textup{)} 
extends uniquely to a function $g\in V_{\mu}$.
\end{cor}
%%%%%%%%%%%%%%%%%%%%%%%%%%%%%%%%%%%%%%%%%%%%%%%%%%%%%%%%%%%%%%%%%%

In the remainder of this section we introduce two explicit functions
which live in $V_{\mu}$. We define for $\gamma\in \CC^*$,
\begin{equation}\label{eigenfunction}
\mu(\gamma)=-1-a^2+a(\gamma+\gamma^{-1}),
\end{equation}
and we write
$X_n=\bigl(bcx^2/q\bigr)^n\bigl(-q/bcx,-1/x;q\bigr)_ng_{\gamma}(xq^{-n})$
where $g_{\gamma}$ is a 
function satisfying $\bigl(Lg_{\gamma}\bigr)(xq^{-n})=
\mu(\gamma)g_{\gamma}(xq^{-n})$. 
Then the $X_n$ satisfy the three term recurrence relation \cite[(2.1)]{GIM}
with parameters $(A,B,C,D,z)$ in \cite[(2.1)]{GIM} 
given by $A=-q/bcx$, $B=-1/x$, $C=-q/acx$, $D=-q/abx$ and
\[z=\frac{abcx^2}{q}(\gamma+\gamma^{-1})=
\frac{q}{ABCD\lambda_{\pm}}+\lambda_{\pm},\qquad
\lambda_{\pm}=\frac{abcx^2}{q}\gamma^{\pm 1}.
\]
In \cite{GIM} 
it is was shown that the three term recurrence relation \cite[(2.1)]{GIM}
is the one satisfied by the 
associated continuous dual $q$-Hahn polynomials. Furthermore,
several explicit solutions of the recurrence relation \cite[(2.1)]{GIM}
were derived explicitly in \cite[section 2]{GIM}.
In particular, the explicit solution \cite[(2.13)]{GIM} 
of the recurrence relation 
\cite[(2.1)]{GIM} implies that
\begin{equation}\label{phiB}
\phi_{\gamma}(x)=\phi_{\gamma}(x;a,b,c)=
{}_3\phi_2\left( \begin{array}{c}
             a\gamma, a/\gamma, -1/x \\ 
             ab, ac
            \end{array} ; q,-bcx \right),\quad |bcx|<1
\end{equation}
satisfies $\bigl(L\phi_{\gamma}\bigr)(x)=\mu(\gamma)\phi_{\gamma}(x)$
for $x\in (-1,\infty(z))_q$ with $|x|<q/bc$. Observe that
$\phi_{\gamma}$ is well defined since $(a,b,c)\in V$.

The explicit solution \cite[(2.13)]{GIM} of the recurrerence
relation \cite[(2.1)]{GIM} implies that
\begin{equation}\label{psiB}
\psi_{\gamma}(x;a,b,c)
={}_3\phi_2\left( \begin{array}{c}
             a\gamma, a/\gamma, -q/bcx \\ 
             qa/b, qa/c
            \end{array} ; q,-qx \right), \quad |qx|<1
\end{equation}
satisfies $\bigl(L\psi_{\gamma}\bigr)(x)=\mu(\gamma)\psi_{\gamma}(x)$
for $x\in (-1,\infty(z))_q$ with $|x|<1$ when $(a,b,c)\in V_z^{gen}$. 

Since $bc<1$, we have that $\phi_{\gamma}$ and 
$\psi_{\gamma}$ are well defined on $I_-$ and that they are solutions of 
$Lg=\mu(\gamma)g$ on $(-1,0]_q$. Since $\phi_{\gamma}$ and $\psi_{\gamma}$
are solutions of $\bigl(Lg\bigr)(x)=\mu(\gamma)g(x)$ for $x\in I_+$ with
$x<1$, it follows that $\phi_{\gamma}$ and $\psi_{\gamma}$ uniquely extend
to functions on $I$ satisfying the eigenvalue equation $Lg=\mu(\gamma)g$
on $(-1,\infty(z))_q$, cf. 
the proof of Lemma \ref{generaleigenfunction}{\bf (i)}.
For given $x\in I$, $\phi_{\gamma}(x)$ and $\psi_{\gamma}(x)$ depend
analytically on  $\gamma\in\CC^*$ and are invariant under 
$\gamma\leftrightarrow \gamma^{-1}$. 

The extension of $\phi_{\gamma}$ to a function on $I$
can also be obtained by the
transformation formula \cite[(3.2.10)]{GR}, which yields  
\begin{equation}\label{ancontB}
\phi_{\gamma}(x)= \frac{\bigl(a\gamma, bc,-abcx/\gamma;q\bigr)_{\infty}}
{\bigl(ab, ac, -bcx;q\bigr)_{\infty}}
{}_3\phi_2\left( \begin{array}{c}
             b/\gamma, c/\gamma, -bcx \\ 
             bc, -abcx/\gamma
            \end{array} ; q,a\gamma \right).
\end{equation}
This gives a single valued analytic continuation for
$\phi_{\gamma}(x)$ for $\gamma\in\CC^*$ with $|\gamma|<a^{-1}$ to  
$x\in \CC\setminus (-\infty, -1/bc]$. 
Observe that $(-\infty,-1/bc]\cap I=\emptyset$ since $(a,b,c)\in V$.
If $a<1$, then the analytic
continuation for $|\gamma|\geq a^{-1}$ can be obtained by 
replacing $\gamma$ by $\gamma^{-1}$ in the right hand side of
\eqref{ancontB}. If $a\geq 1$, then $b<1$, and we can use \cite[(3.2.7)]{GR}
to obtain
\begin{equation}\label{ancontB2}
\phi_{\gamma}(x)= \frac{\bigl(b\gamma, -abcx/\gamma;q\bigr)_{\infty}}
{\bigl(ab, -bcx;q\bigr)_{\infty}}
{}_3\phi_2\left( \begin{array}{c}
             a/\gamma, c/\gamma, -acx \\ 
             ac, -abcx/\gamma
            \end{array} ; q,b\gamma \right),
\end{equation}
which yields then the explicit expression for the analytic continuation of
$\phi_{\gamma}$ ($\gamma\in \CC^*$) in a similar manner as for $a<1$.
A similar remark holds for the solution $\psi_{\gamma}$.
%%%%%%%%%%%%%%%%%%%%%%%%%%%%%%%%%%%%%%%%%%%%%%%%%%%%%%%%%%%%%%%%%%%
\begin{rem}\label{symmetryrem}
If $f(\cdot)=f(\cdot;a,b,c)$
satisfies the eigenvalue equation $Lf=\mu(\gamma)f$, then 
$g(x)=f(bcx/q;a,q/b,q/c)$ satisfies the same eigenvalue equation.
It is easy to check using the explicit expressions \eqref{phiB}
and \eqref{psiB} for $\phi_{\gamma}$ and $\psi_{\gamma}$ that 
the solutions $\phi_{\gamma}$ and $\psi_{\gamma}$ are interchanged by this
symmetry, i.e. 
\begin{equation}\label{symmetryphi}
\psi_{\gamma}(x;a,b,c)=
\phi_{\gamma}(bcx/q;a,q/b,q/c).
\end{equation} 
\end{rem}
%%%%%%%%%%%%%%%%%%%%%%%%%%%%%%%%%%%%%%%%%%%%%%%%%%%%%%%%%%%%%%%%%%%%

{}From \eqref{phiB} we see that $\phi_{\gamma}(-1)=1$ and
by direct computation we have 
\begin{equation}\label{dphi}
\bigl(D_q\phi_{\gamma}\bigr)(x;a,b,c)=\frac{bc\mu(\gamma)}
{(1-q)(1-ab)(1-ac)}\phi_{\gamma}(x;qa,b,c),\qquad x\in I.
\end{equation}
It follows that
\begin{equation}\label{valuephi}
\bigl(D_q\phi_{\gamma}\bigr)(-1)=
\frac{bc\mu(\gamma)}{(1-q)(1-ab)(1-ac)},
\end{equation}
so 
$\bigl(L\phi_{\gamma}\bigr)(-1)=\mu(\gamma)\phi_{\gamma}(-1)$ 
by \eqref{special}. We conclude that
$\phi_{\gamma}|_{I_-}\in V_{\mu(\gamma)}^-$ for all $\gamma\in \CC^*$.
For $\psi_{\gamma}$ we have the formula
\begin{equation}\label{dpsi}
\bigl(D_q\psi_{\gamma}\bigr)(x;a,b,c)=
\frac{q\mu(\gamma)}{(1-q)(1-qa/b)(1-qa/c)}
\psi_{\gamma}(x;qa,b,c).
\end{equation}
We stress already the fact that in general
we have
$\bigl(L\psi_{\gamma}\bigr)(-1)\not=
\mu(\gamma)\psi_{\gamma}(-1)$, 
so that $\psi_{\gamma}|_{I_-}\in \hat{V}_{\mu(\gamma)}\setminus
V_{\mu(\gamma)}$, see Corollary \ref{psifoutind}{\bf (iii)}.

%%%%%%%%%%%%%%%%%%%%%%%%%%%%%%%%%%%%%%%%%%%%%%%%%%%%%%%%%%%%%%%%%%%%%%%
\begin{lem}\label{oknul}
Let $\gamma\in\CC^*$. 
The function $\phi_{\gamma}\in {\mathcal{F}}(I)$ is
continuously differentiable at the origin, i.e. 
$\phi_{\gamma}\in V_{\mu(\gamma)}$. 
In fact, we have
\begin{equation}\label{nul}
\begin{split}
\phi_{\gamma}(0)&={}_2\phi_2\left( \begin{array}{c}
             a\gamma, a/\gamma \\ 
             ab, ac
            \end{array} ; q, bc \right),\\
\phi_{\gamma}'(0)&=\frac{bc\mu(\gamma)}{(1-q)(1-ab)(1-ac)}
{}_2\phi_2\left( \begin{array}{c}
            qa\gamma, qa/\gamma \\ 
             qab, qac
            \end{array} ; q, bc \right).
\end{split}
\end{equation}
For $(a,b,c)\in V_z^{gen}$ we have
$\psi_{\gamma}\in V_{\mu(\gamma)}$ and 
\[
\psi_{\gamma}(0;a,b,c)=\phi_{\gamma}(0;a,q/b,q/c),\qquad
\psi_{\gamma}'(0;a,b,c)=\frac{bc}{q}\phi_{\gamma}'(0;a,q/b,q/c),
\]
where we have extended the definition of $\phi_{\gamma}(0)$
and $\phi_{\gamma}'(0)$ to generic parameters $(a,b,c)\in \CC^{\times 3}$
by analytic continuation of the right hand sides of \eqref{nul}.
\end{lem}
%%%%%%%%%%%%%%%%%%%%%%%%%%%%%%%%%%%%%%%%%%%%%%%%%%%%%%%%%%%%%%%%%%%%%%%%
\begin{proof}
The proof for $\phi_{\gamma}$ is a direct consequence of \eqref{dphi} and
the explicit expression  \eqref{phiB} for $\phi_{\gamma}$.
The proof for $\psi_{\gamma}$ follows then
from \eqref{psiB} and \eqref{dpsi}.
\end{proof}
%%%%%%%%%%%%%%%%%%%%%%%%%%%%%%%%%%%%%%%%%%%%%%%%%%%%%%%%%%%%%%%

In section \ref{section5} we will evaluate the Wronskian 
$W\bigl(\psi_{\gamma},\phi_{\gamma}\bigr)$ explicitly,
see Proposition \ref{Wpsiphiprop}. In particular, this will give 
explicit criteria on the spectral parameter $\gamma$ for which 
we have $\hbox{dim}(V_{\mu(\gamma)})=2$, i.e. for which 
Corollary \ref{extensionpossible} is applicable.
But first we need to study yet another solution of the eigenvalue
equation $Lg=\mu(\gamma)g$, the so called asymptotic 
solution.

%%%%%%%%%%%%%%%%%%%%%%%%%%%%%%%%%%%%%%%%%%%%%%%%%%%%%%%%%%%%%%%%%%%%%%
%%                                                                  %%
%%               The asymptotic solution                            %%
%%                                                                  %%
%%%%%%%%%%%%%%%%%%%%%%%%%%%%%%%%%%%%%%%%%%%%%%%%%%%%%%%%%%%%%%%%%%%%%%

\section{The asymptotic solution}\label{A} 

In this section we determine the asymptotic solution $\Phi_{\gamma}$
of the eigenvalue equation $L_+g=\mu(\gamma)g$ on $I_+$. We furthermore 
determine the $c$-function expansion of $\phi_{\gamma}$ on $I_+$, i.e. we 
write $\phi_{\gamma}|_{I_+}$ explicitly as linear combination 
of the asymptotic solutions $\Phi_{\gamma}$ and $\Phi_{\gamma^{-1}}$.

We define singular sets by
\begin{equation}\label{Ssing}
S_{sing}^{+}=\{\pm q^{-\frac{1}{2}k}\}_{k\in\NN},\qquad
S_{sing}=\{\pm q^{\frac{1}{2}k}\}_{k\in \ZZ},
\end{equation}
and we write $S_{reg}^{+}$, $S_{reg}$ for the 
complements of these singular sets in $\CC^*$.
The $c$-function expansion of $\phi_{\gamma}$ on $I_+$ which we derive
in this section, holds for $\gamma\in S_{reg}$.

We consider for $y>0$ 
the minimal solution \cite[(2.26)]{GIM} of the three term 
recurrence relation \cite[(2.1)]{GIM} with corresponding
parameters $(A,B,C,D,\lambda_-)$ given by
\[
A=-q/bcy,\quad B=-q/aby,\quad C=-q/acy,\quad D=-1/y,\quad
\lambda_-=\frac{abcy^2}{q}\gamma.
\]
This leads 
for $\gamma\in S_{reg}^+$ with $|\gamma|<a^{-1}$ to the following 
explicit solution $\Phi^y_{\gamma}(x)=\Phi^y_{\gamma}(x;a,b,c)$ 
of $Lf=\mu(\gamma)f$ on $[0,\infty(y))_q$, 
\begin{equation}\label{PhiB}
\begin{split}
\Phi_{\gamma}^y(x)=&
\frac{\bigl(-q\gamma/ax,-q^2\gamma/abcx,a\gamma;q\bigr)_{\infty}}
{\bigl(-q/abx,-q/acx,q\gamma^2;q\bigr)_{\infty}}\\
&.\bigl(a\gamma\bigr)^{-k}{}_3\phi_2\left( \begin{array}{c}
            q\gamma/a, -q/abx, -q/acx \\ 
             -q\gamma/ax, -q^2\gamma/abcx 
            \end{array} ; q,a\gamma \right), \quad x=yq^k.
\end{split}
\end{equation}
For given $x\in [0,\infty(y))_q$, 
we have a single valued analytic extension of 
$\Phi^y_{\gamma}(x)$ to $\gamma\in S_{reg}^+$.
Indeed, for $x\in [0,\infty(y))_q$ with $x>bc/q$ 
we can apply \cite[(3.2.7)]{GR}
to arrive at
\begin{equation}\label{PhiBalt}
\begin{split}
\Phi^y_{\gamma}(x)=&
\frac{\bigl(-q/bcx,-q\gamma/ax;q\bigr)_{\infty}}
{\bigl(-q/abx,-q/acx;q\bigr)_{\infty}}\\
&.\bigl(a\gamma\bigr)^{-k}{}_3\phi_2\left( \begin{array}{c}
            q\gamma/a, b\gamma, c\gamma \\ 
             -q\gamma/ax, q\gamma^2 
            \end{array} ; q,-q/bcx\right), \quad x=yq^k,
\end{split}
\end{equation}
which gives the analytic extension in this case.
For $x\in [0,\infty(y))_q$ with $y\leq bc/q$ 
we can use the eigenvalue equation 
$L\Phi^y_{\gamma}=\mu(\gamma)\Phi^y_{\gamma}$ on $[0,\infty(y))_q$ 
to derive the analytic
extension for $\Phi^y_{\gamma}(x)$ from the analytic extension of 
$\Phi^y_{\gamma}(u)$ with $u>bc/q$. 
It follows from \eqref{PhiBalt} that 
the asymptotics to infinity of $\Phi^y_{\gamma}$ is given by
\begin{equation}\label{asymptoticsPhiB}
\Phi_{\gamma}^y(yq^{-m})=\bigl(a\gamma\bigr)^m\bigl(1+{\mathcal{O}}(q^m)\bigr),
\qquad m\to\infty.
\end{equation}
%%%%%%%%%%%%%%%%%%%%%%%%%%%%%%%%%%%%%%%%%%%%%%%%%%%%%%%%%%%%%%%%%%%%%%%%%%%
\begin{Def}\label{asymptoticsolutiondef}
Let $\gamma\in S_{reg}^+$.
We call $\Phi_{\gamma}(\cdot;a,b,c)=\Phi^z_{\gamma}(\cdot;a,b,c)$
the asymptotic solution of $L_+f=\mu(\gamma)f$ on $I_+=[0,\infty(z))_q$.
\end{Def}
%%%%%%%%%%%%%%%%%%%%%%%%%%%%%%%%%%%%%%%%%%%%%%%%%%%%%%%%%%%%%%%%%%%%%%%%%%%

Definition \ref{asymptoticsolutiondef} is justified by the following lemma.
%%%%%%%%%%%%%%%%%%%%%%%%%%%%%%%%%%%%%%%%%%%%%%%%%%%%%%%%%%%%%%%
\begin{lem}\label{l2}
Let $\mu\in \CC^*$ and 
$K=\{zq^{-k}\}_{k\in\ZZ_+}\subset I_+$.
Set
\[M_{\mu}=\{ h: I_+\to\CC  \, | \, L_+h=\mu h \,\, \hbox{ on } I_+
\,\, \hbox{ and }\,\,  \|h_K\|^2<\infty\},
\]
where $h_K\in {\mathcal{F}}(I)$ is defined to be equal to
zero on $I\setminus K$ and to be equal to $h$ on $K$.
Then $M_{\mu}=\hbox{span}\{\Phi_{\gamma}\}$ 
for all $\mu\in \CC\setminus \RR$,  
where $\gamma\in \CC^*$ is the unique non-zero complex number such that
$\mu=\mu(\gamma)$ and $|\gamma|<1$.
\end{lem} 
%%%%%%%%%%%%%%%%%%%%%%%%%%%%%%%%%%%%%%%%%%%%%%%%%%%%%%%%%%%%%
\begin{proof} 
First of all, observe that if $\gamma\in \CC^*$ and $|\gamma|<1$, then 
$\gamma\in S_{reg}^+$, hence $\Phi_{\gamma}$ is well defined.

Let $\mu\in\CC\setminus\RR$ and let 
$\gamma\in\CC^*$, $|\gamma|<1$ be such that $\mu=\mu(\gamma)$.
Using the asymptotics \eqref{pinfty} for the weights of 
the inner product $\langle .,. \rangle$,
it follows from \eqref{asymptoticsPhiB} that $\Phi_{\gamma}\in
M_{\mu(\gamma)}$ because $|\gamma|<1$.

It remains to show that $\hbox{dim}(M_{\mu})=1$. For the proof 
we use some well known results from 
the theory of the classical moment problem, see for instance 
\cite{Ak} or \cite{Simon}. 
If $f:I_+\to\CC$ satisfies $L_+f=\mu f$ on $I_+$,
then by setting 
\[
a_k=f(zq^{-k}){\sqrt{\frac{(1-q)zq^{-k}}{p(zq^{-k})}}}, \qquad k\in
\ZZ,
\]
we see that $\bigl(a_k\bigr)_{k\in\ZZ}$ satisfies the 
recurrence relation 
\begin{equation}\label{recassdualqHahn}
\alpha_{k-1}a_{k-1}+\beta_ka_k+\alpha_ka_{k+1}=\mu a_k
\end{equation}
with
\[ \alpha_k=a\sqrt{\left(1+\frac{q^{k+1}}{abz}\right)
\left(1+\frac{q^{k+1}}{acz}\right)
\left(1+\frac{q^{k+1}}{bcz}\right)\left(1+\frac{q^k}{z}\right)}
\]
and $\beta_k=-A(zq^{-k})-B(zq^{-k})$. 
It follows from \cite[Corollary 4.5]{Simon} and from the
fact that the sequence $\bigl(\alpha_k\bigr)_{k\in\ZZ_+}$ is bounded
that the Hamburger moment problem corresponding to the recurrence
relation \eqref{recassdualqHahn} for $k\in\ZZ_+$ is determined.
By \cite[Theorem 3]{Simon} this implies that
$\hbox{dim}(M_{\mu})=1$, as desired.
\end{proof}
%%%%%%%%%%%%%%%%%%%%%%%%%%%%%%%%%%%%%%%%%%%%%%%%%%%%%%%%%%%%%%%%%%%

%%%%%%%%%%%%%%%%%%%%%%%%%%%%%%%%%%%%%%%%%%%%%%%%%%%%%%%%%%%%%%%%%%%%%
It follows from Lemma \ref{l2} that 
\begin{equation}\label{symmetryPhi}
\Phi_{\gamma}^{bcz/q}(bcx/q;a,q/b,q/c)=\Phi_{\gamma}(x;a,b,c),\qquad
x\in I_+
\end{equation}
for parameters $(a,b,c)\in V$ such that $(a,q/b,q/c)\in V$.
Indeed, both sides of \eqref{symmetryPhi}, considered as function of
$x\in I_+$, are solutions of
$L_+f=\mu(\gamma)f$ on $I_+$ by Remark \ref{symmetryrem}, 
and they
have the same asymptotics to infinity by \eqref{asymptoticsPhiB}.
Formula \eqref{symmetryPhi} is also obvious from the explicit expression
\eqref{PhiB} for $\Phi^y_{\gamma}$.

%%%%%%%%%%%%%%%%%%%%%%%%%%%%%%%%%%%%%%%%%%%%%%%%%%%%%%%%%%%%%%%%%%
\begin{lem}\label{dep}
For $\gamma\in S_{reg}$ and
$x\in I_+$, we have 
\[
W\bigl(\Phi_{\gamma},\Phi_{\gamma^{-1}}\bigr)(x)=aK
\bigl(\gamma-\gamma^{-1}\bigr)\not=0, 
\]
where $K$ is the positive
constant defined by \eqref{Kconstant}. In particular, 
$\{\Phi_{\gamma}, \Phi_{\gamma^{-1}}\}$ is a basis of
$V_{\mu(\gamma)}^+$ when $\gamma\in S_{reg}$.
\end{lem}
%%%%%%%%%%%%%%%%%%%%%%%%%%%%%%%%%%%%%%%%%%%%%%%%%%%%%%%%%%%%%%%%%%%%%
\begin{proof}
The explicit expression for the Wronskian follows by computing the limit 
$\lim_{m\to\infty}W\bigl(\Phi_{\gamma},\Phi_{\gamma^{-1}}\bigr)(zq^{-m})$
using the first expression of \eqref{Wronskian} and the formulas 
\eqref{rinfty} and \eqref{asymptoticsPhiB}.
Since the Wronskian is non-zero, it follows that $\Phi_{\gamma}$ and
$\Phi_{\gamma^{-1}}$ are linear independent, hence they form a basis
of $V_{\mu(\gamma)}^+$ by Lemma 
\ref{generaleigenfunction}{\bf (i)}.
\end{proof}
%%%%%%%%%%%%%%%%%%%%%%%%%%%%%%%%%%%%%%%%%%%%%%%%%%%%%%%%%%%%%%%%%%%

It follows from Lemma \ref{dep} that 
$\phi_{\gamma}|_{I_+}$ and $\psi_{\gamma}|_{I_+}$ 
can be written uniquely as linear combination of $\Phi_{\gamma}$ 
and $\Phi_{\gamma^{-1}}$ for $\gamma\in S_{reg}$.
The corresponding coefficients can be expressed in terms
of the $c$-function $c(\gamma)=c(\gamma;a,b,c;z)$, which is defined by
\begin{equation}\label{cB}
c(\gamma)=\frac{1}{\bigl(ab, ac;q\bigr)_{\infty}\theta(-bcz)}
\frac{\bigl(a/\gamma, b/\gamma, c/\gamma;q\bigr)_{\infty}
\theta(-q/abcz\gamma)}{\bigl(1/\gamma^2;q\bigr)_{\infty}}.
\end{equation}
%%%%%%%%%%%%%%%%%%%%%%%%%%%%%%%%%%%%%%%%%%%%%%%%%%%%%%%%%%%%%%%%%%%%
\begin{prop}\label{conncoef}
Let $\gamma\in S_{reg}$.
Then we have 
\[\phi_{\gamma}(x)=c(\gamma)\Phi_{\gamma}(x)+
c(\gamma^{-1})\Phi_{\gamma^{-1}}(x), \qquad x\in I_+.
\]
The same formula holds for $\psi_{\gamma}(x)$ when $(a,b,c)\in V_z^{gen}$, 
with $c(\gamma)$ replaced by 
\begin{equation}\label{tildecB}
\begin{split}
\tilde{c}(\gamma;a,b,c;z)&=c(\gamma;a,q/b,q/c;bcz/q)\\
&=\frac{1}{\bigl(qa/b, qa/c;q\bigr)_{\infty}\theta(-1/z)}
\frac{\bigl(a/\gamma, q/b\gamma, q/c\gamma;q\bigr)_{\infty}
\theta(-1/az\gamma)}{\bigl(1/\gamma^2;q\bigr)_{\infty}}.
\end{split}
\end{equation}
\end{prop}
%%%%%%%%%%%%%%%%%%%%%%%%%%%%%%%%%%%%%%%%%%%%%%%%%%%%%%%%%%%%%%%%%%%%
\begin{proof}
We first prove the connection coefficient
formula for $\phi_{\gamma}$. Observe that $c(\gamma^{\pm 1})$ 
is well defined for $\gamma\in S_{reg}$ since $(a,b,c)\in V$.
We fix $x=zq^k$ with $k\in \ZZ$ such that $q/acx<1$. Furthermore, 
we assume that $a<1$ and we fix
$\gamma\in \CC\setminus (-\infty,0]$ such that $\gamma\in S_{reg}$
and $a<|\gamma|<1/a$. 

By the assumptions on $x$ and $\gamma$, we may apply
the three term recurrence relation \cite[(3.3.3)]{GR} with
$a\to -bcx$, $b\to b/\gamma$,
$c\to c/\gamma$, $d\to -abcx/\gamma$ and $e\to bc$. We arrive at
\begin{equation*}
\begin{split}
&{}_3\phi_2\left( \begin{array}{c}
             -bcx, b/\gamma, c/\gamma \\ 
             bc, -abcx/\gamma
            \end{array} ; q,a\gamma \right)\\
&=\frac{\bigl(c\gamma,
b\gamma,-q/bx\gamma,-q\gamma/abcx;q\bigr)_{\infty}}
{\bigl(bc,-q/abx,-q/bcx, \gamma^2;q\bigr)_{\infty}}
{}_3\phi_2\left( \begin{array}{c}
             c/\gamma, a/\gamma, q/b\gamma \\ 
             -q/bx\gamma, q/\gamma^2
            \end{array} ; q,-q/acx \right)\\
&\quad-\frac{\bigl(-q\gamma/abcx,-q\gamma/ax,b/\gamma,c/\gamma,a/\gamma
;q\bigr)_{\infty}\theta(-abcx\gamma/q)}
{\bigl(-abcx/q\gamma,bc,-q/acx,-q/abx,-q/bcx;q\bigr)_{\infty}
\theta(\gamma^2)}\\
&\qquad\qquad\qquad\qquad\qquad\qquad\qquad
\qquad.{}_3\phi_2\left( \begin{array}{c}
             q\gamma/a, -q/acx, -q/abx \\ 
             -q^2\gamma/abcx, -q\gamma/ax
            \end{array} ; q,a\gamma \right).
\end{split}
\end{equation*}
The ${}_3\phi_2$ on the left hand side is the
${}_3\phi_2$ in the expression \eqref{ancontB} for $\phi_{\gamma}$
and the second ${}_3\phi_2$ on the right hand side is the
${}_3\phi_2$ in the expression \eqref{PhiB} for $\Phi_{\gamma}$.
Again by the assumptions on $x$ and $\gamma$, we may
rewrite the first ${}_3\phi_2$ on the right hand side using
\cite[(3.2.10)]{GR} with 
$a\to c/\gamma$, $b\to a/\gamma$, $c\to q/b\gamma$,
$d\to q/\gamma^2$ and $e\to -q/bx\gamma$. This gives
\begin{equation*}
\begin{split}
&{}_3\phi_2\left( \begin{array}{c}
             c/\gamma, a/\gamma, q/b\gamma \\ 
             -q/bx\gamma, q/\gamma^2
            \end{array} ; q,-q/acx \right)\\
&\qquad\qquad=
\frac{\bigl(a/\gamma,-q^2/abcx\gamma,-q/ax\gamma;q\bigr)_{\infty}}
{\bigl(q/\gamma^2,-q/bx\gamma,-q/acx;q\bigr)_{\infty}}
{}_3\phi_2\left( \begin{array}{c}
             q/a\gamma, -q/abx, -q/acx \\ 
             -q^2/abcx\gamma, -q/ax\gamma
            \end{array} ; q,a/\gamma \right).
\end{split}
\end{equation*}
The ${}_3\phi_2$ on the right hand side is the
${}_3\phi_2$ in the expression \eqref{PhiB} for $\Phi_{\gamma^{-1}}$.
So substituting this formula in the three term recurrence relation,
and simplifying the formulas using in particular the functional
relation
\begin{equation}\label{functheta}
\theta(q^kx)=
\begin{cases}
q^{-k(k-1)/2}(-x)^{-k}\theta(x),\qquad &k\in \ZZ_+,\\
q^{k(k+1)/2}(-x)^{-k}\theta(x), &k\in -\ZZ_+
\end{cases}
\end{equation}
for the Jacobi theta function, we arrive at the desired result
for restricted choices of $x$, $a$ and $\gamma$. 
The extension to all $x\in I_+$ is made
using the fact that the left hand side and the right hand
side of the $c$-function expansion are solutions 
of the eigenvalue equation
$L_+f=\mu(\gamma)f$ on $I_+$. Finally, the restrictions on $a$ and $\gamma$ 
can be removed by analytic continuation.

The proof of the connection coefficient formula for $\psi_{\gamma}$
follows by analytic continuation 
from the connection coefficient formula for $\phi_{\gamma}$
using \eqref{symmetryphi} and \eqref{symmetryPhi}.
\end{proof}
%%%%%%%%%%%%%%%%%%%%%%%%%%%%%%%%%%%%%%%%%%%%%%%%%%%%%%%%%%%%%%%%%%%%%

For $(a,b,c)\in V_z^{gen}$ and
$\gamma\in S_{reg}^+$ with $|\gamma|<a^{-1}$ we define
$\Phi^-_{\gamma}(x)$ for $x\in (\infty(-1),0]_q$ by
\begin{equation}\label{PhiB-}
\begin{split}
\Phi^-_{\gamma}(x)=&
\frac{\bigl(-q^2\gamma/abcx,-q\gamma/ax,a\gamma;q\bigr)_{\infty}}
{\bigl(-q/abx,-q/acx,q\gamma^2;q\bigr)_{\infty}}\\
&.\bigl(a\gamma\bigr)^{-k}{}_3\phi_2\left( \begin{array}{c}
            q\gamma/a, -q/abx, -q/acx \\ 
             -q^2\gamma/abcx, -q\gamma/ax 
            \end{array} ; q,a\gamma \right), \qquad x=-q^k.
\end{split}
\end{equation}
Observe that $\Phi^-_{\gamma}$ is obtained by taking $y=-1$ in
the definition of the eigenfunction $\Phi^{y}_{\gamma}$, see \eqref{PhiB}.

%%%%%%%%%%%%%%%%%%%%%%%%%%%%%%%%%%%%%%%%%%%%%%%%%%%%%%%%%%%%%%%%%%%%
\begin{lem}\label{collaps}
Let $(a,b,c)\in V_z^{gen}$ and $\gamma\in S_{reg}^+$ with $|\gamma|<a^{-1}$.
Then
\[\Phi^-_{\gamma}(x)=\frac{\bigl(q\gamma/a,q\gamma/b,
q\gamma/c;q\bigr)_{\infty}}{\bigl(q/ab,q/ac,q\gamma^2;q\bigr)_{\infty}}
\phi_{\gamma}(x),\qquad x\in I_-.
\]
\end{lem}
%%%%%%%%%%%%%%%%%%%%%%%%%%%%%%%%%%%%%%%%%%%%%%%%%%%%%%%%%%%%%%%%%%%%%%%%
\begin{proof}
Using the minimal solution \cite[(2.26)]{GIM} of the three term recurrence
relation \cite[(2.1)]{GIM}, it follows that 
$\Phi^-_{\gamma}$ is a solution of $Lf=\mu(\gamma)f$ on $(\infty(-1),0]_q$,
where $L$ is the second order $q$-difference operator defined by
\eqref{equiv2}. In particular, we have 
$A(-1)(\Phi^-_{\gamma}(-q)-\Phi^-_{\gamma}(-1))=
\mu(\gamma)\Phi^-_{\gamma}(-1)$ since $B(-1)=0$. It follows that
$\Phi^-_{\gamma}|_{I_-}\in V_{\mu(\gamma)}^-$. 
 
We have 
$\hbox{dim}(V_{\mu(\gamma)}^-)=1$ by Lemma \ref{generaleigenfunction}{\bf (i)}
and $0\not=\phi_{\gamma}|_{I_-}\in V_{\mu(\gamma)}^-$. It follows that
$\Phi^-_{\gamma}|_{I_-}=C_{\gamma}\phi_{\gamma}|_{I_-}$ for a unique constant
$C_{\gamma}$. The explicit expression for $C_{\gamma}$ can be found
using $\phi_{\gamma}(-1)=1$ and by evaluating $\Phi^-_{\gamma}(-1)$
using the $q$-Gauss sum \cite[(1.5.1)]{GR}. 
\end{proof}
%%%%%%%%%%%%%%%%%%%%%%%%%%%%%%%%%%%%%%%%%%%%%%%%%%%%%%%%%%%%%%%%%%%%%

Observe that the Wronskian $W(\phi_{\gamma}, \Phi^-_{\gamma})$ is identically
zero on $I_-$ by Lemma \ref{collaps}. On the other hand, the Wronskian
$W(\phi_{\gamma},\Phi_{\gamma})$ on $I_+$ is non-zero for generic $\gamma$
by Lemma \ref{dep} and Proposition \ref{conncoef}.
Suppose now that  for a given (generic) $\gamma\in\CC^*$ with $|\gamma|<a^{-1}$ 
we have an extension of 
$\Phi_{\gamma}\in V_{\mu(\gamma)}^+$ to a function $\Phi_{\gamma}\in
V_{\mu(\gamma)}$. Then the Wronskian $W(\phi_{\gamma},\Phi_{\gamma})$ 
is constant on $I$ by Lemma \ref{generaleigenfunction}{\bf (iii)}
and Lemma \ref{oknul}, so $\Phi_{\gamma}|_{I_-}$ is {\it not} 
a multiple of $\Phi^-_{\gamma}$ on $I_-$.

This shows that the extension of $\Phi_{\gamma}\in 
V_{\mu(\gamma)}^+$ to a function $\Phi_{\gamma}\in V_{\mu(\gamma)}$
can not be obtained by taking the explicit expression for $\Phi_{\gamma}$
and extending its definition in an obvious manner to the whole $q$-interval
$I$. We derive in the next section the extension of 
$\Phi_{\gamma}\in V_{\mu(\gamma)}^+$ 
to a function $\Phi_{\gamma}\in V_{\mu(\gamma)}$ for $(a,b,c)\in V_z^{gen}$
by writing $\Phi_{\gamma}$ explicitly as a linear 
combination of the two eigenfunctions $\phi_{\gamma}\in V_{\mu(\gamma)}$ and 
$\psi_{\gamma}\in V_{\mu(\gamma)}$.

%%%%%%%%%%%%%%%%%%%%%%%%%%%%%%%%%%%%%%%%%%%%%%%%%%%%%%%%%%%%%%%%%%
%%                                                              %%
%%           The extension of the asymptotic solution           %%
%%                                                              %%
%%%%%%%%%%%%%%%%%%%%%%%%%%%%%%%%%%%%%%%%%%%%%%%%%%%%%%%%%%%%%%%%%%

\section{The extension of the asymptotic solution}\label{section5}

In this section we first evaluate the Wronskian 
$W\bigl(\psi_{\gamma},\phi_{\gamma}\bigr)$ explicitly in product form.
{}From this evaluation we derive 
explicit conditions on the spectral parameter $\gamma$
for which we have $\hbox{dim}(V_{\mu(\gamma)})=2$. For these values of the
spectral parameter, there exists 
a unique extension of $\Phi_{\gamma}$ which lives in $V_{\mu(\gamma)}$
by Corollary \ref{extensionpossible}.
We will give this extension explicitly as a linear combination of $\phi_{\gamma}$
and $\psi_{\gamma}$.

Observe that if we compute the Wronskian by taking the limit to $0$,
\[W(\psi_{\gamma},\phi_{\gamma})(0)=
\frac{(1-q)^2}{bc}\bigl(\psi_{\gamma}'(0)\phi_{\gamma}(0)-
\psi_{\gamma}(0)\phi_{\gamma}'(0)\bigr),
\]
and by applying Lemma \ref{oknul}, we get an explicit expression 
as a linear combination of products of ${}_2\phi_2$'s, 
while evaluation of $W(\psi_{\gamma},\phi_{\gamma})(x)$ at $x=-1$
using \eqref{valuephi} and $\phi_{\gamma}(-1)=1$ 
gives an expression of the Wronskian
as a linear combination of two ${}_3\phi_2$'s. 
{}From both these expressions, the (non)vanishing properties of the Wronskian
$W(\psi_{\gamma},\phi_{\gamma})$ are hard to derive.
In the next proposition we evaluate the Wronskian
by substitution of the $c$-function expansions for $\phi_{\gamma}$
and $\psi_{\gamma}$.

Recall that the Wronskian $W\bigl(\psi_{\gamma},\phi_{\gamma}\bigr)
$ is constant on $I$ by Lemma
\ref{generaleigenfunction}{\bf (iii)} and Lemma \ref{oknul}. 
%%%%%%%%%%%%%%%%%%%%%%%%%%%%%%%%%%%%%%%%%%%%%%%%%%%%%%
\begin{prop}\label{Wpsiphiprop}
Let $\gamma\in\CC^*$ and $(a,b,c)\in V_z^{gen}$, then
\begin{equation}\label{Wpsiphi}
W\bigl(\psi_{\gamma},\phi_{\gamma}\bigr)=
(1-q)\frac{\bigl(a\gamma,a/\gamma;q\bigr)_{\infty}\theta(bc)}
{\bigl(ab,ac,qa/b,qa/c;q\bigr)_{\infty}}.
\end{equation}
\end{prop}
%%%%%%%%%%%%%%%%%%%%%%%%%%%%%%%%%%%%%%%%%%%%%%%%%%%%%%%%%
\begin{proof}
We first assume that $\gamma\in S_{reg}$.
Using Proposition \ref{conncoef} we then have
\begin{equation}\label{st1}
W\bigl(\psi_{\gamma},\phi_{\gamma}\bigr)=
\tilde{c}(\gamma)c(\gamma^{-1})
W\bigl(\Phi_{\gamma}, \Phi_{\gamma^{-1}}\bigr)(z)
+\tilde{c}(\gamma^{-1})c(\gamma)
W\bigl(\Phi_{\gamma^{-1}}, \Phi_{\gamma}\bigr)(z),
\end{equation}
where $c(\cdot)$ and $\tilde{c}(\cdot)$ are given by \eqref{cB}
and \eqref{tildecB}, respectively.
By Lemma \ref{dep}, the right hand
side of \eqref{st1} is equal to
\begin{equation}\label{st2}
W\bigl(\psi_{\gamma},\phi_{\gamma}\bigr)=-aK
\left(c(\gamma){\tilde{c}}(\gamma^{-1})-
c(\gamma^{-1}){\tilde{c}}(\gamma)\right)\bigl(\gamma-\gamma^{-1}\bigr),
\end{equation}
where $K$ is the positive constant defined by \eqref{Kconstant}.
It follows by direct computation that
\begin{equation*}
\begin{split}
c(\gamma){\tilde{c}}(\gamma^{-1})-
c(\gamma^{-1}){\tilde{c}}(\gamma)=
&\frac{1}{\bigl(ab,ac,qa/b,qa/c;q\bigr)_{\infty}\theta(-q/bcz,-1/z)}\\
.\frac{\bigl(a\gamma,a/\gamma;q\bigr)_{\infty}}
{\bigl(\gamma^2, \gamma^{-2};q\bigr)_{\infty}}
&\left\{\theta\bigl(q\gamma/b, c/\gamma, -q/abcz\gamma,
  -\gamma/az\bigr)\right.\\
&\left.\quad-\theta\bigl(q/b\gamma, c\gamma,-q\gamma/abcz,-1/az\gamma\bigr)
\right\}.
\end{split}
\end{equation*}
Now we can apply the $\theta$-product identity 
\begin{equation}\label{thetaidentity}
\theta(x\lambda,
x/\lambda,\mu\nu,\mu/\nu)-\theta(x\nu,x/\nu,\lambda\mu,\mu/\lambda)=
\frac{\mu}{\lambda}\theta(x\mu,x/\mu,\lambda\nu,\lambda/\nu),
\end{equation}
see \cite[Exercise 2.16]{GR},
with parameter values 
\[x=c\sqrt{\frac{q}{bc}}, 
\quad \lambda=\gamma\sqrt{\frac{q}{bc}},
\quad \mu=-\frac{1}{az}{\sqrt{\frac{q}{bc}}},
\quad \nu=\frac{1}{\gamma}\sqrt{\frac{q}{bc}},
\]
to obtain
\begin{equation}\label{st3}
\begin{split}
c(\gamma){\tilde{c}}(\gamma^{-1})-
c(\gamma^{-1}){\tilde{c}}(\gamma)=
&\frac{-1}{\bigl(ab,ac,qa/b,qa/c;q\bigr)_{\infty}\theta(-q/bcz,-1/z)}\\
&.\frac{1}{az\gamma}\frac{\bigl(a\gamma,a/\gamma;q\bigr)_{\infty}}
{\bigl(\gamma^2, \gamma^{-2};q\bigr)_{\infty}}
\theta(-q/abz,-acz,q/bc,\gamma^2).
\end{split}
\end{equation}
The proposition now follows for $\gamma\in S_{reg}$ by substitution of 
\eqref{st3} in \eqref{st2} and using $\theta(x)=\theta(q/x)$.
By continuity in $\gamma$, it follows that \eqref{Wpsiphi} holds for all
$\gamma\in \CC^*$. 
\end{proof}
%%%%%%%%%%%%%%%%%%%%%%%%%%%%%%%%%%%%%%%%%%%%%%%%%%%%%%%%%%%%%%%%%%%

In the following corollary we give a characterization 
of the eigenfunction $\phi_{\gamma}$ of $L$
for generic values of the spectral parameter $\gamma$.
In section \ref{PID} we show that the eigenfunction $\phi_{\gamma}$
plays the role of the spherical function for the big $q$-Jacobi 
function transform.
%%%%%%%%%%%%%%%%%%%%%%%%%%%%%%%%%%%%%%%%%%%%%%%%%%%%%%%%%%%%%%%%%%
\begin{cor}\label{psifoutind}
Let $\gamma\in \CC^*$ and $(a,b,c)\in V_z^{gen}$.\\ 
{\bf (i)} $\{\phi_{\gamma}, \psi_{\gamma}\}$ is a
linear basis of $V_{\mu(\gamma)}$  if 
\begin{equation}\label{con4}
\gamma\not\in \{aq^n,a^{-1}q^{-n}\}_{n\in \ZZ_+}.
\end{equation}
{\bf (ii)} Assume that \eqref{con4} is satisfied.
Let $f\in {\mathcal{F}}(I)$ be a function satisfying

{\bf --} $Lf=\mu(\gamma)f$ on $I$;

{\bf --} $f$ is continuously differentiable at the origin;

{\bf --} $f(-1)=1$.\\
Then $f=\phi_{\gamma}$.\\
{\bf (iii)} If \eqref{con4} is satisfied, then 
$\bigl(L\psi_{\gamma}\bigr)(-1)\not=\mu(\gamma)\psi_{\gamma}(-1)$.
\end{cor}
%%%%%%%%%%%%%%%%%%%%%%%%%%%%%%%%%%%%%%%%%%%%%%%%%%%%%%%%%%%%%%%%%%%
\begin{proof}
We fix $\gamma\in\CC^*$ satisfying \eqref{con4}.
Then the Wronskian $W(\psi_{\gamma},\phi_{\gamma})$ is non-zero since
$(a,b,c)\in V_z^{gen}$,
see Proposition \ref{Wpsiphiprop}. Hence $\{\phi_{\gamma},
\psi_{\gamma}\}$  
is a linear basis of $V_{\mu(\gamma)}$ by 
Proposition \ref{generaleigenfunctionprop} and Lemma \ref{oknul}.

For the proof of {\bf (ii)}, we observe first that $\phi_{\gamma}$
satisfies the three properties as stated in {\bf (ii)}, see section \ref{eigenf}.
If $f$ is another function satisfying the same three properties,
then $f|_{I_-}=\phi_{\gamma}|_{I_-}$ since 
$f(-1)=\phi_{\gamma}(-1)=1$ and
$\hbox{dim}\bigl(V_{\mu(\gamma)}^-\bigr)=1$ by Lemma \ref{generaleigenfunction}. 
Hence $f=\phi_{\gamma}$ on $I$ by Corollary \ref{extensionpossible}.

Finally, for the proof of {\bf (iii)}, observe that 
$\psi_{\gamma}|_{I_-}$ and  $\phi_{\gamma}|_{I_-}$ are linearly independent
by part {\bf (i)} of the corollary 
and by Proposition \ref{generaleigenfunctionprop}.
Since $\phi_{\gamma}|_{I_-}\in V_{\mu(\gamma)}^-$ by part {\bf (ii)} of 
the corollary and
$\hbox{dim}\bigl(V_{\mu(\gamma)}^-\bigr)=1$  
by Lemma \ref{generaleigenfunction}{\bf (i)}, it follows that
$\psi_{\gamma}|_{I_-}\in \hat{V}_{\mu(\gamma)}^-\setminus V_{\mu(\gamma)}^-$.
So $\bigl(L\psi_{\gamma}\bigr)(-1)\not=\mu(\gamma)\psi_{\gamma}(-1)$,
as desired.
\end{proof}
%%%%%%%%%%%%%%%%%%%%%%%%%%%%%%%%%%%%%%%%%%%%%%%%%%%%%%%%%%%%%%%%%%%%

Set $\gamma_n=aq^n$ ($n\in \ZZ_+$). We have the following 
description of the excluded set $\{\gamma_n^{\pm 1}\}_{n\in\ZZ_+}$
in Corollary \ref{psifoutind}. 
%%%%%%%%%%%%%%%%%%%%%%%%%%%%%%%%%%%%%%%%%%%%%%%%%%%%%%%%%%%%%%%%%%%%%
\begin{prop}\label{polcase}
Let $S_{pol}\subset \CC^*$ be the set of spectral parameters
$\gamma\in \CC^*$
for which the eigenvalue equation $\bigl(Lf\bigr)(x)=
\mu(\gamma)f(x)$ on $I$ has a solution $f\not=0$ 
which is polynomial in $x$. Then
$S_{pol}=\{\gamma_n^{\pm 1}\}_{n\in\ZZ_+}$. The non-zero polynomial
eigenfunction corresponding to the eigenvalue $\mu(\gamma_n^{\pm 1})$ 
is the big $q$-Jacobi
polynomial of degree $n$ and is explicitly given by 
\begin{equation}\label{connpol}
\begin{split}
\phi_{\gamma_n}(x)&=\frac{\bigl(qa/b,qa/c;q\bigr)_n}
{\bigl(ab,ac;q\bigr)_n}\left(\frac{bc}{q}\right)^n\psi_{\gamma_n}(x)\\
&=\frac{\bigl(qa/c;q\bigr)_n}{\bigl(ac;q\bigr)_n}
\left(\frac{-c}{aq^{(n+1)/2}}\right)^n
{}_3\phi_2\left( \begin{array}{c}
            q^{-n}, -abx, q^{n}a^2 \\ 
             ab, qa/c 
            \end{array} ; q,q \right).
\end{split}
\end{equation}
\end{prop}
%%%%%%%%%%%%%%%%%%%%%%%%%%%%%%%%%%%%%%%%%%%%%%%%%%%%%%%%%%%%%%%%%%%
\begin{proof}
This is well known, see for instance \cite{AA} and \cite[section
7.3]{GR}. The connection between $\phi_{\gamma_n}$ and 
$\psi_{\gamma_n}$ with the ${}_3\phi_2$ in the last equality of
\eqref{connpol} follows from \cite[(3.2.5)]{GR}.
\end{proof}  
%%%%%%%%%%%%%%%%%%%%%%%%%%%%%%%%%%%%%%%%%%%%%%%%%%%%%%%%%%%%%%
\begin{rem}
Proposition \ref{Wpsiphi} and Proposition \ref{polcase} 
are essential ingredients for a functional-analytic derivation
of the orthogonality relations and the quadratic norm evaluations
for the big $q$-Jacobi polynomials, see section \ref{polsection}.
\end{rem}
%%%%%%%%%%%%%%%%%%%%%%%%%%%%%%%%%%%%%%%%%%%%%%%%%%%%%%%%%%%%%%%%%%%

We write $S_{pol}^{\pm}=\{\gamma_n^{\pm 1}\}_{n\in\ZZ_+}$, where
$\gamma_n=aq^n$. 
Observe that $S_{pol}^+\subset S_{pol}\subset S_{reg}\subset S_{reg}^+$ 
when $(a,b,c)\in V_z^{gen}$.

%%%%%%%%%%%%%%%%%%%%%%%%%%%%%%%%%%%%%%%%%%%%%%%%%%%%%%%%%%%%%%%%%%%
\begin{prop}\label{extensionPhi}
Let $\gamma\in S_{reg}^+\setminus S_{pol}^+$ and $(a,b,c)\in V_z^{gen}$.
Then 
\begin{equation}\label{connpsiphi}
\Phi_{\gamma}(x)=K(\gamma)\phi_{\gamma}(x)+
\tilde{K}(\gamma)\psi_{\gamma}(x),\qquad x\in I_+,
\end{equation}
with $K(\gamma)=K(\gamma;a,b,c;z)$ given by
\[
K(\gamma)=\frac{\bigl(ab,ac,q\gamma/b,q\gamma/c;q\bigr)_{\infty}}
{\bigl(q\gamma^2,a/\gamma;q\bigr)_{\infty}}
\frac{\theta(-bcz,-az/\gamma)}{\theta(bc,-abz,-acz)},
\]
and $\tilde{K}(\gamma)=\tilde{K}(\gamma;a,b,c;z)=
K(\gamma;a,q/b,q/c;bcz/q)$.
\end{prop}
%%%%%%%%%%%%%%%%%%%%%%%%%%%%%%%%%%%%%%%%%%%%%%%%%%%%%%%%%%%%%%%%%%%%%%
\begin{proof}
We first prove \eqref{connpsiphi} for 
$\gamma\in S_{reg}\setminus S_{pol}$.
By Corollary \ref{psifoutind}{\bf (i)} and Proposition
\ref{generaleigenfunctionprop}, there exist unique 
$K(\gamma), \tilde{K}(\gamma)\in \CC$ such that \eqref{connpsiphi}
holds for all $x\in I_+$. 
These coefficients can be expressed in terms of Wronskians by
\begin{equation}\label{exprW}
K(\gamma)=
\frac{W(\psi_{\gamma},\Phi_{\gamma})(z)}{W(\psi_{\gamma},\phi_{\gamma})},
\qquad 
\tilde{K}(\gamma)=
\frac{W(\Phi_{\gamma},\phi_{\gamma})(z)}{W(\psi_{\gamma},\phi_{\gamma})}.
\end{equation}
By Lemma \ref{dep} and Proposition \ref{conncoef} 
we have
\begin{equation}\label{laat}
\begin{split}
W(\psi_{\gamma},\Phi_{\gamma})(z)&=
aK\tilde{c}(\gamma^{-1})\bigl(\gamma^{-1}-\gamma\bigr),\\
W(\Phi_{\gamma},\phi_{\gamma})(z)&=
aKc(\gamma^{-1})\bigl(\gamma-\gamma^{-1}\bigr),
\end{split}
\end{equation}
where $K$ is the positive constant defined by \eqref{Kconstant}.
Formula \eqref{connpsiphi} now follows by substituting \eqref{laat}, the
explicit formula \eqref{Wpsiphi} 
for the Wronskian $W(\psi_{\gamma},\phi_{\gamma})$, 
and the explicit expressions \eqref{cB} and 
\eqref{tildecB} for $c(\gamma)$ and $\tilde{c}(\gamma)$ in \eqref{exprW}, 
and using the theta function identities 
$\theta(x)=\theta(q/x)$ and \eqref{functheta}. 

It follows now by continuity in $\gamma$ 
that \eqref{connpsiphi} is valid for 
$\gamma\in S_{reg}^+\setminus S_{pol}^+$.
\end{proof}
%%%%%%%%%%%%%%%%%%%%%%%%%%%%%%%%%%%%%%%%%%%%%%%%%%%%%%%%%%%%%%%%%%

%%%%%%%%%%%%%%%%%%%%%%%%%%%%%%%%%%%%%%%%%%%%%%%%%%%%%%%%%%%%%%%%%%%%%%
\begin{rem} 
Fix $\gamma\in S_{reg}\setminus S_{pol}$ and suppose that
$(a,b,c)\in V_z^{gen}$.
The proof of Proposition \ref{Wpsiphiprop} implies the matrix equation
\begin{equation}\label{matrixequation}
\begin{pmatrix}
c(\gamma) & c(\gamma^{-1})\\
\tilde{c}(\gamma) & \tilde{c}(\gamma^{-1})
\end{pmatrix}
=\begin{pmatrix}
K(\gamma)  & \tilde{K}(\gamma)\\
K(\gamma^{-1}) & \tilde{K}(\gamma^{-1})
\end{pmatrix}^{-1}.
\end{equation}
Indeed, \eqref{matrixequation} follows easily from the explicit formula
\eqref{st3} for the determinant of the matrix on the left hand side of 
\eqref{matrixequation}. Let $x\in I_+$, then the matrix equation 
\eqref{matrixequation} implies that the two connection coefficient formulas
\[
\phi_{\gamma}(x)=c(\gamma)\Phi_{\gamma}(x)+
c(\gamma^{-1})\Phi_{\gamma^{-1}}(x),\quad
\psi_{\gamma}(x)=\tilde{c}(\gamma)\Phi_{\gamma}(x)+
\tilde{c}(\gamma^{-1})\Phi_{\gamma^{-1}}(x)
\]
are equivalent to the two connection coefficient formulas
$\Phi_{\gamma^{\pm 1}}(x)=K(\gamma^{\pm 1})\phi_{\gamma}(x)+
\tilde{K}(\gamma^{\pm 1})\psi_{\gamma}(x)$.
\end{rem}
%%%%%%%%%%%%%%%%%%%%%%%%%%%%%%%%%%%%%%%%%%%%%%%%%%%%%%%%%%%%%%%%%

%%%%%%%%%%%%%%%%%%%%%%%%%%%%%%%%%%%%%%%%%%%%%%%%%%%%%%%%%%%%%%%%%%%
\begin{cor}\label{imp}
Let $\gamma\in S_{reg}^+\setminus S_{pol}$ and $(a,b,c)\in V_z^{gen}$.
The unique extension of $\Phi_{\gamma}\in V_{\mu(\gamma)}^+$
to a function $\Phi_{\gamma}\in V_{\mu(\gamma)}$ 
is given by
\begin{equation}\label{extensionPhi!}
\Phi_{\gamma}(x)=K(\gamma)\phi_{\gamma}(x)+\tilde{K}(\gamma)\psi_{\gamma}(x),
\qquad x\in I,
\end{equation}
with  $K(\gamma)$, $\tilde{K}(\gamma)$ as defined in Proposition 
\ref{Wpsiphiprop}.
\end{cor}
%%%%%%%%%%%%%%%%%%%%%%%%%%%%%%%%%%%%%%%%%%%%%%%%%%%%%%%%%%%%%%%%%%%
\begin{proof}
Immediate from  Corollary \ref{extensionpossible},
Corollary \ref{psifoutind}{\bf (i)} and Proposition \ref{extensionPhi}. 
\end{proof}
%%%%%%%%%%%%%%%%%%%%%%%%%%%%%%%%%%%%%%%%%%%%%%%%%%%%%%%%%%%%%%%%%%%%

In section \ref{A} we have seen that $\Phi_{\gamma}(x)$ depends analytically
on $\gamma\in S_{reg}^+$ for all $x\in I_+$. We have the following analogous
result for the extension \eqref{extensionPhi!}
of $\Phi_{\gamma}\in V_{\mu(\gamma)}^+$.
 
%%%%%%%%%%%%%%%%%%%%%%%%%%%%%%%%%%%%%%%%%%%%%%%%%%%%%%%%%%%%%%%%%%%
\begin{lem}\label{reg}
Let $\gamma\in S_{pol}$ and $(a,b,c)\in V_z^{gen}$. 
There exists a unique extension of 
$\Phi_{\gamma}\in V_{\mu(\gamma)}^+$ to a function
$\Phi_{\gamma}\in V_{\mu(\gamma)}$ such that
$\tilde{\gamma}\mapsto \Phi_{\tilde{\gamma}}(x)$ is analytic at 
$\tilde{\gamma}=\gamma$ for all $x\in I$.

If $\gamma=\gamma_n^{-1}\in S_{pol}^-$, then the extension of 
$\Phi_{\gamma}$
is given by \eqref{extensionPhi!}. If $\gamma=\gamma_n\in 
S_{pol}^+$, then the extension of $\Phi_{\gamma_n}$ is given by
\begin{equation}\label{anext}
\Phi_{\gamma_n}(x)=
\left(\underset{\gamma=\gamma_n}{\hbox{Res}}K(\gamma)\right)
\frac{\partial\phi_{\gamma}}{\partial\gamma}(x)|_{\gamma=\gamma_n}
+M_n\phi_{\gamma_n}(x)+ \left(\underset{\gamma=\gamma_n}{\hbox{Res}}
\tilde{K}(\gamma)\right)
\frac{\partial\psi_{\gamma}}{\partial\gamma}(x)|_{\gamma=\gamma_n}
\end{equation}
for $x\in I$, with $M_n$ given by the existing limit 
$M_n=\lim_{\gamma\to\gamma_n}M_n(\gamma)$, where
\begin{equation}\label{Mngamma}
M_n(\gamma)=\frac{\bigl(qa/b,qa/c;q\bigr)_n}{\bigl(ab,ac;q\bigr)_n}
\left(\frac{bc}{q}\right)^nK(\gamma)+\tilde{K}(\gamma).
\end{equation}
\end{lem}
%%%%%%%%%%%%%%%%%%%%%%%%%%%%%%%%%%%%%%%%%%%%%%%%%%%%%%%%%%%%%%%%%%%%
\begin{proof}
The proof for $\gamma\in S_{pol}^-$ is trivial
since $K(\tilde{\gamma})$, $\tilde{K}(\tilde{\gamma})$, 
$\phi_{\tilde{\gamma}}(x)$ and $\psi_{\tilde{\gamma}}(x)$ are regular
at $\tilde{\gamma}=\gamma$.

For $\gamma_n\in S_{pol}^+$ ($n\in\ZZ_+$), observe that $K(\gamma)$ and
$\tilde{K}(\gamma)$ have simple poles at $\gamma=\gamma_n$ and that
$\phi_{\gamma}(x)$ and $\psi_{\gamma}(x)$ are regular at 
$\gamma=\gamma_n$. It follows from \eqref{extensionPhi!} and the first 
equality of \eqref{connpol} that the singularity of $\Phi_{\gamma}(x)$
at $\gamma=\gamma_n$ is removable for $x\in I$ if the (at most simple)
singularity of $M_n(\gamma)$ at $\gamma=\gamma_n$ is removable. This can be
checked by direct computions using the theta function identities
$\theta(x)=\theta(q/x)$ and \eqref{functheta}.

It follows that $\lim_{\gamma\to\gamma_n}M_n(\gamma)$ exists
and that \eqref{anext} holds. The extension \eqref{anext} of
$\Phi_{\gamma_n}$ lies in $V_{\mu(\gamma_n)}$
because $\phi_{\gamma_n}$ and the derivatives
of $\phi_{\gamma}$ and $\psi_{\gamma}$ at $\gamma=\gamma_n$ are
continuously differentiable at the origin, cf. Lemma \ref{oknul}.
\end{proof}
%%%%%%%%%%%%%%%%%%%%%%%%%%%%%%%%%%%%%%%%%%%%%%%%%%%%%%%%%%%%%%%%%%%%%%

For future reference, we collect here the main results concerning
the asymptotic solution of the eigenvalue equation $Lf=\mu(\gamma)f$.

%%%%%%%%%%%%%%%%%%%%%%%%%%%%%%%%%%%%%%%%%%%%%%%%%%%%%%%%%%%%%%%%%%%%
\begin{thm}\label{samenvatting}
Let $x\in I$ and $(a,b,c)\in V_z^{gen}$.

{\bf (i)} $\Phi_{\gamma}\in V_{\mu(\gamma)}$ for all $\gamma\in
S_{reg}^+$.

{\bf (ii)} $\Phi_{\gamma}(x)$ is analytic at $\gamma\in S_{reg}^+$.

{\bf (iii)} $\Phi_{\gamma}\in {\mathcal{D}}$ for $\gamma\in \CC^*$
with $|\gamma|<1$.

{\bf (iv)} $\phi_{\gamma}(x)=c(\gamma)\Phi_{\gamma}(x)+
c(\gamma^{-1})\Phi_{\gamma^{-1}}(x)$ for all $\gamma\in S_{reg}$.

{\bf (v)} Let $\gamma\in S_{reg}^+$. 
The Wronskian $W(\gamma)=W\bigl(\Phi_{\gamma},\phi_{\gamma}\bigr)$
is constant on $I$, and 
$W(\gamma)=aKc(\gamma^{-1})\bigl(\gamma-\gamma^{-1}\bigr)$, where $K$
is the positive constant defined by \eqref{Kconstant}.
\end{thm}
%%%%%%%%%%%%%%%%%%%%%%%%%%%%%%%%%%%%%%%%%%%%%%%%%%%%%%%%%%%%%%%%%%%%%%%%
\begin{proof}
{\bf (i)} and {\bf (ii)} follow from Corollary \ref{imp} and Lemma \ref{reg}.\\
{\bf (iii)} Let $\gamma\in \CC^*$ with $|\gamma|<1$, then 
$\gamma\in S_{reg}^+$, hence $\Phi_{\gamma}$ is well defined. 
Then $\Phi_{\gamma}\in {\mathcal{D}}$ follows from 
part {\bf (i)} of the theorem and from Lemma \ref{l2}.\\
{\bf (iv)} We first prove the connection coefficient formula for
$\gamma\in S_{reg}\setminus S_{pol}$. Then the connection 
coefficient formula is
valid for all $x\in I_+$, see Proposition \ref{conncoef}. Since 
$\phi_{\gamma},\Phi_{\gamma^{\pm 1}}\in V_{\mu(\gamma)}$, it follows
by the uniqueness property of extensions of eigenfunctions, see 
Corollary \ref{extensionpossible} and Corollary \ref{psifoutind}{\bf (i)},
that the connection coefficient formula is valid for all $x\in I$.
The connection coefficient formula holds then for all $\gamma\in S_{reg}$
by continuity.\\
{\bf (v)} This follows from Lemma \ref{generaleigenfunction}{\bf (iii)},
part {\bf (iv)} of the theorem, and Lemma \ref{dep}.
\end{proof}
%%%%%%%%%%%%%%%%%%%%%%%%%%%%%%%%%%%%%%%%%%%%%%%%%%%%%%%%%%%%%%%%%%%%%%%
\begin{rem}
Observe that the statements of Theorem \ref{samenvatting}{\bf (i)} 
and {\bf (iii)} are not in contradiction with the self-adjointness of
$(L,{\mathcal{D}})$, see Proposition \ref{s}.
Indeed, let $\gamma\in \CC\setminus \RR$
with $|\gamma|<1$, then $\gamma\in S_{reg}^+$ and
$\mu(\gamma)\in \CC\setminus\RR$. By 
Theorem \ref{samenvatting}{\bf (i)} and {\bf (iii)} we have that 
$\Phi_{\gamma}\in {\mathcal{D}}$ and 
$\bigl(L\Phi_{\gamma}\bigr)(x)=\mu(\gamma)\Phi_{\gamma}(x)$ 
for all $x\in (-1,\infty(z))_q$, 
but the eigenvalue equation does not hold in the end-point $x=-1$.
In fact, the self-adjointness
of $(L,{\mathcal{D}})$ forces that $\bigl(L\Phi_{\gamma}\bigr)(-1)\not=
\mu(\gamma)\Phi_{\gamma}(-1)$. Another proof of this inequality 
can be given using the non-vanishing of the 
Wronskian $W(\phi_{\gamma},\Phi_{\gamma})$,
cf. the proof of Corollary \ref{psifoutind}{\bf (iii)}.
\end{rem}
%%%%%%%%%%%%%%%%%%%%%%%%%%%%%%%%%%%%%%%%%%%%%%%%%%%%%%%%%%%%%%%%%%%

%%%%%%%%%%%%%%%%%%%%%%%%%%%%%%%%%%%%%%%%%%%%%%%%%%%%%%%%%%%%%
%%                                                         %%
%%                 The Green function                      %%
%%                                                         %%
%%%%%%%%%%%%%%%%%%%%%%%%%%%%%%%%%%%%%%%%%%%%%%%%%%%%%%%%%%%%%
\section{The Green function}

In this section we define the Green function of the self-adjoint
operator $(L,{\mathcal{D}})$. With the proper extensions of the
asymptotic expansion $\Phi_{\gamma}$ now at hand, see 
Theorem \ref{samenvatting}, the construction of the Green function is
a straightforward extension of the construction given by Kakehi 
\cite{K} and Kakehi, Masuda and Ueno \cite{KMU} 
for the little $q$-Jacobi function transform. We use some standard terminology
and results for unbounded self-adjoint operators and their spectral measures
for which we refer to Dunford and Schwartz \cite[Chapter XII]{DS} and
Rudin \cite[Chapter 13]{R}.

Let $(a,b,c)\in V_z^{gen}$ and $\mu\in\CC\setminus\RR$. 
Let $\gamma\in\CC^*$ be the
unique non-zero complex number such that $\mu=\mu(\gamma)$ and $|\gamma|<1$.
Observe that $\gamma\in\CC\setminus\RR$, hence $\gamma\in S_{reg}$
and $W(\gamma)\not=0$, see Theorem \ref{samenvatting}{\bf (v)}.
We define the Green kernel $K_{\gamma}(x,y)$ for $x,y\in I$ by
\begin{equation}\label{K}
K_{\gamma}(x,y)=
\begin{cases}
W(\gamma)^{-1}\Phi_{\gamma}(x)\phi_{\gamma}(y), 
\qquad &y\leq x,\\
W(\gamma)^{-1}\phi_{\gamma}(x)\Phi_{\gamma}(y),
\qquad &y\geq x.
\end{cases}
\end{equation}
Observe that $K_{\gamma}(x,\cdot), K_{\gamma}(\cdot,x)\in {\mathcal{D}}$ 
for all $x\in I$ in view of Lemma \ref{oknul} and 
Theorem \ref{samenvatting}{\bf (i)} and {\bf (iii)},
so we have a well defined linear map ${\mathcal{H}}\to {\mathcal{F}}(I)$
mapping $f\in {\mathcal{H}}$ to 
\begin{equation}\label{Greenfunction}
G_f(x,\gamma)=\langle f, {\overline{K_{\gamma}(x,\cdot)}}\rangle, \qquad x\in I.
\end{equation}
Written out explicitly, we arrive at the formula
\begin{equation}\label{expl}
G_f(x,\gamma)=W(\gamma)^{-1}
\left(\Phi_{\gamma}(x)\int_{-1}^xf(y)\phi_{\gamma}(y)\frac{d_qy}{p(y)}
+\phi_{\gamma}(x)\int_{x}^{\infty(z)}f(y)\Phi_{\gamma}(y)
\frac{d_qy}{p(y)}\right).
\end{equation}
By the self-adjointness of $(L,{\mathcal{D}})$ we have that
the resolvent $\bigl(L-\mu.\hbox{Id}\bigr)^{-1}$ is a one to one, 
continuous map from ${\mathcal{H}}$ onto 
${\mathcal{D}}$ for all $\mu\in\CC\setminus \RR$.

%%%%%%%%%%%%%%%%%%%%%%%%%%%%%%%%%%%%%%%%%%%%%%%%%%%%%%%%%%%%%%%%%%%
\begin{prop}\label{G}
Let $(a,b,c)\in V_z^{gen}$, $f\in {\mathcal{H}}$ and 
$\mu=\mu(\gamma)\in \CC\setminus\RR$ with $\gamma\in \CC^*$ and
$|\gamma|<1$. Then 
$G_f(\cdot,\gamma)=\bigl(L-\mu.\hbox{Id}\bigr)^{-1}f$. 
In particular, $G_f(\cdot,\gamma)\in {\mathcal{D}}$.
\end{prop}
%%%%%%%%%%%%%%%%%%%%%%%%%%%%%%%%%%%%%%%%%%%%%%%%%%%%%%%%%%%%%%%%%%%%%
\begin{proof}
We first prove that 
\begin{equation}\label{check} 
\bigl(\bigl(L-\mu(\gamma)\bigr)G_f(\cdot,\gamma)\bigr)(x)=f(x), 
\qquad 
\forall f\in {\mathcal{H}}.
\end{equation}
For the proof we need to consider the
two cases $x\in (-1,\infty(z))_q$ and $x=-1$ seperately.\\
{\it Case 1:} $x\in (-1,\infty(z))_q$.

Observe that the product rule \eqref{productrule}
for the $q$-derivative $D_q$
and Lemma \ref{self-adjointform} 
imply the following product rule for $L$:
\begin{equation*}
\begin{split}
\bigl(L(fg)\bigr)(x)&=\bigl(Lf\bigr)(x)g(x)+f(x)\bigl(Lg\bigr)(x)
+qp(x)r(x)\bigl(D_qf\bigr)(x)\bigl(D_qg\bigr)(x)\\
&\qquad\qquad\qquad
+p(x)r(q^{-1}x)\bigl(D_qf\bigr)(q^{-1}x)\bigl(D_qg\bigr)(q^{-1}x).
\end{split}
\end{equation*}
Combined with the easily verified formulas
\begin{equation*}
D_q\left(x\mapsto \int_{-1}^xf(y)d_qy\right)=f(x),\qquad
D_q\left(x\mapsto\int_{x}^{\infty(z)}f(y)d_qy\right)=-f(x)
\end{equation*} 
and the definition \eqref{expl} of $G_f(x,\gamma)$, we obtain
\begin{equation*}
\begin{split}
&W(\gamma)\left(\bigl(L-
\mu(\gamma)\bigr)G_f(\cdot,\gamma)\right)(x)=\\
&\quad\qquad\Phi_{\gamma}(x)p(x)\left(D_q(p^{-1}
rf\phi_{\gamma})\right)(q^{-1}x)
-\phi_{\gamma}(x)p(x)\left(D_q(p^{-1}rf\Phi_{\gamma})\right)(q^{-1}x)\\
&\quad\qquad+p(x)\left(p^{-1}rf\bigl(D_q\Phi_{\gamma}\bigr)
\phi_{\gamma}\right)(q^{-1}x)
-p(x)\left(p^{-1}rf\bigl(D_q\phi_{\gamma}\bigr)
\Phi_{\gamma}\right)(q^{-1}x)\\
&\quad\qquad+q\left(rf\bigl(D_q\Phi_{\gamma}\bigr)\phi_{\gamma}\right)(x)
-q\left(rf\bigl(D_q\phi_{\gamma}\bigr)\Phi_{\gamma}\right)(x)
\end{split}
\end{equation*}
for $f\in {\mathcal{H}}$.
By writing out the $D_q$-terms in this formula, 
we see that the right hand side reduces to 
$W(\Phi_{\gamma},\phi_{\gamma})(x)f(x)$, which in turn is equal to
$W(\gamma)f(x)$ by Theorem
\ref{samenvatting}{\bf (v)}. This completes
the proof of \eqref{check} for $x\in (-1,\infty(z))_q$.\\
{\it Case 2:} $x=-1$.

Observe that the Green function in the end-point $-1$ reduces to
\[G_f(-1,\gamma)=\int_{-1}^{\infty(z)}f(y)\Phi_{\gamma}(y)\frac{d_qy}{p(y)}
\]
since $\phi_{\gamma}(-1)=1$. Using the product rule \eqref{productrule}, 
we then have that
\begin{equation*}
\begin{split}
W(\gamma)&\bigl(D_qG_f(\cdot,\gamma)\bigr)(-1)=
\Phi_{\gamma}(-q)\phi_{\gamma}(-1)\frac{f(-1)}{p(-1)}\\
&\qquad\qquad-\phi_{\gamma}(-q)\Phi_{\gamma}(-1)\frac{f(-1)}{p(-1)}
+\bigl(D_q\phi_{\gamma}\bigr)(-1)G_f(-1,\gamma).
\end{split}
\end{equation*}
Combined with the expression of $L$ in 
the end-point $x=-1$ in terms
of the $q$-derivative and the functions $p(\cdot)$ and $r(\cdot)$, 
see Lemma \ref{self-adjointform}, and using
the fact that $(L\phi_{\gamma})(-1)=\mu(\gamma)\phi_{\gamma}(-1)$, we
obtain
\[
W(\gamma)\left(\bigl(L-\mu(\gamma)\bigr)G_f(\cdot,\gamma)\right)(-1)=
W(\Phi_{\gamma},\phi_{\gamma})(-1)f(-1).
\] 
Hence \eqref{check} is valid for the end-point $x=-1$ since 
$W(\Phi_{\gamma},\phi_{\gamma})(-1)=W(\gamma)$ by 
Theorem \ref{samenvatting}{\bf (v)}.

It remains to show how \eqref{check} leads to a proof of the proposition.
Let ${\mathcal{D}}_{fin}$ 
be the set of functions $f:I\to\CC$ with finite support. Then 
${\mathcal{D}}_{fin}\subset {\mathcal{D}}\subset {\mathcal{H}}$
as dense subspaces. Let $f\in {\mathcal{D}}_{fin}$, then 
$G_f(\cdot,\gamma)$ is continuously differentiable at the origin since 
$\phi_{\gamma}, \Phi_{\gamma}\in V_{\mu(\gamma)}$. Furthermore,
$G_f(x,\gamma)$ is a constant multiple of $\Phi_{\gamma}(x)$ for $x\gg 0$, 
hence 
$G_f(\cdot,\gamma)\in {\mathcal{H}}$.
By \eqref{check} it follows that 
$LG_f(\cdot,\gamma)\in {\mathcal{H}}$,
hence $G_f(\cdot,\gamma)\in {\mathcal{D}}$.
Combined with \eqref{check}, this proves that
\begin{equation}\label{check2}
\bigl((L-\mu(\gamma).\hbox{Id})^{-1}f\bigr)(x)=G_f(x,\gamma)=
\langle f,{\overline{K_{\gamma}(x,\cdot)}}\rangle,\qquad 
\forall f\in {\mathcal{D}}_{fin}
\end{equation}
for all $x\in I$. By continuity, \eqref{check2} is valid for all 
$f\in {\mathcal{H}}$. This completes the proof of the proposition.
\end{proof}
%%%%%%%%%%%%%%%%%%%%%%%%%%%%%%%%%%%%%%%%%%%%%%%%%%%%%%%%%%%%%%%%%%
\begin{rem}\label{phibetter2}
Proposition \ref{G} is in general not   
valid when $\phi_{\gamma}$ is replaced by $\psi_{\gamma}$ and $W(\gamma)$
is replaced by the constant value of $W(\Phi_{\gamma},\psi_{\gamma})$
on $I$ in the definition of the Green function
$G_f(\cdot,\gamma)$. Indeed, for the proof of 
$\left(\bigl(L-\mu(\gamma)\bigr)G_f(\cdot,\gamma)\right)(-1)=f(-1)$ in the
previous proposition, we use the fact that
$\phi_{\gamma}$ is a solution of $\bigl(Lf\bigr)(x)=\mu(\gamma)f(x)$
in the end-point $x=-1$. This property fails in general to be true for
$\psi_{\gamma}$, see Corollary \ref{psifoutind}{\bf (iii)}.
\end{rem}
%%%%%%%%%%%%%%%%%%%%%%%%%%%%%%%%%%%%%%%%%%%%%%%%%%%%%%%%%%%%%%%%%%

Proposition \ref{G} plays a crucial role in determining the explicit
form of the resolution of the identity $L=\int_{\RR}tdE(t)$
for the self-adjoint operator $(L,{\mathcal{D}})$ on ${\mathcal{H}}$ 
since the spectral measure $E$ is related to the resolvent of $L$
by
\begin{equation}\label{inversion}
\begin{split}
&\langle E\bigl((\mu_1,\mu_2)\bigr)f,g\rangle\\
&\,\,=\lim_{\delta\downarrow 0}\lim_{\epsilon\downarrow 0}
\frac{1}{2\pi i}\int_{\mu_1+\delta}^{\mu_2-\delta}
\left(\langle\bigl(L-(\mu+i\epsilon)\bigr)^{-1}f,g\rangle-
\langle\bigl(L-(\mu-i\epsilon)\bigr)^{-1}f,g\rangle\right)d\mu,
\end{split}
\end{equation}
where $\mu_1<\mu_2$ and $f\in {\mathcal{D}}$, $g\in {\mathcal{H}}$,
see \cite[Theorem XII.2.10]{DS}.
In the following two sections we use Proposition \ref{G}
and \eqref{inversion} to give an explicit description of the  
continuous and discrete contributions to the spectral resolution $E$.

%%%%%%%%%%%%%%%%%%%%%%%%%%%%%%%%%%%%%%%%%%%%%%%%%%%%%%%%%%%%%%%%
%%                                                            %%
%%           The continuous spectrum                          %%
%%                                                            %%
%%%%%%%%%%%%%%%%%%%%%%%%%%%%%%%%%%%%%%%%%%%%%%%%%%%%%%%%%%%%%%%%

\section{The continuous spectrum}\label{contmass}

We start this section by proving that the closed
interval $[-(1+a)^2,-(1-a)^2]$ is contained in the continuous
spectrum $\sigma_c(L)$ of $(L,{\mathcal{D}})$. 
In section \ref{pointmass} we will see in fact that this interval is
exactly equal to $\sigma_c(L)$ for $(a,b,c)\in V_z^{gen}$. 
Furthermore, we compute the spectral projection
$P_c:=E\bigl([-(1+a)^2,-(1-a)^2]\bigr)$ explicitly, and
give the Plancherel formula and inversion formula for the continuous
part of the big $q$-Jacobi function transform.

For $n\in \NN$ and $x\in I$ we set 
$\varphi_\gamma^{(n)}(x)=\phi_{\gamma}^{(n)}(x)/\|\phi_{\gamma}^{(n)}\|$,
where $\phi_{\gamma}^{(n)}\in {\mathcal{H}}$ is defined by
\begin{equation*}
\phi_{\gamma}^{(n)}(x)=
\begin{cases}
\phi_{\gamma}(x)\quad  &\hbox{ if } \,\,\, x\in I\setminus 
[zq^{-n-1},\infty(z))_q,\\
0 &\hbox{ if } \,\,\, x\in [zq^{-n-1},\infty(z))_q.
\end{cases}
\end{equation*}

%%%%%%%%%%%%%%%%%%%%%%%%%%%%%%%%%%%%%%%%%%%%%%%%%%%%%%%%%%%%%%%%%
\begin{lem}\label{contspectrum}
Let $\mu\in [-(1+a)^2,-(1-a)^2]$ and $\theta\in [0,\pi]$ such that
$\mu=\mu(e^{i\theta})$. Then 
$\bigl(\varphi_{e^{i\theta}}^{(n)}\bigr)_{n\in\NN}$
is a generalized eigenfunction of $\bigl(L,{\mathcal{D}}\bigr)$ 
with generalized eigenvalue $\mu$. 
Furthermore, the continuous spectrum $\sigma_c(L)$ of the 
self-adjoint operator
$\bigl(L, {\mathcal{D}}\bigr)$ 
contains the interval $[-(1+a)^2, -(1-a)^2]$.
\end{lem}
%%%%%%%%%%%%%%%%%%%%%%%%%%%%%%%%%%%%%%%%%%%%%%%%%%%%%%%%%%%%%%%%%
\begin{proof}
Recall that $\bigl(\varphi_{e^{i\theta}}^{(n)}\bigr)_{n\in\NN}$
is a generalized eigenfunction of $(L,{\mathcal{D}})$ with generalized
eigenvalue $\mu(e^{i\theta})$ if 
$\varphi_{e^{i\theta}}^{(n)}\in {\mathcal{D}}$, 
$\|\varphi_{e^{i\theta}}^{(n)}\|=1$ and 
$\lim_{n\to\infty}\bigl(L\varphi_{e^{i\theta}}^{(n)}-
\mu(e^{i\theta})\varphi_{e^{i\theta}}^{(n)}\bigr)=0$ in ${\mathcal{H}}$. 
This can be checked for all $\theta\in [0,\pi]$ by
an elementary computation
using Proposition \ref{psifoutind}{\bf (ii)}.

Let $\mu\in [-(1+a)^2,-(1-a)^2]$. Then $\mu$ is part of the spectrum 
$\sigma(L)$ of $(L,{\mathcal{D}})$ since there exists a generalized 
eigenfunction of $(L,{\mathcal{D}})$ with generalized eigenvalue $\mu$.
Then $\mu$ is in the continuous
spectrum $\sigma_c(L)$ 
of $(L,{\mathcal{D}})$ or in the point spectrum $\sigma_p(L)$ of
$(L,{\mathcal{D}})$ since a self-adjoint operator does not have
residual spectrum, see \cite[Theorem 13.27]{R}.

It remains to show that $\mu\not\in \sigma_p(L)$. 
We have to distinguish between
the cases $\mu\in \bigl(-(1+a)^2,-(1-a)^2\bigr)$ and
$\mu=\mu(\pm 1)=-(1\mp a)^2$.\\
{\it Case 1:} $\mu\in \bigl(-(1+a)^2,-(1-a)^2\bigr)$, i.e. 
$\mu=\mu(e^{i\theta})$ with $\theta\in (0,\pi)$.

It is a straightforward consequence of the asymptotic 
behaviour of $\Phi_{\gamma}$ to infinity, see \eqref{asymptoticsPhiB},
that any non-zero linear combination of
the basis elements $\{\Phi_{e^{i\theta}}, \Phi_{e^{-i\theta}}\}$ of 
$V_{\mu(e^{i\theta})}^+$ does not lie in $M_{\mu(e^{i\theta})}$
(see Lemma \ref{l2} for
the definition of $M_{\mu}$), hence $\mu\not\in \sigma_p(L)$.\\
{\it Case 2:} $\mu=\mu(\pm 1)=-(1\mp a)^2$.
 
Observe that $\frac{\partial\Phi_{\gamma}}{\partial\gamma}|_{\gamma=\pm 1}
\in V_{\mu(\pm 1)}^+$ since
$\frac{\partial\mu}{\partial\gamma}(\gamma)|_{\gamma=\pm 1}=0$.
Using that the 
asymptotics to infinity of $\frac{\partial\Phi_{\gamma}}{\partial\gamma}$
is given by
\[\frac{\partial\Phi_{\gamma}}{\partial\gamma}(zq^{-m})=m\gamma^{-1}
(a\gamma)^m\bigl(1+{\mathcal{O}}(q^m)\bigr),\qquad m\to\infty,\]
we obtain that
\[W\bigl(\Phi_{\pm 1}, 
\frac{\partial\Phi_{\gamma}}{\partial\gamma}|_{\gamma=\pm 1}\bigr)(x)=
-aK\not=0,\qquad x\in I_+,
\]
compare with the proof of Lemma \ref{dep}. Combined with
Lemma \ref{generaleigenfunction}{\bf (i)}, it follows that
$\{\Phi_{\pm 1}, 
\frac{\partial\Phi_{\gamma}}{\partial\gamma}|_{\gamma=\pm 1}\}$
is a basis of $V_{\mu(\pm 1)}^+$.
It is easy to see, using asymptotics 
to infinity, that any non-zero linear combination of these basis elements
does not lie in $M_{\mu(\pm 1)}$, hence $\mu(\pm 1)\not\in \sigma_p(L)$.
\end{proof}
%%%%%%%%%%%%%%%%%%%%%%%%%%%%%%%%%%%%%%%%%%%%%%%%%%%%%%%%%%%%%%%%%%%%

Let ${\mathcal{D}}_{fin}\subset{\mathcal{H}}$
be the linear subspace of functions $f:I\to\CC$ with finite 
support. Observe that ${\mathcal{D}}_{fin}\subset {\mathcal{D}}\subset
{\mathcal{H}}$ as dense subspaces.
We define the big $q$-Jacobi function transform by
\begin{equation}\label{F}
\bigl({\mathcal{F}}f\bigr)(\gamma)=
\langle f,\phi_{\gamma}\rangle=\int_{-1}^{\infty(z)}
f(x){\overline{\phi_{\gamma}(x)}}\frac{d_qx}{p(x)}, \qquad f\in 
{\mathcal{D}}_{fin},\,\, \gamma\in \CC^*.
\end{equation}
In this section, we will regard the big $q$-Jacobi function transform 
${\mathcal{F}}f$ of $f\in {\mathcal{D}}_{fin}$ 
as a function on the unit circle 
$\TT$. Observe
that ${\mathcal{F}}f$ is $W$-invariant, where $W=\{\pm 1\}$ acts by
$\bigl(\ep.g\bigr)(\gamma)=g\bigl(\gamma^{\ep}\bigr)$, $\ep\in W$.

We define an absolutely continuous measure on $\TT$ by
\begin{equation}\label{mu}
d\nu(\gamma)=\frac{1}{4\pi iK}\frac{d\gamma}{c(\gamma)c(\gamma^{-1})\gamma}
=\frac{1}{4\pi iK}\frac{d\gamma}{|c(\gamma)|^2\gamma},
\end{equation}
where $K$ is the positive constant defined by \eqref{Kconstant}.
Observe that the measure $d\nu$ is a well defined measure on
$\TT$ since $(a,b,c)\in V$. In particular, possible zeros of the denominator
of the weight function $1/|c(\cdot)|^2$ are compensated by zeros of
the numerator. 

We first show that ${\mathcal{F}}$ extends uniquely to a partial isometry
${\mathcal{F}}: {\mathcal{H}}\to L_W^2\bigl(\TT,d\nu\bigr)$, where
$L_W^2\bigl(\TT,d\nu\bigr)$ is the $L^2$-space of $W$-invariant functions
with respect to $d\nu$, when $(a,b,c)\in V_z^{gen}$. 
We start with the following crucial consequence 
of \eqref{inversion} and Proposition \ref{G}.
%%%%%%%%%%%%%%%%%%%%%%%%%%%%%%%%%%%%%%%%%%%%%%%%%%%%%%%%%%%%%%%%%%%%
\begin{prop}\label{contE}
Let $(a,b,c)\in V_z^{gen}$. Choose 
$-(1+a)^2\leq \mu_1< \mu_2\leq -(1-a)^2$ and let 
$0\leq \theta_2<\theta_1\leq\pi$ such that $\mu_j=\mu(e^{i\theta_j})$
\textup{(}$j=1,2$\textup{)}. Then
\begin{equation}\label{contpart}
\langle E\bigl((\mu_1,\mu_2)\bigr)f,g\rangle=\frac{1}{2\pi K}
\int_{\theta_2}^{\theta_1}
\bigl({\mathcal{F}}f\bigr)(e^{i\theta})
{\overline{\bigl({\mathcal{F}}g\bigr)(e^{i\theta})}}
\frac{d\theta}{|c(e^{i\theta})|^2}
\end{equation}
for all $f,g\in {\mathcal{D}}_{fin}$. 
\end{prop}
%%%%%%%%%%%%%%%%%%%%%%%%%%%%%%%%%%%%%%%%%%%%%%%%%%%%%%%%%%%%%%%%%%
\begin{proof}
For $\mu\in\CC\setminus \RR$ we write $\gamma[\mu]\in\CC\setminus\RR$ for the
unique complex number with modulus less than $1$ such that
$\mu(\gamma[\mu])=\mu$. Fix $f,g\in {\mathcal{D}}_{fin}$.
By a straightforward computation using Proposition
\ref{G}, we have for $\mu\in \RR$ and $\epsilon >0$
that
\begin{equation}\label{crucfor}
\begin{split}
\langle \left(L-\bigl(\mu\pm i\epsilon\bigr)\right)^{-1}f,g\rangle=
&{\underset{\stackrel{(x,y)\in I\times I}{x\leq y}}{\iint}}
\frac{\phi_{\gamma[\mu\pm i\epsilon]}(x)
\Phi_{\gamma[\mu\pm i\epsilon]}(y)}{W\bigl(\gamma[\mu\pm
i\epsilon]\bigr)}\bigl(1-\frac{1}{2}\delta_{x,y}\bigr)\\
&\qquad.\bigl(f(x){\overline{g(y)}}+f(y){\overline{g(x)}}
\bigr)\frac{d_qx}{p(x)}\frac{d_qy}{p(y)},
\end{split}
\end{equation}
where $\delta_{x,y}$ is the Kronecker-delta.
Let $-(1+a)^2<\mu_1\leq\mu\leq\mu_2<-(1-a)^2$ and let
$0<\theta_2\leq\theta\leq\theta_1<\pi$ such that $\mu=\mu(e^{i\theta})$ and
$\mu_j=\mu(\theta_j)$ for $j=1,2$.
Then we have
\begin{equation}\label{gammalim} 
\lim_{\epsilon\downarrow 0}\gamma[\mu\pm i\epsilon]=e^{\mp i\theta}.
\end{equation}
Using the connection coefficients formula given in 
Theorem \ref{samenvatting}{\bf (iv)} and using
the fact that $W(\gamma)^{-1}$ and $\Phi_{\gamma^{\pm 1}}(x)$
($x\in I$) are regular at $\gamma\in 
\TT\setminus \{\pm 1\}$, we obtain
\begin{equation*}
\begin{split}
\lim_{\ep\downarrow 0}\left(\frac{\phi_{\gamma[\mu+i\ep]}(x)
\Phi_{\gamma[\mu+i\ep]}(y)}{W\bigl(\gamma[\mu+i\ep]\bigr)}-\right.
&\left.\frac{\phi_{\gamma[\mu-i\ep]}(x)\Phi_{\gamma[\mu-i\ep]}(y)}
{W\bigl(\gamma[\mu-i\ep]\bigr)}\right)\\
=&\frac{\phi_{e^{i\theta}}(x)\Phi_{e^{-i\theta}}(y)}
{W\bigl(e^{-i\theta}\bigr)}-\frac{\phi_{e^{i\theta}}(x)\Phi_{e^{i\theta}}(y)}
{W\bigl(e^{i\theta}\bigr)}\\
=&\frac{1}{aK}\frac{\phi_{e^{i\theta}}(x)\phi_{e^{i\theta}}(y)}
{c\bigl(e^{i\theta}\bigr)c\bigl(e^{-i\theta}\bigr)
\bigl(e^{-i\theta}-e^{i\theta}\bigr)}.
\end{split}
\end{equation*}
It follows now by symmetrization of the double $q$-Jackson integral that
\begin{equation*}
\begin{split}
\lim_{\epsilon\downarrow 0}&\left(\langle 
\left(L-\bigl(\mu+i\epsilon\bigr)\right)^{-1}f,g\rangle-
\langle \left(L-\bigl(\mu-i\epsilon\bigr)\right)^{-1}f,g\rangle\right)\\
&=\frac{1}{aK}{\underset{\stackrel{(x,y)\in I\times I}{x\leq y}}{\iint}}
\frac{\phi_{e^{i\theta}}(x)\phi_{e^{i\theta}}(y)
\bigl(f(x){\overline{g(y)}}+f(y){\overline{g(x)}}\bigr)
\bigl(1-\frac{1}{2}\delta_{x,y}\bigr)}{|c(e^{i\theta})|^2
\bigl(e^{-i\theta}-e^{i\theta}\bigr)}\frac{d_qx}{p(x)}\frac{d_qy}{p(y)}\\
&=\frac{1}{aK}
\frac{\bigl({\mathcal{F}}f\bigr)(e^{i\theta})
{\overline{\bigl({\mathcal{F}}g\bigr)(e^{i\theta})}}}
{|c\bigl(e^{i\theta}\bigr)|^2\bigl(e^{-i\theta}-e^{i\theta}\bigr)}.
\end{split}
\end{equation*}
The proposition follows now for all $-(1+a)^2<\mu_1<\mu_2<-(1-a)^2$
using \eqref{inversion}
and changing the integration variable to $\theta$ using 
the map $\theta\mapsto \mu(e^{i\theta})$, see \eqref{eigenfunction}. 
The result now also holds when $\mu_1=-(1+a)^2$
or $\mu_2=-(1-a)^2$ since the spectral measure $E$ is countably additive.
\end{proof}
%%%%%%%%%%%%%%%%%%%%%%%%%%%%%%%%%%%%%%%%%%%%%%%%%%%%%%%%%%%%%%%%%%%%%%%%

%%%%%%%%%%%%%%%%%%%%%%%%%%%%%%%%%%%%%%%%%%%%%%%%%%%%%%%%%%%%%%%%%%%%%%%%
\begin{cor}\label{gevolg1}
The big $q$-Jacobi function transform ${\mathcal{F}}$ uniquely
extends 
to a continuous linear mapping ${\mathcal{F}}_c: {\mathcal{H}}\to 
L_W^2\bigl(\TT,d\nu\bigr)$. 

If $(a,b,c)\in V_z^{gen}$, then 
${\mathcal{F}}_c: {\mathcal{H}}\to L_W^2\bigl(\TT,d\nu\bigr)$
factorizes through the orthogonal projection 
$P_c=E\bigl([-(1+a)^2,-(1-a)^2]\bigr)$ \textup{(}i.e. 
${\mathcal{F}}_c={\mathcal{F}}_c\circ P_c$\textup{)}, 
and the restriction of ${\mathcal{F}}_c$ 
to the range ${\mathcal{R}}(P_c)$ of $P_c$ is an isometric isomorphism
onto the range ${\mathcal{R}}({\mathcal{F}}_c)\subset L_W^2\bigl(\TT,d\nu\bigr)$
of ${\mathcal{F}}_c$.
\end{cor}
%%%%%%%%%%%%%%%%%%%%%%%%%%%%%%%%%%%%%%%%%%%%%%%%%%%%%%%%%%%%%%%%%%%%%%%%
\begin{proof}
We assume first that $(a,b,c)\in V_z^{gen}$.
In view of Proposition \ref{contpart} applied to the special case
$\mu_1=\mu(-1)=-(1+a)^2$ and $\mu_2=\mu(1)=-(1-a)^2$, it suffices 
to observe that $E\bigl(\{\mu(\pm 1)\}\bigr)=0$, which is a consequence
of Lemma \ref{contspectrum} and \cite[Theorem 13.27]{R}.

Finally, observe that the 
inequality $\|{\mathcal{F}}f\|_2^2\leq \|f\|^2$ for 
$f\in {\mathcal{D}}_{fin}$ and $(a,b,c)\in V_z^{gen}$, where $\|.\|_2$
is the norm of $L_W^2\bigl(\TT,d\nu\bigr)$, holds 
for all $(a,b,c)\in V$ by continuity. 
Hence ${\mathcal{F}}$ can be uniquely extended to
a continuous linear map ${\mathcal{F}}_c: {\mathcal{H}}\to 
L_W^2\bigl(\CC^*,d\nu\bigr)$
for all parameters $(a,b,c)\in V$.
\end{proof}
%%%%%%%%%%%%%%%%%%%%%%%%%%%%%%%%%%%%%%%%%%%%%%%%%%%%%%%%%%%%%%%%%

%%%%%%%%%%%%%%%%%%%%%%%%%%%%%%%%%%%%%%%%%%%%%%%%%%%%%%%%%%%%%%%%%
\begin{Def}
We call ${\mathcal{F}}_c$ the continuous part of the big
$q$-Jacobi function transform.
\end{Def}
%%%%%%%%%%%%%%%%%%%%%%%%%%%%%%%%%%%%%%%%%%%%%%%%%%%%%%%%%%%%%%%%%

Observe that the limit $\langle f,\phi_{\gamma}\rangle_l=
\lim_{k,m\to\infty}\langle f,\phi_{\gamma}\rangle_{k;l,m}$
exists for all $\gamma\in \TT$ since 
$\phi_{\gamma}$ is continuously differentiable at the origin, see
Lemma \ref{oknul}. It follows that
\begin{equation}\label{pointwise}
\bigl({\mathcal{F}}_cf\bigr)(\gamma)=
\lim_{l\to-\infty}\langle f,\phi_{\gamma}\rangle_l,\qquad
\gamma\in\TT \,\,\,\hbox{ a.e.}
\end{equation}
for $f\in {\mathcal{H}}$, where $\langle .,. \rangle_l$  
($l\in\ZZ$) is the truncated
inner product
\begin{equation}\label{trunc}
\langle f,g\rangle_l=\int_{-1}^{zq^l}f(x){\overline{g(x)}}\frac{d_qx}
{p(x)}.
\end{equation}

In the remainder of this section we show that ${\mathcal{F}}_c$ is surjective
and we give an explicit
formula for the inverse of the isometric isomorphism
${\mathcal{F}}_c: {\mathcal{R}}(P_c)\to L_W^2\bigl(\TT,d\nu)$.
The methods we employ are similar to the ones employed by G{\"o}tze
\cite{G} and by Braaksma and Meulenbeld \cite{BM} for the classical
Jacobi function transform, and by Kakehi \cite{K} and Kakehi, Masuda and
Ueno \cite{KMU} for the little $q$-Jacobi function transform.

%%%%%%%%%%%%%%%%%%%%%%%%%%%%%%%%%%%%%%%%%%%%%%%%%%%%%%%%%%%%%%%%%%
\begin{lem}\label{Wphiphi}
Let $\gamma,\delta\in \TT$ with $\mu(\gamma)\not=\mu(\delta)$.
For $k\in \ZZ_+$, $l,m\in \ZZ$ with $l<m$, we have 
\[ \langle \phi_{\gamma},\phi_{\delta}\rangle_{k;l,m}=
\frac{W(\phi_{\gamma},\phi_{\delta})(zq^{l-1})
-W(\phi_{\gamma},\phi_{\delta})(zq^m)
+W(\phi_{\gamma},\phi_{\delta})(-q^k)}
{\mu(\gamma)-\mu(\delta)}.
\]
\end{lem}
%%%%%%%%%%%%%%%%%%%%%%%%%%%%%%%%%%%%%%%%%%%%%%%%%%%%%%%%%%%%%%%%%%%
\begin{proof}
For $\gamma,\delta\in\TT$ we have that 
$\bigl(L\phi_{\gamma}\bigr){\overline{\phi_{\delta}}}-
\phi_{\gamma}{\overline{\bigl(L\phi_{\delta}\bigr)}}=
\bigl(\mu(\gamma)-\mu(\delta)\bigr)
\phi_{\gamma}{\overline{\phi_{\delta}}}$ on $I$.
The proof follows now from Lemma \ref{stoktermen} since 
$\phi_{\gamma}:I\to\CC$ is real valued for $\gamma\in \TT$.
\end{proof}
%%%%%%%%%%%%%%%%%%%%%%%%%%%%%%%%%%%%%%%%%%%%%%%%%%%%%%%%%%%%%%%%%%

We define a linear map 
${\mathcal{G}}_c: L_W^2\bigl(\TT,d\nu\bigr)\to {\mathcal{F}}(I)$ by
\begin{equation}\label{Finvers}
\bigl({\mathcal{G}}_cg\bigr)(x)=\int_{\TT}g(\gamma)\phi_{\gamma}(x)d\nu(\gamma),
\qquad x\in I.
\end{equation}
Observe that 
\begin{equation}
\langle {\mathcal{G}}_cg_1, {\mathcal{G}}_cg_2
\rangle_{k;l,m}=
\int_\TT\int_\TT g_1(\gamma){\overline{g_2(\gamma')}}
\langle \phi_{\gamma}, \phi_{\gamma'}\rangle_{k;l,m}d\nu(\gamma)
d\nu(\gamma')
\end{equation}
for $k\in \ZZ_+$, $l<m$ in $\ZZ$ and
$g_1,g_2\in L_W^2\bigl(\TT,d\nu\bigr)$.
%%%%%%%%%%%%%%%%%%%%%%%%%%%%%%%%%%%%%%%%%%%%%%%%%%%%%%%%%%%%%%%%%%%
\begin{lem}\label{reduction}
The limit $\langle {\mathcal{G}}_cg_1, {\mathcal{G}}_cg_2\rangle_l
=\lim_{k,m\to\infty}\langle {\mathcal{G}}_cg_1, 
{\mathcal{G}}_cg_2\rangle_{k;l,m}$ exists for all
$g_1,g_2\in L_W^2\bigl(\TT,d\nu\bigr)$, and
\begin{equation*}
\begin{split}
\langle {\mathcal{G}}_cg_1,
{\mathcal{G}}_cg_2\rangle_l=
&\int_\TT\int_\TT g_1(\gamma){\overline{g_2(\gamma')}}
\langle \phi_{\gamma}, \phi_{\gamma'}\rangle_ld\nu(\gamma)d\nu(\gamma')\\
=&\int_\TT\int_\TT g_1(\gamma){\overline{g_2(\gamma')}}
\frac{W\bigl(\phi_{\gamma},\phi_{\gamma'}\bigr)\bigl(zq^{l-1}\bigr)}
{\bigl(\mu(\gamma)-\mu(\gamma')\bigr)}
d\nu(\gamma)d\nu(\gamma').
\end{split}
\end{equation*}
\end{lem}
%%%%%%%%%%%%%%%%%%%%%%%%%%%%%%%%%%%%%%%%%%%%%%%%%%%%%%%%%%%%%%%%%%
\begin{proof}
Observe that there exists a constant $K_l>0$ 
such that $|\langle \phi_{\gamma},\phi_{\gamma'}\rangle_{k;l,m}|
\leq K_l$ for all $\gamma,\gamma'\in \TT$ and for all
$k,m\in\ZZ_+$ with $m>l$. The first equality follows then from
Lebesgue's dominated convergence theorem.

The second equality follows from Lemma \ref{Wphiphi}, using that
$W\bigl(\phi_{\gamma},\phi_{\gamma'}\bigr)(0^+)=
W\bigl(\phi_{\gamma},\phi_{\gamma'}\bigr)(0^-)$ since $\phi_{\gamma}$ is 
continuously differentiable at the origin, see Lemma \ref{oknul}.
\end{proof}
%%%%%%%%%%%%%%%%%%%%%%%%%%%%%%%%%%%%%%%%%%%%%%%%%%%%%%%%%%%%%%%%%%

Finally, we determine the limit 
$\langle {\mathcal{G}}_cg_1,{\mathcal{G}}_cg_2\rangle=\lim_{l\to-\infty}
\langle {\mathcal{G}}_cg_1,{\mathcal{G}}_cg_2\rangle_l$.
The result is as follows.
%%%%%%%%%%%%%%%%%%%%%%%%%%%%%%%%%%%%%%%%%%%%%%%%%%%%%%%%%%%%%%%%%
\begin{prop}\label{weaklimit}
The limit $\langle {\mathcal{G}}_cg_1, {\mathcal{G}}_cg_2\rangle=
\lim_{l\to-\infty}\langle {\mathcal{G}}_cg_1, {\mathcal{G}}_cg_2\rangle_l$ 
exists for all $g_1,g_2\in L_W^2\bigl(\TT,d\nu\bigr)$, and
\begin{equation}\label{i}
\langle {\mathcal{G}}_cg_1, {\mathcal{G}}_cg_2\rangle=
\int_{\TT}g_1(\gamma){\overline{g_2(\gamma)}}d\nu(\gamma).
\end{equation}
In particular, ${\mathcal{G}}_c :L_W^2\bigl(\TT,d\nu\bigr)\to {\mathcal{H}}$ 
is an isometric isomorphism
onto the range ${\mathcal{R}}\bigl({\mathcal{G}}_c\bigr)\subset {\mathcal{H}}$
of ${\mathcal{G}}_c$.
\end{prop}
%%%%%%%%%%%%%%%%%%%%%%%%%%%%%%%%%%%%%%%%%%%%%%%%%%%%%%%%%%%%%%%%%%%
\begin{proof}
We only sketch the proof, since it is similar to the
little $q$-Jacobi case, see \cite[Proposition 6.1]{KMU} and 
\cite[Proposition 7.4]{K}. Let
$C_W(\TT)$ be the algebra of complex valued, 
continuous, $W$-invariant functions on $\TT$.
We fix $g_1,g_2\in C_W(\TT)$ such that $\{-1,1\}$ is not in the support of 
$g_1$ and $g_2$. It suffices to give a proof of \eqref{i} for such functions
$g_1$ and $g_2$.
We start with the second expression of $\langle {\mathcal{G}}_{l}g_1,
{\mathcal{G}}_lg_2\rangle_l$ in Lemma
\ref{reduction}, and replace $\phi_{\gamma}$ by its $c$-function
expansion, see Theorem \ref{samenvatting}{\bf (iv)}.
Using the estimate
\[\sup_{\delta\leq \theta\not=\theta'\leq\pi-\delta}
\left|\frac{R_k(e^{\pm i\theta})-R_k(e^{\pm i\theta'})}
{e^{i\theta}-e^{i\theta'}}\right|={\mathcal{O}}(k(qa)^k),\qquad k\to\infty
\]
for $0<\delta<\pi/2$, where 
$R_k(\gamma)=\Phi_{\gamma}(zq^{-k})-(a\gamma)^k$ (cf. 
\eqref{asymptoticsPhiB}), which can easily be proved using the mean 
value theorem, we may replace in 
the expression of $\lim_{l\to-\infty}\langle {\mathcal{G}}_cg_1,
{\mathcal{G}}_cg_2\rangle_l$ the 
function $\Phi_{\gamma}(zq^{l-1})$ by
its asymptotic value $(a\gamma)^{1-l}$ at $\infty$. 
Combined with \eqref{rinfty} it follows now from the bounded convergence 
theorem that
\[
\lim_{m\to\infty}\langle {\mathcal{G}}_cg_1,{\mathcal{G}}_cg_2\rangle_{1-m}
=\frac{-a}{4\pi^2K}
\lim_{m\to\infty}\int_{0}^\pi
\int_{0}^\pi\frac{g_1(e^{i\theta}){\overline{g_2(e^{i\theta'})}}
s_m(\theta,\theta')d\theta d\theta'}{|c(e^{i\theta})|^2
|c(e^{i\theta'})|^2(\mu(e^{i\theta'})-\mu(e^{i\theta}))}
\]
provided that the limit in the right hand side of the equality exists, with
\[s_m(\theta,\theta')=\sum_{\epsilon,\xi\in\{\pm 1\}}
c(e^{i\epsilon\theta})c(e^{i\xi\theta'})e^{i(m-1)(\epsilon\theta+\xi\theta')}
(e^{i\epsilon\theta}-e^{i\xi\theta'}).
\]
Since $c(\gamma)$ is continuous and non-zero on $\TT\setminus \{\pm 1\}$,
the Riemann-Lebesgue lemma implies that 
the contributions of the sums of $s_m$ with $\ep\xi>0$ tend to 
zero in the limit, so $s_m(\theta,\theta')$ may be replaced by
\begin{equation*}
\begin{split}
t_m(\theta,\theta')=
&c(e^{i\theta})c(e^{-i\theta'})e^{i(m-1)(\theta-\theta')}(e^{i\theta}-
e^{-i\theta'})\\
&+c(e^{-i\theta})c(e^{i\theta'})e^{i(m-1)(\theta'-\theta)}(e^{-i\theta}-
e^{i\theta'})\\
&=-4c(e^{i\theta})c(e^{-i\theta'})\sin\left(\frac{\theta+\theta'}{2}\right)
\sin\left(\frac{(2m-1)(\theta-\theta')}{2}\right)\\
&+\bigl(c(e^{-i\theta})c(e^{i\theta'})-c(e^{i\theta})c(e^{-i\theta'})\bigr)
e^{i(m-1)(\theta'-\theta)}(e^{-i\theta}-e^{i\theta'}).
\end{split}
\end{equation*}
Applying the Riemann-Lebesgue lemma again and using
\[\mu(e^{i\theta'})-\mu(e^{i\theta})=2a\bigl(\cos(\theta')-\cos(\theta)\bigr)=
4a\sin\left(\frac{\theta+\theta'}{2}\right)
\sin\left(\frac{\theta-\theta'}{2}\right),
\]
we arrive at
\[\lim_{m\to\infty}\langle{\mathcal{G}}_cg_1,
{\mathcal{G}}_cg_2\rangle_{1-m}=\lim_{m\to\infty}\frac{1}{4\pi^2K}
\int_0^\pi\int_0^\pi 
\frac{g_1(e^{i\theta}){\overline{g_2(e^{i\theta'})}}}{c(e^{-i\theta})
c(e^{i\theta'})}D_m(\theta,\theta')d\theta d\theta'
\]
provided that the limit in the right hand side of the equality exists,
where $D_m$ is the Dirichlet kernel,
\[D_m(\theta,\theta')= \frac{\sin\bigl((2m-1)(\theta-\theta')/2\bigr)}
{\sin\bigl((\theta-\theta')/2\bigr)}.\]
The result follows now from the well known $L^2$-properties of the Dirichlet
kernel.
\end{proof}
%%%%%%%%%%%%%%%%%%%%%%%%%%%%%%%%%%%%%%%%%%%%%%%%%%%%%%%%%%%%%%%%%%%

Recall the notation $P_c$ for the orthogonal projection 
$E\bigl([-(1+a)^2,-(1-a)^2]\bigr)$. 

%%%%%%%%%%%%%%%%%%%%%%%%%%%%%%%%%%%%%%%%%%%%%%%%%%%%%%%%%%%%%%%%%%%%%%%%
\begin{prop}\label{mainprop}

{\bf (i)} ${\mathcal{F}}_c|_{\mathcal{R}({\mathcal{G}}_c)}:
\mathcal{R}({\mathcal{G}}_c)\rightarrow L_W^2\bigl(\TT,d\nu\bigr)$
is a surjective isometric isomorphism. Its inverse is given by
 ${\mathcal{G}}_c: L_W^2\bigl(\TT,d\nu\bigr)\rightarrow 
 \mathcal{R}({\mathcal{G}}_c)$.

{\bf (ii)}
We have ${\mathcal{R}}({\mathcal{G}}_c)={\mathcal{R}}(P_c)$
for $(a,b,c)\in V_z^{gen}$.
\end{prop}
%%%%%%%%%%%%%%%%%%%%%%%%%%%%%%%%%%%%%%%%%%%%%%%%%%%%%%%%%%%%%%%%%%%%%%
\begin{proof} 
{\bf (i)} It follows from Corollary \ref{gevolg1}, Lemma \ref{reduction} 
and Proposition \ref{weaklimit} that 
\begin{equation*}
\begin{split}
\int_{\TT}\bigl({\mathcal{F}}_c\bigl({\mathcal{G}}_cf\bigr)\bigr)(\gamma)
{\overline{g(\gamma)}}d\nu(\gamma)
=&\lim_{l\to-\infty}\lim_{k,l\to\infty}
\int_{\TT}\langle {\mathcal{G}}_cf, \phi_{\gamma}\rangle_{k;l,m}
{\overline{g(\gamma)}}d\nu(\gamma)\\
=&\lim_{l\to-\infty}\lim_{k,l\to\infty}\langle {\mathcal{G}}_cf, 
{\mathcal{G}}_cg\rangle_{k;l,m}=\int_{\TT}f(\gamma){\overline{g(\gamma)}}
d\nu(\gamma)
\end{split}
\end{equation*}
for all $f,g\in L_W^2\bigl(\TT,d\nu\bigr)$.
Hence ${\mathcal{F}}_c\circ
{\mathcal{G}}_c$ is the identity on $L^2_W\bigl(\TT,d\nu\bigr)$ and 
${\mathcal{R}}\bigl({\mathcal{F}}_c\bigr)=L_W^2\bigl(\TT,d\nu\bigr)$.
Consequently, we have ${\mathcal{G}}_c\circ 
{\mathcal{F}}_c|_{{\mathcal{R}}({\mathcal{G}}_c)}=
\hbox{Id}_{{\mathcal{R}}({\mathcal{G}}_c)}$.

{\bf (ii)}
%It remains to show that ${\mathcal{R}}({\mathcal{G}}_c)={\mathcal{R}}(P_c)$.  
Assume that $(a,b,c)\in V_z^{gen}$.
Let $f\in {\mathcal{R}}({\mathcal{G}}_c)$.
Observe that $\|f\|=\|{\mathcal{F}}_cf\|_2$ by the previous paragraph and
by Proposition \ref{weaklimit}, 
where $\|.\|_2$ is the norm of 
$L_W^2\bigl(\TT,d\nu\bigr)$.
By Corollary \ref{gevolg1}, this implies that
$f\in {\mathcal{R}}(P_c)$.

Let $f\in {\mathcal{R}}(P_c)$. We have seen that 
${\mathcal{F}}_c\bigl({\mathcal{R}}({\mathcal{G}}_c)\bigr)=
L_W^2\bigl(\TT,d\nu\bigr)$, hence there exists a function
$g\in {\mathcal{R}}({\mathcal{G}}_c)$ such that ${\mathcal{F}}_cf=
{\mathcal{F}}_cg$. Now $f-g\in {\mathcal{R}}(P_c)$ by the previous paragraph, 
and ${\mathcal{F}}_c|_{{\mathcal{R}}(P_c)}$ is 
injective by Corollary \ref{gevolg1}, 
hence $f=g\in {\mathcal{R}}({\mathcal{G}}_c)$, as desired.

\end{proof}
%%%%%%%%%%%%%%%%%%%%%%%%%%%%%%%%%%%%%%%%%%%%%%%%%%%%%%%%%%%%%%%%%

%%%%%%%%%%%%%%%%%%%%%%%%%%%%%%%%%%%%%%%%%%%%%%%%%%%%%%%%%%%%
%%                                                        %%
%%              The point spectrum                        %%
%%                                                        %%
%%%%%%%%%%%%%%%%%%%%%%%%%%%%%%%%%%%%%%%%%%%%%%%%%%%%%%%%%%%%

\section{The point spectrum}\label{pointmass}

In this section we determine the resolution of the identity
$E$ on $\RR\setminus [-(1+a)^2,-(1-a)^2]$. 
Observe that the kernel $K_{\gamma}(x,y)$ 
of the Green function is meromorphic as function of $\gamma\in {\mathbb{D}}$, 
where
${\mathbb{D}}=\{ v\in \CC \, | \, 0<|v|< 1 \}$ is the punctured open
unit disc in the complex plane. Since $\phi_{\gamma}(x)$ and
$\Phi_{\gamma}(x)$ are analytic at $\gamma\in {\mathbb{D}}$, see Theorem 
\ref{samenvatting}{\bf (ii)} for $\Phi_{\gamma}$, we have that
the poles of $K_{\gamma}(x,y)$ coincide with the 
poles of $W(\gamma)^{-1}$, which in turn coincide with the zeros
of the map $\gamma\mapsto c(\gamma^{-1})$ 
$(\gamma\in {\mathbb{D}})$, where
$c(\cdot)$ is the $c$-function defined by \eqref{cB}.
In particular, the poles of $K_{\gamma}(x,y)$ in 
$\gamma\in {\mathbb{D}}$
are independent of $x,y\in I$ and are given by the set
\begin{equation}\label{Sd}
S=\left\{\frac{1}{q^ke} \,\, \left|\right. \,\, e\in \{a,b,c\}, k\in \ZZ_+: 
\,\, \frac{1}{q^{k}e}<1
\right\} \cup \left\{-\frac{abcz}{q^{k+1}}
 \,\, | \,\, k\in \ZZ: \,\,
\frac{abcz}{q^{k+1}}<1 \right\}.
\end{equation}
Observe that these poles of $K_{\gamma}(x,y)$ are simple since 
at most one of the parameters $a,b,c$ is $\geq 1$.
Furthermore, observe that $S\subset S_{reg}$ when $(a,b,c)\in V_z^{gen}$.

%%%%%%%%%%%%%%%%%%%%%%%%%%%%%%%%%%%%%%%%%%%%%%%%%%%%%%%%%%%%%%%%%%%%
\begin{lem}\label{Ltwo}
The function $\phi_{\tilde{\gamma}}$ \textup{(}$\tilde{\gamma}\in S$\textup{)}
is an eigenfunction of $(L,{\mathcal{D}})$ with eigenvalue 
$\mu(\tilde{\gamma})$. In particular, $\mu(S)$ is contained in
the point spectrum $\sigma_p(L)$ of $(L,{\mathcal{D}})$.
\end{lem}
%%%%%%%%%%%%%%%%%%%%%%%%%%%%%%%%%%%%%%%%%%%%%%%%%%%%%%%%%%%%%%%%%%%%%
\begin{proof}
Let $\tilde{\gamma}\in S$. Then 
$\phi_{\tilde{\gamma}}\in V_{\mu(\tilde{\gamma})}$ by Lemma \ref{oknul}.
Furthermore, if $\tilde{\gamma}\not\in S_{reg}$, then 
\begin{equation}\label{phiatpositive}
\phi_{\tilde{\gamma}}|_{I_+}=
c(\tilde{\gamma})\Phi_{\tilde{\gamma}}\in M_{\mu(\tilde{\gamma})},
\end{equation}
where the first equality follows from Proposition \ref{conncoef}.
Since $c(\gamma)$ and $\Phi_{\gamma}\in V_{\mu(\gamma)}^+$
are regular at $\gamma\in S\cap S_{reg}$, it follows that 
\eqref{phiatpositive} is valid for all $\tilde{\gamma}\in S$.
Hence $\phi_{\tilde{\gamma}}\in {\mathcal{D}}$ and 
$L\phi_{\tilde{\gamma}}=\mu(\tilde{\gamma})\phi_{\tilde{\gamma}}$
on $I$ for all $\tilde{\gamma}\in S$, as desired.
\end{proof}
%%%%%%%%%%%%%%%%%%%%%%%%%%%%%%%%%%%%%%%%%%%%%%%%%%%%%%%%%%%%%%%%%%%%%%%

%%%%%%%%%%%%%%%%%%%%%%%%%%%%%%%%%%%%%%%%%%%%%%%%%%%%%%%%%%%%%%%%%%%%%%%
\begin{rem}\label{phiateverywhere}
If $(a,b,c)\in V_z^{gen}$, then $S\subset S_{reg}$ and
$\phi_{\tilde{\gamma}}|_{I_+}=c(\tilde{\gamma})\Phi_{\tilde{\gamma}}$
by Theorem \ref{samenvatting}{\bf (iv)}, so 
formula \eqref{phiatpositive} extends for these parameter values
to the whole $q$-interval $I$.
\end{rem}
%%%%%%%%%%%%%%%%%%%%%%%%%%%%%%%%%%%%%%%%%%%%%%%%%%%%%%%%%%%%%%%%%%%%%%%%%

Let $\tilde{\gamma}\in S$.
It follows from Lemma \ref{Ltwo} that 
the linear functional 
$f\mapsto \bigl({\mathcal{F}}f\bigr)(\tilde{\gamma})$ with
$f\in {\mathcal{D}}_{fin}$ is bounded, where ${\mathcal{F}}$ is the 
big $q$-Jacobi function 
transform defined by \eqref{F}, hence it uniquely extends to
a continuous linear functional on ${\mathcal{H}}$. It is given explicitly by
\begin{equation}
\bigl({\mathcal{F}}f\bigr)(\tilde{\gamma})=
\langle f,\phi_{\tilde{\gamma}}\rangle,\qquad
f\in {\mathcal{H}}.
\end{equation}
We define the weight $d\nu\bigl(\{\tilde{\gamma}\}\bigr)$ for 
$\tilde{\gamma}\in S$ and $(a,b,c)\in V_z^{gen}$ by
\begin{equation}\label{discrweight}
d\nu\bigl(\{\tilde{\gamma}\}\bigr)=\frac{1}{K}
\underset{\gamma=\tilde{\gamma}}{\hbox{Res}}
\left(\frac{-1}{\gamma c(\gamma)c(\gamma^{-1})}\right)
=\frac{1}{K}\underset{\gamma=\tilde{\gamma}^{-1}}{\hbox{Res}}
\left(\frac{1}{\gamma c(\gamma)c(\gamma^{-1})}\right),
\end{equation}
where $K$ is the positive constant defined by \eqref{Kconstant}.
Observe that the poles $\tilde{\gamma}\in S$ of 
the meromorphic function 
$\gamma\mapsto\bigl(c(\gamma)c(\gamma^{-1})\bigr)^{-1}$ are 
simple for parameters $(a,b,c)\in V_z^{gen}$.
We derive explicit expressions for the
discrete weights $d\nu\bigl(\{\tilde{\gamma}\}\bigr)$ by
relating the continuous part of the Plancherel measure with
the weight function $\Delta(x)=\Delta(x;t_0,t_1,t_2,t_3)$ of the Askey-Wilson 
polynomials, 
\begin{equation}\label{AskeyWilsonmeasure}
\Delta(x)=\frac{\bigl(x^2,1/x^2;q\bigr)_{\infty}}
{\prod_{j=0}^3\bigl(t_jx,t_j/x;q\bigr)_{\infty}}.
\end{equation}
Observe that
\begin{equation}\label{relPlanchAW}
\frac{1}{Kc(\gamma)c(\gamma^{-1})\gamma}=M
\frac{\Delta(\gamma;e,f,-q/efgz,-efgz)}{\bigl(g\gamma, 
g/\gamma;q\bigr)_{\infty}\gamma},
\end{equation}
where $\{e,f,g\}$ is an arbitrary permutation of $\{a,b,c\}$ 
(taking multiplicity into account),
and where $M=M(a,b,c;z)$ is the positive constant
\begin{equation}\label{M}
M=\frac{1}{K}\bigl(ab,ac;q\bigr)_{\infty}^2\theta(-bcz,-bcz)
=\frac{\bigl(ab,ac;q\bigr)_{\infty}^2\theta(-abz,-acz,-bcz)}
{(1-q)z\theta(-1/z)}.
\end{equation}
For a simple pole $eq^k$ of $\Delta(x;t_0,t_1,t_2,t_3)$, where 
$e\in \{t_j\}_{j=0}^3$
and $k\in\ZZ_+$, we have the explicit formula
\begin{equation}\label{discAW}
\begin{split}
{\underset{\gamma=eq^k}{\hbox{Res}}}&\left(\frac{\Delta(x;t_0,t_1,t_2,t_3)}{x}
\right)=
\frac{\bigl(e^{-2};q\bigr)_{\infty}}
{\bigl(q,ef,f/e,eg,g/e,eh,h/e;q\bigr)_{\infty}}\\
&\qquad\qquad\qquad\qquad.\frac{\bigl(e^2,ef,eg,eh;q\bigr)_k}
{\big(q,qe/f,qe/g,qe/h;q\bigr)_k}\frac{(1-e^2q^{2k})}{(1-e^2)}
\left(\frac{q}{efgh}\right)^k,
\end{split}
\end{equation}
see \cite[(7.5.22)]{GR}, where $\{f,g,h\}$ is such that $\{e,f,g,h\}=
\{t_0,t_1,t_2,t_3\}$
(taking multiplicity into account). It is well known that the right hand
side of \eqref{discAW} is well defined and positive for 
real parameters $t_i$ such that the  
$t_it_j$ ($i\not=j$) are strictly less than one.

Combined with \eqref{discrweight}, \eqref{relPlanchAW} and 
\eqref{qshifttransform}, 
we obtain for $e\in \{a,b,c\}$ and $k\in \ZZ_+$ such that $1/q^ke<1$,
\begin{equation}\label{normposdisc}
\begin{split}
d\nu\bigl(\{ 1/q^ke\}\bigr)=
&M\frac{\bigl(e^{-2};q\bigr)_{\infty}}
{\bigl(q,ef,f/e,eg,g/e;q\bigr)_{\infty}\theta(-fgz,-e^2fgz)}\\
&.\frac{\bigl(e^2,ef,eg;q\bigr)_k}
{\bigl(q,qe/f,qe/g;q\bigr)_k}\frac{(1-e^2q^{2k})}{(1-e^2)}
\left(\frac{-q^{(k+1)/2}}{fg}\right)^k,
\end{split}
\end{equation}
where $f,g$ are such that $\{e,f,g\}=\{a,b,c\}$ (taking multiplicity
into account). By the positivity of the residues of the Askey-Wilson measure,
it follows that the right hand side of \eqref{normposdisc}
is well defined, regular and strictly positive for parameters 
$(a,b,c)\in V$.
 
Similarly, we obtain for $k\in \ZZ$ such that $abcz/q^{k+1}<1$, 
\begin{equation}\label{normstrange}
\begin{split}
&d\nu\bigl(\{-abcz/q^{k+1}\}\bigr)=\\
&\quad=\frac{M}{\bigl(q,q,-q/abz,-q/acz,-q/bcz,
-a^2bcz/q,-ab^2cz/q,-abc^2z/q;q\bigr)_{\infty}}\\
&\qquad.\frac{\bigl(-q/abz,-q/acz,-q/bcz;q\bigr)_{k}}
{\bigl(-1/a^2bc,-1/ab^2c,-1/abc^2;q\bigr)_{k}}
\left(q^{2k}-\frac{a^2b^2c^2z^2}{q^2}\right)
\left(\frac{-q^{(k+3)/2}}{a^2b^2c^2z}\right)^k.
\end{split}
\end{equation}
The right hand side of \eqref{normstrange}
is well defined, regular and strictly positive for 
parameters $(a,b,c)\in V$.
We will sometimes abuse notation 
by writing $\nu\bigl(\{\tilde{\gamma}\}\bigr)$ for the discrete
weight $d\nu\bigl(\{\gamma\}\bigr)$ ($\tilde{\gamma}\in S$).

%%%%%%%%%%%%%%%%%%%%%%%%%%%%%%%%%%%%%%%%%%%%%%%%%%%%%%%%%%%%%%%%%%%%%%
\begin{prop}\label{respositief}
Let $(a,b,c)\in V_z^{gen}$.

{\bf (i)}
For $\mu_1<\mu_2<-(1+a)^2$ or $-(1-a)^2<\mu_1<\mu_2$ such that
$\mu(S)\cap \bigl(\mu_1,\mu_2\bigr)=\emptyset$, 
we have $E\bigl((\mu_1,\mu_2)\bigr)=0$.

{\bf (ii)} For $\tilde{\gamma}\in S$ and $f,g\in {\mathcal{H}}$ we have
\[\langle E\bigl(\{\mu(\tilde{\gamma})\}\bigr)f,g\rangle
=\bigl({\mathcal{F}}f\bigr)(\tilde{\gamma})
{\overline{\bigl({\mathcal{F}}g\bigr)(\tilde{\gamma})}}
d\nu\bigl(\{\tilde{\gamma}\}\bigr).
\]
\end{prop}
%%%%%%%%%%%%%%%%%%%%%%%%%%%%%%%%%%%%%%%%%%%%%%%%%%%%%%%%%%%%%%%%%%%%%%
\begin{proof}
Throughout the proof, we use the notations introduced in the
proof of Proposition \ref{contE}.

{\bf (i)} Let $f,g\in {\mathcal{D}}_{fin}$ and fix $\mu_1<\mu_2$ satisfying
the properties as stated in {\bf (i)}. It
suffices to prove that $\langle E\bigl((\mu_1,\mu_2)\bigr)f,g
\rangle=0$. 
Observe that \eqref{crucfor} is still valid in the present setting, 
but that the analogue of \eqref{gammalim} is now given by
\begin{equation}\label{gammalimdis}
\lim_{\epsilon \downarrow 0}\gamma[\mu\pm i\epsilon]=\gamma,
\end{equation}
where $\mu_1<\mu<\mu_2$ and where $\gamma\in {\mathbb{D}}$ is the unique
element satisfying $\mu=\mu(\gamma)$. 
By the condition $\mu(S)\cap \bigl(\mu_1,\mu_2\bigr)=\emptyset$, we have
for $\gamma\in {\mathbb{D}}$ satisfying $\mu_1<\mu(\gamma)<\mu_2$ that
\begin{equation}\label{ha}
\begin{split}
\lim_{\epsilon\downarrow 0}\langle 
\left(L-\bigl(\mu(\gamma)\pm i\epsilon\bigr)\right)^{-1}f,g\rangle=
\underset{\stackrel{(x,y)\in I\times I}{x\leq y}}{\iint}&
\frac{\phi_{\gamma}(x)\Phi_{\gamma}(y)}{W(\gamma)}
\bigl(1-\frac{1}{2}\delta_{x,y}\bigr)\\
&.\bigl(f(x){\overline{g(y)}}+f(y){\overline{g(x)}}
\bigr)\frac{d_qx}{p(x)}
\frac{d_qy}{p(y)}.
\end{split}
\end{equation}
By the bounded convergence theorem we may interchange the limit
$\epsilon\downarrow 0$ and the integration in \eqref{inversion}. Combined
with \eqref{ha}, this gives  
$\langle E\bigl((\mu_1,\mu_2)\bigr)f,g\rangle=0$, as desired.

{\bf (ii)} Let $\tilde{\gamma}\in S$ and $f,g\in {\mathcal{D}}_{fin}$.
Choose arbitrary $\mu_1<\mu_2$ such that 
$\mu(S)\cap \bigl(\mu_1,\mu_2\bigr)=\mu(\tilde{\gamma})$ and such that
$[\mu_1,\mu_2]\cap [-(1+a)^2,-(1-a)^2]=\emptyset$. Then by the first part of
the proposition, we have $E\big(\{\mu(\tilde{\gamma})\}\bigr)=
E\bigl((\mu_1,\mu_2)\bigr)$. 
We compute now $\langle E\bigl((\mu_1,\mu_2)\bigr)f,g\rangle$ 
using \eqref{inversion}.
We substitute \eqref{crucfor} in \eqref{inversion} and change the integration
parameter $\mu$ in \eqref{inversion} to the $\gamma$-parameter using 
\eqref{eigenfunction}.
Observe that for $\mu_1<\mu<\mu_2$ and $\epsilon>0$, 
we have that $\gamma[\mu+i\epsilon]$ (respectively $\gamma[\mu-i\epsilon]$)
lies in the lower (respectively upper) half plane
of $\CC$. Since $W(\gamma)^{-1}$ has a simple pole in $\tilde{\gamma}$,
it follows by Cauchy's Theorem that
\begin{equation*}
\begin{split}
\langle E\bigl((\mu_1,\mu_2)\bigr)f,g\rangle
=&\underset{\stackrel{(x,y)\in I\times I}{x\leq y}}{\iint}
\bigl(f(x){\overline{g(y)}}+f(y){\overline{g(x)}}\bigr)
\bigl(1-\frac{1}{2}\delta_{x,y}\bigr)\\
&.\phi_{\tilde{\gamma}}(x)\Phi_{\tilde{\gamma}}(y)
\left(\frac{a}{\tilde{\gamma}}\bigl(\tilde{\gamma}^{-1}-\tilde{\gamma}\bigr)
\underset{\gamma=\tilde{\gamma}}{\hbox{Res}}\bigl(W(\gamma)^{-1}\bigr)
\right)\frac{d_qx}{p(x)}\frac{d_qy}{p(y)},
\end{split}
\end{equation*}
where the factor 
$\frac{a}{\tilde{\gamma}}(\tilde{\gamma}^{-1}-\tilde{\gamma})$
arises from changing the integration variable in \eqref{inversion}
to $\gamma$ using the map $\gamma\mapsto \mu(\gamma)$,
and from the fact that one has to change sign in order to get a positive
oriented curve around $\tilde{\gamma}$. 
The proof is now completed using the explicit expression of 
$W(\gamma)$, using Remark \ref{phiateverywhere} and by symmetrizing the double
$q$-Jackson integral.
\end{proof}
%%%%%%%%%%%%%%%%%%%%%%%%%%%%%%%%%%%%%%%%%%%%%%%%%%%%%%%%%%%%%%%%%%%%

%%%%%%%%%%%%%%%%%%%%%%%%%%%%%%%%%%%%%%%%%%%%%%%%%%%%%%%%%%%%%%%%%%%%%%
\begin{cor}
Let $(a,b,c)\in V_z^{gen}$. 
The spectrum $\sigma(L)$ of the self-adjoint operator 
$\bigl(L,{\mathcal{D}}\bigr)$ is given by 
\[
\sigma(L)=[-(1+a)^2,-(1-a)^2]\cup\mu(S),\]
where $\sigma_c(L)=[-(1+a)^2,-(1-a)^2]$ is the continuous spectrum
and $\sigma_p(L)=\mu(S)$ is the point spectrum.
\end{cor}
%%%%%%%%%%%%%%%%%%%%%%%%%%%%%%%%%%%%%%%%%%%%%%%%%%%%%%%%%%%%%%%%%%%%%%
\begin{proof}
It follows from \cite[Theorem 13.27]{R},
Proposition \ref{contE}, Proposition \ref{mainprop},
Lemma \ref{Ltwo} and Proposition \ref{respositief} that 
$\sigma(L)=
[-(1+a)^2,-(1-a)^2]\cup\mu(S)$. Observe now that
$[-(1+a)^2,-(1-a)^2]\subset \sigma_c(L)$ by Corollary 
\ref{contspectrum} and $\mu(S)\subset \sigma_p(L)$ by 
Lemma \ref{Ltwo}.
It follows that $\sigma_c(L)=[-(1+a)^2,-(1-a)^2]$ and $\sigma_p(L)=\mu(S)$,
since $\sigma(L)$ is the disjoint union of $\sigma_c(L)$ and 
$\sigma_p(L)$, see \cite[Theorem 13.27]{R}. 
\end{proof}
%%%%%%%%%%%%%%%%%%%%%%%%%%%%%%%%%%%%%%%%%%%%%%%%%%%%%%%%%%%%%

%%%%%%%%%%%%%%%%%%%%%%%%%%%%%%%%%%%%%%%%%%%%%%%%%%%%%%%%%%%%%%%%%%%%%
\begin{cor}\label{posPlanch}
The functions $\phi_{\tilde{\gamma}} \in
{\mathcal{D}}\subset {\mathcal{H}}$ \textup{(}$\tilde{\gamma}\in S$\textup{)} 
are mutually orthogonal in ${\mathcal{H}}$. Their quadratic norms are given by
$\|\phi_{\tilde{\gamma}}\|^2=\nu\bigl(\{\tilde{\gamma}\}\bigr)^{-1}$.
\end{cor}
%%%%%%%%%%%%%%%%%%%%%%%%%%%%%%%%%%%%%%%%%%%%%%%%%%%%%%%%%%%%%%%%%%%%%
\begin{proof}
Orthogonality is clear since the functions $\phi_{\tilde{\gamma}}$
($\tilde{\gamma}\in S$) are eigenfunctions of the self-adjoint operator
$(L,{\mathcal{D}})$ with mutually different eigenvalues $\mu(\tilde{\gamma})$
($\tilde{\gamma}\in S$), see Lemma \ref{Ltwo}. 

It remains to derive the explicit expression for the 
quadratic norm $\|\phi_{\tilde{\gamma}}\|^2$. 
We first assume that $(a,b,c)\in V_z^{gen}$.
Observe that $E\bigl(\{\mu(\tilde{\gamma})\}\bigr)\phi_{\tilde{\gamma}}=
\nu\bigl(\{\tilde{\gamma}\}\bigr)\|\phi_{\tilde{\gamma}}\|^2
\phi_{\tilde{\gamma}}$ for $\tilde{\gamma}\in S$ by Proposition 
\ref{respositief}{\bf (ii)}. Since  
$\nu\bigl(\{\tilde{\gamma}\}\bigr)\|\phi_{\tilde{\gamma}}\|^2\not=0$ and
$E\bigl(\{\mu(\tilde{\gamma})\}\bigr)$ is a projection, it follows that 
$\|\phi_{\tilde{\gamma}}\|^2=\nu\bigl(\{\tilde{\gamma}\}\bigr)^{-1}$
for $\tilde{\gamma}\in S$, as desired.
The result now follows for $(a,b,c)\in V$ by continuity.
\end{proof}
%%%%%%%%%%%%%%%%%%%%%%%%%%%%%%%%%%%%%%%%%%%%%%%%%%%%%%%%%%%%%%%%%%%%%

%%%%%%%%%%%%%%%%%%%%%%%%%%%%%%%%%%%%%%%%%%%%%%%%%%%%%%%%%%%%%%%%%%%%%%
%%                                                                  %%
%%    The Plancherel formula and the inversion formula and the dual %%
%%      orthogonality relations                                     %%
%%                                                                  %%
%%%%%%%%%%%%%%%%%%%%%%%%%%%%%%%%%%%%%%%%%%%%%%%%%%%%%%%%%%%%%%%%%%%%%% 

\section{The Plancherel formula, inversion formula and the dual 
orthogonality relations}\label{PID}

The explicit knowledge of the resolution of the identity
$E$ of the self-adjoint operator $(L,{\mathcal{D}})$ on ${\mathcal{H}}$ 
leads directly to the Plancherel formula and the inversion formula
for the big $q$-Jacobi function transform, which we formulate in this section 
explicitly. We show that
the dual orthogonality relations imply
orthogonality relations for the continuous dual $q^{-1}$-Hahn polynomials
with respect to a one-parameter family of non-extremal weight functions.
Furthermore, the dual orthogonality relations give 
explicit sets of 
functions which complete the continuous dual $q^{-1}$-Hahn polynomials
to orthogonal bases of the corresponding $L^2$-space.

We define the
measure $d\nu(\cdot)=d\nu(\cdot;a,b,c;z)$ on $\CC^*$ by
\begin{equation}
\int_{\CC^*}f(\gamma)d\nu(\gamma)=\int_{\TT}f(\gamma)d\nu(\gamma)
+\sum_{\tilde{\gamma}\in S}f(\tilde{\gamma})d\nu(\{\tilde{\gamma}\}),
\end{equation}
where the measure $d\nu(\gamma)$ on $\TT$ is defined by \eqref{mu},
$S\subset {\mathbb{D}}$ is the discrete set defined by \eqref{Sd}
and the point mass $d\nu(\{\tilde{\gamma}\})$ for $\tilde{\gamma}\in S$
is defined by \eqref{normposdisc} for $\tilde{\gamma}>0$ and
\eqref{normstrange} for $\tilde{\gamma}<0$.

The measure $d\nu(\cdot)$ on $\CC^*$ 
is well defined for $(a,b,c)\in V$.  
Indeed, for the absolutely continuous part of the measure
the conditions on the parameters are such that the possible
zeros of the denominator of the corresponding weight function $1/|c(\cdot)|^2$
are compensated by zeros of the numerator.
For the discrete part of the measure it follows from the explicit
expressions \eqref{normposdisc} and \eqref{normstrange}
that $d\nu(\{\tilde{\gamma}\})$ is well defined and strictly positive
for all $\tilde{\gamma}\in S$. 

Let $L^2_{W}\bigl(\CC^*,d\nu\bigr)$ be the Hilbert space
of $W$-invariant $L^2$-functions with respect to
the measure $d\nu$. 
We define the big $q$-Jacobi function
transform for functions $f\in {\mathcal{H}}$
such that $f(zq^{-k})=0$ for $k\gg 0$ by
\begin{equation}\label{FourierB}
g(\gamma):=\bigl({\mathcal{F}}f\bigr)(\gamma)=\langle f,\phi_{\gamma}\rangle=
\int_{-1}^{\infty(z)}f(x){\overline{\phi_{\gamma}(x)}}\frac{d_qx}{p(x)},
\qquad \gamma\in\CC^*,
\end{equation}
and we define for functions $g\in L^2_{W}\bigl(\CC^*,d\nu\bigr)$
satisfying $g(-q^{k}abcz)=0$ for $k\gg 0$,
\begin{equation}\label{FourierinversionB}
f(x):=\bigl({\mathcal{G}}g\bigr)(x)=
\int_{\CC^*}g(\gamma)\phi_{\gamma}(x)d\nu(\gamma),\qquad
x\in I.
\end{equation}
The results in section \ref{contmass} and section \ref{pointmass} 
lead to the following main theorem of this paper.
%%%%%%%%%%%%%%%%%%%%%%%%%%%%%%%%%%%%%%%%%%%%%%%%%%%%%%%%%%%%%%%%%
\begin{thm}[The big $q$-Jacobi function transform]\label{main}
Let $z>0$ and $(a,b,c)\in V$.
The maps ${\mathcal{F}}$ and ${\mathcal{G}}$ uniquely extend
to surjective isometric isomorphisms 
\begin{equation*}
{\mathcal{F}}: {\mathcal{H}}\to  
L^2_{W}\bigl(\CC^*,d\nu\bigr),\qquad
{\mathcal{G}}: L^2_{W}\bigl(\CC^*,d\nu\bigr)\to
{\mathcal{H}}.
\end{equation*}
Furthermore, ${\mathcal{G}}={\mathcal{F}}^{-1}$, hence \eqref{FourierB}
and \eqref{FourierinversionB} give the big $q$-Jacobi function transform
pair \textup{(}interpreted in the suitable $L^2$-sense\textup{)}.
\end{thm}
%%%%%%%%%%%%%%%%%%%%%%%%%%%%%%%%%%%%%%%%%%%%%%%%%%%%%%%%%%%%%%%%%
\begin{proof}
We write $\bigl( .,.\bigr)$ for the
inner product of $L^2_{W}\bigl(\CC^*,d\nu\bigr)$.

Suppose that $(a,b,c)\in V_z^{gen}$ and let $f,g\in {\mathcal{D}}_{fin}$.
Applying Proposition \ref{contE} and Proposition \ref{respositief},
we obtain $\langle f,g\rangle=\langle E(\RR)f,g\rangle=
\bigl( {\mathcal{F}}f, {\mathcal{F}}g\bigr)$. In order to extend this
result to parameters $(a,b,c)$ in $V$, we show that 
$\langle f,g\rangle$ and $\bigl({\mathcal{F}}f,{\mathcal{F}}g\bigr)$
depend continuously on $(a,b,c)\in V$. This is clear for 
$\langle f,g\rangle$, 
while for $\bigl({\mathcal{F}}f,{\mathcal{F}}g\bigr)$ it is clear except
for the term
\[
\sum_{\tilde{\gamma}\in S: \tilde{\gamma}<0}
\bigl({\mathcal{F}}f\bigr)(\tilde{\gamma}){\overline{\bigl(
{\mathcal{F}}g\bigr)(\tilde{\gamma})}}d\nu(\tilde{\gamma})
=
\sum_{k\in\ZZ}\bigl({\mathcal{F}}f\bigr)(-q^kabcz){\overline{\bigl(
{\mathcal{F}}g\bigr)(-q^kabcz)}}d\nu(-q^kabcz),
\]
where we use the convention that $d\nu(-q^kabcz)=0$ if 
$-q^kabcz\not\in S$ for the right hand side. 
The continuity of this term follows
by Lebesgue's dominated convergence theorem, 
using the asymptotics
\begin{equation}\label{asingamma}
\phi_{-q^kabcz}(x)={\mathcal{O}}(x^{-k}),\quad
d\nu\bigl(-q^kabcz\bigr)={\mathcal{O}}(z^kq^{k(k-1)/2}),\qquad
k\to\infty,
\end{equation}
where $x\in I$, which hold uniformly for $(a,b,c)$ in compacta
of $V$. To compute the asymptotics \eqref{asingamma}
for $\phi_{\gamma}(x)$ with $x\in I_+$ we used
\eqref{phiatpositive},
and for $x\in I_-$ we used 
the formula
\begin{equation}\label{pf}
\phi_{\gamma}(x)=
\frac{\bigl(q/ab,q/ac,q\gamma^2;q\bigr)_{\infty}}
{\bigl(q\gamma/a,q\gamma/b,q\gamma/c;q\bigr)_{\infty}}
\Phi_{\gamma}^{-}(x),
\end{equation}
see Lemma \ref{collaps} (observe that the right hand side
of \eqref{pf} is well defined for
$\gamma<0$ with $|\gamma|<a^{-1}$
and can be uniquely extended by continuity to $(a,b,c)\in V$).
 
It follows that 
$\langle f,g\rangle=\bigl({\mathcal{F}}f,{\mathcal{F}}g\bigr)$ 
for $f,g\in {\mathcal{D}}_{fin}$ and $(a,b,c)\in V$,
hence ${\mathcal{F}}$ uniquely extends to an isometric isomorphism 
${\mathcal{F}}: {\mathcal{H}}\to L^2_{W}\bigl(\CC^*,d\nu\bigr)$
onto its image for all $(a,b,c)\in V$.

Let $f$ and $g$ be $W$-invariant continuous functions on $\CC^*$
with compact support. We have
$\bigl({\mathcal{G}}f\bigr)(x)=\bigl({\mathcal{G}}_cf\bigr)(x)+
\bigl({\mathcal{G}}_df\bigr)(x)$ for $x\in I$ with ${\mathcal{G}}_c$ given by
\eqref{Finvers} and
\[
\bigl({\mathcal{G}}_df\bigr)(x)=\sum_{\tilde{\gamma}\in S}
f(\tilde{\gamma})\phi_{\tilde{\gamma}}(x)d\nu\bigl(\{\tilde{\gamma}\}\bigr),
\]
which is a finite sum by the assumptions on $f$.
Assume that $(a,b,c)\in V_z^{gen}$, then it
follows from Proposition \ref{mainprop}{\bf (ii)} that 
${\mathcal{G}}_cf\in {\mathcal{R}}(P_c)$. Furthermore, it follows
from the proof of Corollary \ref{posPlanch} 
that ${\mathcal{G}}_df\in 
{\mathcal{R}}(P_d)=\bigl({\mathcal{R}}(P_c)\bigr)^{\bot}$, where
$P_d=E\bigl(\mu(S)\bigr)$. 
Hence
${\mathcal{G}}f\in {\mathcal{H}}$ and 
$\langle {\mathcal{G}}f, {\mathcal{G}}g\rangle=
\langle {\mathcal{G}}_cf, {\mathcal{G}}_cg\rangle
+\langle {\mathcal{G}}_df,{\mathcal{G}}_dg\rangle$.
It follows from Proposition \ref{mainprop} and Corollary \ref{posPlanch} 
that
\begin{equation}\label{partPlanch}
\langle {\mathcal{G}}_cf, {\mathcal{G}}_cg\rangle=
\int_{\gamma\in \TT}f(\gamma){\overline{g(\gamma)}}d\nu(\gamma),
\qquad 
\langle {\mathcal{G}}_df, {\mathcal{G}}_dg\rangle=
\sum_{\tilde{\gamma}\in S}f(\tilde{\gamma}){\overline{g(\tilde{\gamma})}}
d\nu\bigl(\{\tilde{\gamma}\}\bigr).
\end{equation}
This implies that 
$\langle {\mathcal{G}}f, {\mathcal{G}}g\rangle=\bigl(f,g\bigr)$.

It follows from
Proposition \ref{mainprop}{\bf (i)} and Corollary \ref{posPlanch}
that \eqref{partPlanch} is valid for all $(a,b,c)\in V$.
Furthermore, by continuity arguments, we have 
$\langle {\mathcal{G}}_cf,{\mathcal{G}}_dg\rangle=0$ 
for all $(a,b,c)\in V$. Hence,
${\mathcal{G}}$ uniquely extends to an isometric isomorphism
${\mathcal{G}}: L_{W}^2\bigl(\CC^*,d\nu\bigr)\to {\mathcal{H}}$
onto its image for all $(a,b,c)\in V$.

A direct computation now shows
that 
\[\bigl( {\mathcal{F}}f,g\bigr)=\langle f,{\mathcal{G}}g\rangle,
\qquad \forall f\in {\mathcal{H}}, \,\, \forall g\in 
L_{W}^2\bigl(\CC^*,d\nu\bigr)
\]
for $(a,b,c)\in V$, cf. the proof of Proposition \ref{mainprop}.
This implies that 
${\mathcal{R}}\bigl({\mathcal{F}}\bigr)=L_{W}^2\bigl(\CC^*,d\nu\bigr)$,
${\mathcal{R}}\bigl({\mathcal{G}}\bigr)={\mathcal{H}}$ and
${\mathcal{G}}={\mathcal{F}}^{-1}$, which completes the proof of the theorem.

\end{proof}
%%%%%%%%%%%%%%%%%%%%%%%%%%%%%%%%%%%%%%%%%%%%%%%%%%%%%%%%%%%%%%

In the remainder of this section we derive
a one-parameter family
of non-extremal orthogonality measures for the continuous dual 
$q^{-1}$-Hahn polynomials, and explicit sets of functions
which complement the polynomials to orthogonal bases of the 
corresponding Hilbert spaces.

We write $t_0=1/a$, $t_1=1/b$, $t_2=1/c$. The condition $(a,b,c)\in V$ 
is then equivalent to the conditions
\begin{equation}\label{con5}
t_i>0,\qquad t_it_j>1 \,\, (i\not=j)
\end{equation}
on the parameters $(t_0,t_1,t_2)$.
We define polynomials 
$p_k(\gamma)=p_k(\gamma;t_0,t_1,t_2;q^{-1})$ ($k\in\ZZ_+$) in 
$\gamma+\gamma^{-1}$  by
\begin{equation}\label{contdual}
p_k(\gamma)=\phi_{\gamma}(-q^k;t_0^{-1},t_1^{-1},t_2^{-1})=
{}_3\phi_2\left( \begin{array}{c}
             q^k, t_0\gamma, t_0/\gamma \\ 
             t_0t_1, t_0t_2
            \end{array} ; q^{-1},q^{-1} \right),\qquad
            k\in\ZZ_+.
\end{equation}
Here we use \cite[Exercise 1.4(i)]{GR} to obtain the second equality.
Observe that the second equality of \eqref{contdual} shows
that $\{p_k(\cdot;t_0,t_1,t_2;q^{-1}) \}_{k\in\ZZ_+}$ are exactly the
continuous dual $q^{-1}$-Hahn polynomials, i.e. Askey-Wilson polynomials
in base $q^{-1}$ with one of the four parameters equal to zero.  
For $z>0$, we define a measure 
$d\sigma_z(.)=d\sigma_z(.;t_0,t_1,t_2;q^{-1})$ on $\CC^*$ by
\begin{equation}
\int_{\CC^*}f(\gamma)d\sigma_z(\gamma)=
\frac{1}{M}\int_{\CC^*}f(\gamma)d\nu(\gamma;t_0^{-1},t_1^{-1},t_2^{-1};z).
\end{equation}
Explicitly, we have
\begin{equation}
\int_{\CC^*}f(\gamma)d\sigma_z(\gamma)=
\frac{1}{4\pi i}\int_{\TT}f(\gamma)w_z(\gamma)\frac{d\gamma}{\gamma}
+\sum_{\tilde{\gamma}\in S_z}f(\tilde{\gamma})
\underset{\gamma=\tilde{\gamma}}{\hbox{Res}}\left(\frac{w_z(\gamma)}{\gamma}
\right),
\end{equation}
where
\begin{equation*}
\begin{split}
S_z=S^{-1}=&\left\{\frac{q^k}{e} \, | \, e\in \{t_0,t_1,t_2\},  k\in \ZZ_+: 
\frac{q^k}{e}>1\right\}\\
&\qquad\quad\cup
\left\{ \frac{-q^{k}t_0t_1t_2}{z} \, | \, k\in \ZZ: 
\frac{q^{k}t_0t_1t_2}{z}>1 \right\}
\end{split}
\end{equation*}
and with weight function $w_z(.)=w_z(.;t_0,t_1,t_2;q^{-1})$ given by
\begin{equation}
w_z(\gamma)=\frac{\bigl(\gamma^2,1/\gamma^{2};q\bigr)_{\infty}}
{\theta\bigl(-z\gamma/t_0t_1t_2,
-z/t_0t_1t_2\gamma\bigr)\prod_{j=0}^2
\bigl(\gamma/t_j, 1/t_j\gamma;q\bigr)_{\infty}}
\end{equation}
for parameters $(t_0,t_1,t_2)$ such that the poles of $w_z(\gamma)$ at
$\gamma\in S_z$ are simple. 
Finally, we define $W$-invariant
functions $r_k^z(\gamma)=r_k^z(\gamma;t_0,t_1,t_2;q^{-1})$ 
for $k\in\ZZ$ by
\begin{equation}\label{contdualextra}
\begin{split}
r_k^z(\gamma)&=\phi_{\gamma}(zq^k;t_0^{-1},t_1^{-1},t_2^{-1})\\ 
&=\frac{\bigl(\gamma/t_0, 1/t_1t_2, 
-q^kz/t_0t_1t_2\gamma;q\bigr)_{\infty}}
{\bigl(1/t_0t_1, 1/t_0t_2, -q^kz/t_1t_2;q\bigr)_{\infty}}
{}_3\phi_2\left( \begin{array}{c}
             1/t_1\gamma, 1/t_2\gamma, -q^kz/t_1t_2 \\ 
             1/t_1t_2, -q^kz/t_0t_1t_2\gamma
            \end{array} ; q,\gamma/t_0 \right)
\end{split}
\end{equation}
where the second equality holds for 
$\gamma\in \CC^*$ with $|\gamma/t_0|<1$, see \eqref{ancontB}.

%%%%%%%%%%%%%%%%%%%%%%%%%%%%%%%%%%%%%%%%%%%%%%%%%%%%%%%%%%%%%%%%%%%%%%
\begin{thm}\label{Hahntheorem}
Let $z>0$ and fix parameters $t_i$ satisfying the conditions 
\eqref{con5}. Then,
$\{p_k\}_{k\in\ZZ_+}\cup \{r_k^z\}_{k\in\ZZ}$ is an orthogonal
basis of the Hilbert space $L_W^2\bigl(\CC^*,d\sigma_z\bigr)$.
The quadratic norms of the basis elements are given by
\begin{equation*}
\begin{split}
\int_{\CC^*}|p_k(\gamma)|^2d\sigma_z(\gamma)&=\frac{\theta\bigl(-z\bigr)}
{\theta\bigl(-z/t_0t_1,-z/t_0t_2, -z/t_1t_2\bigr)}\\
&\qquad.  \frac{1}{\bigl(q,1/t_0t_1,1/t_0t_2,1/t_1t_2;q\bigr)_{\infty}}
\frac{\bigl(q,1/t_1t_2;q\bigr)_k}{\bigl(1/t_0t_1,1/t_0t_2;q\bigr)_k}q^{-k},\\
\int_{\CC^*}|r_k^z(\gamma)|^2d\sigma_z(\gamma)&=
\frac{\theta\bigl(-1/z\bigr)}
{\theta\bigl(-z/t_0t_1,-z/t_0t_2,-z/t_1t_2\bigr)}\\
&\qquad. \frac{1}{\bigl(1/t_0t_1,1/t_0t_2;q\bigr)_{\infty}^2}
\frac{\bigl(-q^kz/t_0t_1,-q^kz/t_0t_2;q\bigr)_{\infty}}
{\bigl(-q^kz/t_1t_2,-q^{1+k}z;q\bigr)_{\infty}}q^{-k}.
\end{split}
\end{equation*}
\end{thm}
%%%%%%%%%%%%%%%%%%%%%%%%%%%%%%%%%%%%%%%%%%%%%%%%%%%%%%%%%%%%%%%%%%
\begin{proof}
It follows from Theorem \ref{main} that 
\[
\{ \gamma\mapsto \phi_{\gamma}(-q^k) \, | \, k\in\ZZ_+\}\cup
\{\gamma\mapsto \phi_{\gamma}(zq^k)\, | \, k\in\ZZ\}
\]
is
an orthogonal basis of $L_W^2\bigl(\CC^*,d\nu\bigr)$, and that their quadratic
norms are given by
\begin{equation*}
\begin{split}
\int_{\CC^*}|\phi_{\gamma}(-q^k)|^2d\nu(\gamma)&=\frac{p(-q^k)}{(1-q)q^k},
\qquad k\in\ZZ_+,\\
\int_{\CC^*}|\phi_{\gamma}(zq^k)|^2d\nu(\gamma)&=\frac{p(zq^k)}{(1-q)zq^k},
\qquad k\in\ZZ.
\end{split}
\end{equation*}
The theorem follows now immediately by setting 
$t_0=1/a$, $t_1=1/b$ and $t_2=1/c$ and using the explicit
expressions \eqref{pr} and \eqref{M} of the function $p(\cdot)$ and the
constant $M$.
\end{proof}
%%%%%%%%%%%%%%%%%%%%%%%%%%%%%%%%%%%%%%%%%%%%%%%%%%%%%%%%%%%%%%%%%%%%%%%

%%%%%%%%%%%%%%%%%%%%%%%%%%%%%%%%%%%%%%%%%%%%%%%%%%%%%%%%%%%%%%%%%%%%%%%
\begin{rem} 
We remarked in section \ref{section1} that
the big $q$-Jacobi function transform is associated with harmonic analysis
on the $SU(1,1)$ quantum group. Analogous considerations for the quantum
group of plane motions lead to the so-called big $q$-Hankel transform,
see \cite{Koelink2} and \cite{CKK}. The corresponding function theoretic
aspects of the big $q$-Hankel transform are discussed in detail
in \cite{CKK}. 

The dual orthogonality relations for the big $q$-Hankel transform
have a similar interpretation as the dual orthogonality relations for
the big $q$-Jacobi function transform, namely, they give
orthogonality relations for Moak's  
$q$-Laguerre polynomials with respect to a one-parameter family
of non-extremal orthogonality measures, as well as explicit sets of
functions which complement the $q$-Laguerre polynomials
to orthogonal bases of the associated Hilbert spaces, 
see \cite[Theorem 4.1]{CKK}.
\end{rem}
%%%%%%%%%%%%%%%%%%%%%%%%%%%%%%%%%%%%%%%%%%%%%%%%%%%%%%%%%%%%%%%%%%%%%

%%%%%%%%%%%%%%%%%%%%%%%%%%%%%%%%%%%%%%%%%%%%%%%%%%%%%%%%%%%%%%%%%%
%%                                                              %%
%%  Big $q$-Jacobi polynomials: functional analytic approach    %%
%%                                                              %%
%%%%%%%%%%%%%%%%%%%%%%%%%%%%%%%%%%%%%%%%%%%%%%%%%%%%%%%%%%%%%%%%%%

\section{Big $q$-Jacobi polynomials: functional analytic approach}
\label{polsection}

We show in this section how the orthogonality relations and the
quadratic norm evaluations for the big $q$-Jacobi polynomials can
be derived from a functional analytic approach. The arguments are
closely related to the ones used for the big $q$-Jacobi function transform, so
we merely sketch the steps and indicate the main differences.
To avoid confusion with previous notations, we label definitions
in this section with a subscript (or superscript) $\wp$ 
(indicating that it is connected with the polynomial case).

The conditions on the parameters $a,b,c$ are now taken to be
\begin{equation}\label{condparB}
ab,qa/b,ac,qa/c<1,\qquad  bc<0.
\end{equation}

We consider the second order $q$-difference equation $L$ 
acting on the space ${\mathcal{F}}(I_{\wp})$ of complex-valued functions
$f:I_p\to\CC$, where $I_{\wp}=[-1,-q/bc]_q=\{-q^k\}_{k\in\ZZ_+}\cup
\{-q^{k+1}/bc\}_{k\in\ZZ_+}$. At the end-points $x=-1$ and $x=-q/bc$
this should be read as $\bigl(Lf\bigr)(x)=A(x)(f(qx)-f(x))$.
We can write this in self-adjoint form, similarly
as was done in Lemma \ref{self-adjointform} for the non-compact case,
with the same functions $p$ and $r$, see \eqref{pr}.
Observe that $p(x)>0$ for all $x\in I_{\wp}$ by the conditions
\eqref{condparB} on the parameters.

We define ${\mathcal{H}}_{\wp}=\{f:I_{\wp}\to\CC \,\, | \,\,
\|f\|_{\wp}^2=\langle f,f\rangle_{\wp}<\infty\}$, where 
\[
\langle f,g\rangle_{\wp}=\int_{-1}^{-q/bc}f(x){\overline{g(x)}}
\frac{d_qx}{p(x)}.
\]
It is well known that the big $q$-Jacobi polynomials,
which are explicitly given by \eqref{connpol}, form an orthogonal basis
of ${\mathcal{H}}_{\wp}$, see \cite{AA}. This fact and the evaluation of
the corresponding quadratic norms have been derived in \cite{AA}
using the $q$-binomial formula \cite[(1.3.2)]{GR} and the 
$q$-Pfaff-Saalsch{\"u}tz formula \cite[(1.7.2)]{GR}. 
In this section we derive these results by functional analytic methods.

Similarly as in the non-compact setting, we truncate the inner product
by
\[\langle f,g\rangle_{\wp,k,l}=
\left(\int_{-1}^{-q^{k+1}}+\int_{-q^{l+2}/bc}^{-q/bc}\right)
f(x){\overline{g(x)}}\frac{d_qx}{p(x)}
\] 
for $k,l\in\ZZ_+$. The analogue of Lemma \ref{stoktermen} is then
given by 
\[\langle Lf,g\rangle_{\wp,k,l}-\langle f,Lg\rangle_{\wp,k,l}=
W(f,{\overline{g}})(-q^k)-W(f,{\overline{g}})(-q^{l+1}/bc),
\]
with the Wronskian as defined in \eqref{Wronskian}.

Let $\alpha\in\TT$, then we write  ${\mathcal{D}}_{\wp,\alpha}$ for the
functions $f\in {\mathcal{H}}_{\wp}$ such that $Lf\in
{\mathcal{H}}_{\wp}$ and $f(0^+)=\alpha f(0^-)$, $f'(0^+)=\alpha
f'(0^-)$, where now $f(0^-)=\lim_{k\to\infty}f(-q^k)$ and
$f(0^+)=\lim_{k\to\infty}f(-q^k/bc)$, etc.
(cf. section \ref{sdomain}).
By similar arguments as in section \ref{sdomain}, we have
%%%%%%%%%%%%%%%%%%%%%%%%%%%%%%%%%%%%%%%%%%%%%%%%%%%%%%%%%%%%%%%%
\begin{prop}
The operator $(L,{\mathcal{D}}_{\wp,\alpha})$ on ${\mathcal{H}}_{\wp}$
is self-adjoint for all $\alpha\in \TT$. 
\end{prop}
%%%%%%%%%%%%%%%%%%%%%%%%%%%%%%%%%%%%%%%%%%%%%%%%%%%%%%%%%%%%%%%%

We denote
${\mathcal{D}}_{\wp}={\mathcal{D}}_{\wp,1}$.
Now observe that $\phi_{\gamma}, \psi_{\gamma}\in {\mathcal{D}}_{\wp}$
for $\gamma\in\CC^*$ by Lemma \ref{oknul}
(see \eqref{phiB} and \eqref{psiB} for the definition of $\phi_{\gamma}$  
and $\psi_{\gamma}$, respectively). 
Furthermore, from the arguments as given in section \ref{eigenf}
it follows that 
$\bigl(L\phi_{\gamma}\bigr)(x)=\mu(\gamma)\phi_{\gamma}(x)$ for
$x\in [-1,-q/bc)_q=I_{\wp}\setminus \{-q/bc\}$
and $\bigl(L\psi_{\gamma}\bigr)(x)=
\mu(\gamma)\psi_{\gamma}(x)$ for $x\in
(-1,-q/bc]_q=I_{\wp}\setminus \{-1\}$,
where $\mu(\gamma)$ is
given by \eqref{eigenfunction}. The Wronskian 
$W_{\wp}(\gamma)=W(\psi_{\gamma},\phi_{\gamma})\in {\mathcal{F}}(I_{\wp})$ 
can again be seen to be constant on $I_{\wp}$ (cf.
Lemma \ref{generaleigenfunction}{\bf (iii)}), and
an explicit expression of the Wronskian 
$W_{\wp}(\gamma)$ is given by 
\[
W_{\wp}(\gamma)=(1-q)\frac{\bigl(a\gamma,a/\gamma;q\bigr)_{\infty}\theta(bc)}
{\bigl(ab,ac,qa/b,qa/c;q\bigr)_{\infty}}.
\]
Indeed, observe that Proposition \ref{Wpsiphi}
is also valid for the present choice \eqref{condparB} of parameter values
by analytic continuation. 

In particular, the functions
$\phi_{\gamma}$ and $\psi_{\gamma}$ in ${\mathcal{D}}_{\wp}$
are linearly independent if and only if $\gamma\not\in S_{pol}$, where
$S_{pol}=\{\gamma_n^{\pm 1}\}_{n\in\ZZ_+}$, $\gamma_n=aq^n$. Hence
$\bigl(L\phi_{\gamma}\bigr)(-q/bc)\not=\mu(\gamma)\phi_{\gamma}(-q/bc)$
and 
$\bigl(L\psi_{\gamma}\bigr)(-1)\not=\mu(\gamma)\psi_{\gamma}(-1)$
if $\gamma\not\in S_{pol}$, cf. Corollary \ref{psifoutind}{\bf (iii)}.

Let $\mu=\mu(\gamma)\in \CC\setminus\RR$. 
We define the Green kernel $K_{\gamma}^{\wp}(x,y)$
for $x,y\in I_{\wp}$ by
\begin{equation*}
K_{\gamma}^{\wp}(x,y)=
\begin{cases}
W_{\wp}(\gamma)^{-1}\psi_{\gamma}(x)\phi_{\gamma}(y), 
\qquad &y\leq x,\\
W_{\wp}(\gamma)^{-1}\phi_{\gamma}(x)\psi_{\gamma}(y),
\qquad &y\geq x.
\end{cases}
\end{equation*}
We have a well defined linear map ${\mathcal{H}}_{\wp}\to 
{\mathcal{F}}(I_{\wp})$ which maps $f\in {\mathcal{H}}_{\wp}$ to 
\begin{equation}
G_{f}^{\wp}(x,\gamma)=\langle f, 
{\overline{K_{\gamma}^{\wp}(x,\cdot)}}\rangle_{\wp}, 
\qquad x\in I_{\wp}.
\end{equation}
By similar arguments as in the proof of Proposition \ref{G},
we derive that for $f\in {\mathcal{H}}_{\wp}$ and for
$\gamma\in\CC^*$ such that $\mu(\gamma)\in \CC\setminus\RR$,
\begin{equation}\label{inverpol}
G_f^{\wp}(\cdot,\gamma)=\bigl(L-\mu(\gamma).\hbox{Id}\bigr)^{-1}f.
\end{equation}
For the proof of \eqref{inverpol}, we have to consider
case 2 of the proof of Proposition
\ref{G} twice, namely for
the end-point $-1$ as well as for the end-point $-q/bc$.
The arguments go through for both end-points, since $\phi_{\gamma}$
(respectively $\psi_{\gamma}$)
satisfies the eigenvalue equation $(Lf)(x)=\mu(\gamma)f(x)$ in the
end-point $x=-1$ (respectively $x=-q/bc$).

We are now in a position to compute the resolution of the identity
$E_{\wp}$ for the self-adjoint operator $\bigl(L,{\mathcal{D}}_{\wp}\bigr)$
on ${\mathcal{H}}_{\wp}$ using \eqref{inversion}.
For the moment it is convenient to assume that the parameters also satisfy 
the condition $a^2\not\in \{q^k\}_{k\in \ZZ}$. This condition can be
removed later on by continuity.

We keep the notations of the proof of Proposition \ref{contE}.
Choose $f,g\in {\mathcal{H}}_{\wp}$ with finite support.
Then we have for $\mu\in\RR$ and $\epsilon>0$, that
\begin{equation}\label{polcrucfor}
\begin{split}
\langle \left(L-\bigl(\mu\pm i\epsilon\bigr)\right)^{-1}f,g\rangle_{\wp}=
&{\underset{\stackrel{(x,y)\in I_{\wp}\times I_{\wp}}{x\leq y}}{\iint}}
\frac{\phi_{\gamma[\mu\pm i\epsilon]}(x)
\psi_{\gamma[\mu\pm i\epsilon]}(y)}{W_{\wp}\bigl(\gamma[\mu\pm
i\epsilon]\bigr)}\bigl(1-\frac{1}{2}\delta_{x,y}\bigr)\\
&\qquad.\bigl(f(x){\overline{g(y)}}+f(y){\overline{g(x)}}
\bigr)\frac{d_qx}{p(x)}\frac{d_qy}{p(y)}.
\end{split}
\end{equation}
Let $\xi>0$ 
such that $\mu\bigl(S_{pol}\bigr)\cap [-(1+a)^2-\xi, -(1-a)^2+\xi]=\emptyset$.
Using \eqref{gammalim}, \eqref{gammalimdis} and 
the invariance of $\phi_{\gamma}$, $\psi_{\gamma}$ and $W_{\wp}(\gamma)$
under $\gamma\leftrightarrow\gamma^{-1}$, we obtain from \eqref{polcrucfor}
that for all $-(1+a)^2-\xi<\mu<-(1-a)^2+\xi$, 
\[
\lim_{\epsilon\downarrow 0}\left(
\langle \bigl(L-(\mu+i\epsilon)\bigr)^{-1}f,g\rangle_{\wp}
-\langle \bigl(L-(\mu-i\epsilon)\bigr)^{-1}f,g\rangle_{\wp}\right)=0.
\]
It follows now from \eqref{inversion} that
$E_{\wp}\bigl([-(1+a)^2,-(1-a)^2]\bigr)=0$.

For $-\infty<\mu_1<\mu_2<-(1+a)^2$ or 
$-(1-a)^2<\mu_1<\mu_2<\infty$
such that $(\mu_1,\mu_2)\cap \mu(S_{pol})=\emptyset$, we have
$E_{\wp}\bigl((\mu_1,\mu_2)\bigr)=0$, cf. the proof of Proposition 
\ref{respositief}{\bf (i)}. 
Setting $\mu_n=\mu(\gamma_n)$ for $n\in \ZZ_+$, 
we have, due to the simple pole of
$W_{\wp}(\gamma)^{-1}$ in $\gamma=\gamma_n$,
\begin{equation*}
\begin{split}
\langle E_{\wp}\bigl(\{\mu_n\}\bigr)f,g\rangle_{\wp}
=&\underset{\stackrel{(x,y)\in I_{\wp}\times I_{\wp}}{x\leq y}}{\iint}
\bigl(f(x){\overline{g(y)}}+f(y){\overline{g(x)}}\bigr)
\bigl(1-\frac{1}{2}\delta_{x,y}\bigr)\\
&.\phi_{\gamma_n}(x)\psi_{\gamma_n}(y)
\left(\frac{a}{\gamma_n}\bigl(\gamma_n^{-1}-\gamma_n\bigr)
\underset{\gamma=\gamma_n}{\hbox{Res}}\bigl(W_{\wp}(\gamma)^{-1}\bigr)
\right)\frac{d_qx}{p(x)}\frac{d_qy}{p(y)}
\end{split}
\end{equation*}
for $f,g\in {\mathcal{H}}_{\wp}$ with finite support, cf. the proof of
Proposition \ref{respositief}{\bf (ii)}. Using \eqref{connpol} to rewrite
$\psi_{\gamma_n}$ as a multiple of $\phi_{\gamma_n}$, and by
symmetrizing the double $q$-Jackson integral, we obtain
\begin{equation}\label{almostthere}
\langle E_{\wp}(\{\mu_n\})f,g\rangle_{\wp}={\mathcal{N}}_{\wp}(n)^{-1}
\langle f,\psi_{\gamma_n}\rangle_{\wp}\langle \psi_{\gamma_n},g\rangle_{\wp},
\end{equation}
where
\begin{equation*}
\begin{split}
{\mathcal{N}}_{\wp}(n)=&\left(\frac{a}{\gamma_n}
\bigl(\gamma_n^{-1}-\gamma_n\bigr)
\frac{\bigl(qa/b,qa/c;q\bigr)_n}{\bigl(ab,ac;q\bigr)_{n}}
\left(\frac{bc}{q}\right)^n
\underset{\gamma=\gamma_n}{\hbox{Res}}\bigl(W_{\wp}(\gamma)^{-1}\bigr)
\right)^{-1}\\
&=\frac{(1-q)\bigl(q,bc,q/bc,q^{2n+1}a^2;q\bigr)_{\infty}
\bigl(q,q^{n}a^2;q\bigr)_n}
{\bigl(qa/b,qa/c,q^{n}ab,q^{n}ac;q\bigr)_{\infty}
\bigl(qa/b,qa/c;q\bigr)_n}\left(-\frac{q^{(3-n)/2}}{bc}\right)^n.
\end{split}
\end{equation*}
Observe that \eqref{almostthere} holds for all 
$f,g\in {\mathcal{H}}_{\wp}$ by continuity.

We can now immediately recover the orthogonality relations
and quadratic norm evaluations for the big $q$-Jacobi polynomials
by applying well known properties of resolution of identities to
$E_{\wp}$. The result is as follows. 
%%%%%%%%%%%%%%%%%%%%%%%%%%%%%%%%%%%%%%%%%%%%%%%%%%%%%%%%%%%%%%%%%%%
\begin{thm}[\cite{AA}]
If the parameters $(a,b,c)$ satisfy the conditions \eqref{condparB}, then
the polynomials $\{\psi_{\gamma_n}\}_{n\in\ZZ_+}$ 
form an orthogonal
basis of ${\mathcal{H}}_{\wp}$, and their quadratic norms are given by 
\[
\|\psi_{\gamma_n}\|_{\wp}^2={\mathcal{N}}_{\wp}(n),\qquad n\in\ZZ_+.
\]
\end{thm}
%%%%%%%%%%%%%%%%%%%%%%%%%%%%%%%%%%%%%%%%%%%%%%%%%%%%%%%%%%%%%%%%%%

%%%%%%%%%%%%%%%%%%%%%%%%%%%%%%%%%%%%%%%%%%%%%%%%%%%%%%%%%%%%%%%%
%%                                                            %%
%%                       Bibliography                         %%
%%                                                            %%
%%%%%%%%%%%%%%%%%%%%%%%%%%%%%%%%%%%%%%%%%%%%%%%%%%%%%%%%%%%%%%%%

\bibliographystyle{amsplain}

\begin{thebibliography}{000}
\bibitem[1]{AA} G.E. Andrews, R. Askey, {\it Classical orthogonal
    polynomials}, in: {\it Polyn{\^o}mes orthogonaux et applications}, 
Lecture Notes in Math. {\bf 1171}, Springer, Berlin and New York,
pp. 36--62.
\bibitem[2]{Ak} N.I. Akhiezer, {\it The classical moment problem and some 
related questions in analysis}, Hafner (1965).
\bibitem[3]{AW2} R. Askey, J. Wilson, {\it Some basic 
hypergeometric orthogonal 
polynomials that generalize Jacobi polynomials}, 
Mem. Amer. Math. Soc. {\bf 54} (1985), no. 319. 
\bibitem[4]{BM} B.L.J. Braaksma, B. Meulenbeld, {\it Integral transforms
with generalized Legendre functions as kernels}, Compositio Math.
{\bf 18} (1967), pp. 235--287.
\bibitem[5]{CKK} N. Ciccoli, E. Koelink, T.H. Koornwinder, {\it
$q$-Laguerre polynomials and big $q$-Bessel functions and their orthogonality
relations}, preprint (1998), to appear in Methods Appl. Anal., 
math.CA/9805023.
\bibitem[6]{DS} N. Dunford, J.T. Schwartz, {\it Linear operators, part II: 
Spectral theory}, Interscience (1963).
\bibitem[7]{GR} G. Gasper, M. Rahman, {\it Basic 
hypergeometric series}, 
Encyclopedia of Mathematics and its Applications 35, 
Cambridge University Press (1990).
\bibitem[8]{G} F. G{\"o}tze, {\it Verallgemeinerung einer Integral
transformation ovn Mehler-Fock durch den von Kuipers und Meulenbeld
eingef{\"u}hrten Kern $P_k^{m,n}(z)$}, Indag. Math. {\bf 27} (1965),
pp. 396--404. 
\bibitem[9]{GIM} D.P. Gupta, M.E.H. Ismail, D.R. Masson,
{\it Contiguous relations, basic hypergeometric functions, and
orthogonal polynomials. III. Associated continuous dual $q$-Hahn
polynomials}, J. Comp. Appl. Math. {\bf 68} (1996), pp. 115-149.
\bibitem[10]{IM} M.E.H. Ismail, D.R. Masson, {\it $q$-Hermite polynomials,
biorthogonal rational functions, and $q$-beta integrals}, Trans. Amer.
Math. Soc. {\bf 346} (1994), pp. 63--116.
\bibitem[11]{K} T. Kakehi, {\it Eigenfunction expansion associated with
the Casimir operator on the quantum group $SU(1,1)$}, 
Duke Math. J. {\bf 80} (1998), pp. 535--573.
\bibitem[12]{KMU} T. Kakehi, T. Masuda, K. Ueno, {\it
Spectral analysis of a $q$-difference operator which arises from the
quantum $SU(1,1)$ quantum group}, J. Operator Theory {\bf 33} (1995),
pp. 159--196.
\bibitem[13]{Koelink2} H.T. Koelink, {\it The quantum group of plane motions
and basic Bessel functions}, Indag. Mathem. N.S. {\bf 6} (1995), pp. 197--211.
\bibitem[14]{Koel} H.T. Koelink, {\it Askey-Wilson polynomials and the 
quantum $SU(2)$ group: survey and applications}, Acta Appl. Math.
{\bf 44} (1996), pp. 295--352.
\bibitem[15]{KV} H.T. Koelink, J. Verding, {\it Spectral analysis and the 
Haar functional on the quantum $SU(2)$ group}, 
Comm. Math. Phys. {\bf 177} (1996), pp. 399--415.
\bibitem[16]{K1} T.H. Koornwinder, {\it Jacobi functions and analysis 
on noncompact semisimple Lie groups}, in: {\it Special functions: 
Group theoretic aspects and applications}, R.A. Askey, T.H. Koornwinder
and  W. Schempp (eds.), Reidel, 1984, pp. 1--85.
\bibitem[17]{K2} T.H. Koornwinder, {\it Askey-Wilson 
polynomials as zonal spherical 
functions on the $SU(2)$ quantum group}, SIAM J. 
Math. Anal. {\bf 24} (1993), pp.
795--813.
\bibitem[18]{MR} D.R. Masson, J. Repka, {\it Spectral theory of Jacobi
matrices in $l^2(\ZZ)$ and the $su(1,1)$ Lie algebra}, SIAM J. Math. Anal. 
{\bf 22} (1991), pp. 1131--1146.
\bibitem[19]{NM} M. Noumi, K. Mimachi, {\it Quantum $2$-spheres
and big $q$-Jacobi polynomials}, Comm. Math. Phys. {\bf 128} (1990),
pp. 521--531.
\bibitem[20]{R} W. Rudin, {\it Functional analysis},
Tata McGraw-Hill (1973).
\bibitem[21]{Simon} B. Simon, {\it The classical moment problem as a 
self-adjoint finite difference operator}, Adv. Math. {\bf 137} (1998),
pp. 82--203.
\bibitem[22]{VK} L.L. Vaksman, L.I. Korogodsky, {\it Spherical functions on the
quantum group $SU(1,1)$ and a $q$-analogue of the Mehler-Fock formula},
Funct. Anal. Appl. {\bf 25} (1991), pp. 48--49.
\bibitem[23]{VS} L.L. Vaksman, Ya. S. Soibelman, {\it Algebra of functions
on the quantum group $SU(2)$}, Funct. Anal. Appl. {\bf 22} (1988), pp. 170--181.
\end{thebibliography}

\end{document}